\documentclass{article}

\usepackage{graphicx} 
\usepackage{algorithm}
\usepackage{algorithmic}

\usepackage[utf8]{inputenc}
\usepackage{graphicx}
\usepackage{bm}
\usepackage{comment}
\usepackage{psfrag}
\usepackage{latexsym,amsmath,amsfonts,amscd,amsthm}
\usepackage{changebar}
\usepackage{color}
\usepackage{bm}
\usepackage{tikz}
\usepackage{multirow}
\usepackage{xcolor}
\usepackage{hyperref}
\usepackage{setspace}
\usepackage{amssymb}
\usepackage{mathrsfs}
\usepackage{caption}
\usepackage{subcaption}

\usepackage[title]{appendix}
\usepackage{ulem}
\usepackage{tikz}

\usetikzlibrary{arrows,backgrounds,positioning}

\newlength{\spse}
\newtheorem{thm}{Theorem}[section]

\newtheorem{rem}[thm]{Remark}

\newcommand{\BV}{\mathbf{V}}
\newcommand{\BT}{\mathbf{T}}
\newcommand{\BI}{\mathbf{I}}

\newcommand{\BA}{\mathbf{A}}

\newcommand{\BOmega}{\boldsymbol{\Omega}}

\newcommand{\bx}{\mathbf{x}}

\newcommand{\bc}{\mathbf{c}}
\newcommand{\bb}{\mathbf{b}}
\newcommand{\bp}{\mathbf{p}}
\newcommand{\bq}{\mathbf{q}}

\newcommand{\bn}{\boldsymbol{n}}

\newcommand{\half}{\frac{1}{2}}
\newcommand{\br}{\boldsymbol{r}}

\newcommand{\bphi}{\boldsymbol{\phi}}
\newcommand{\bpsi}{\boldsymbol{\psi}}

\newcommand{\BM}{{\bf{M}} }
\newcommand{\BD}{{\bf{D}}}
\newcommand{\BX}{{\bf{X}}}

\newcommand{\BB}{{\bf{B}} }
\newcommand{\BSigma}{{\mathbf{\Sigma}}}

\newcommand{\bz}{{\bf{z}}}
\newcommand{\by}{{\bf{y}}}

\newcommand{\dx}{\Delta x}
\newcommand{\dy}{\Delta y}

\newcommand{\BmcD}{\boldsymbol{\mathcal{D}}}
\newcommand{\BmcE}{\boldsymbol{\mathcal{E}}}

\newcommand{\algorithmicbreak}{\textbf{break}}
\newcommand{\BREAK}{\STATE \algorithmicbreak}

\title{An Inexact Low-Rank Source Iteration for Steady-State Radiative Transfer Equation with Diffusion Synthetic Acceleration}

\author{Wei Guo\footnote{Department of Mathematics and Statistics, Texas Tech University, Lubbock, TX, 70409, USA.},\quad   Zhichao Peng\footnote{Department of Mathematics, The Hong Kong University of Science and Technology, Clear Water Bay, Kowloon, Hong Kong, China. Corresponding author, E-mail: pengzhic@ust.hk.}}

\begin{document}

\maketitle

\begin{abstract}
    We propose an inexact low-rank source iteration with diffusion synthetic acceleration (SI–DSA) for solving the multidimensional steady-state radiative transfer equation (RTE) in the second-order formulation. The angular flux is represented in either a low-rank matrix or hierarchical Tucker tensor (HTT) format, enabling substantial reductions in computational resources. Each SI step is solved using a preconditioned low-rank conjugate gradient (CG) method with a diffusion preconditioner. To further improve efficiency, we introduce an adaptive inexact strategy that dynamically relaxes the inner CG tolerance during early SI iterations. The method exploits the tensor-product structure of the discretized operators to perform all matrix–vector operations in low-rank form. Numerical experiments on 2D2V benchmark problems—including diffusion-dominated, transport-dominated, and multiscale problems—demonstrate that the proposed approach achieves errors on the order of $10^{-4}–10^{-5}$ relative to full-rank reference solutions, while reducing the degrees of freedom by up to two orders of magnitude. In the diffusion-dominated case, the low-rank solver achieves speedups exceeding $90\times$ over its full-rank counterpart and remains competitive in solving challenging transport-dominated and multiscale problems while providing substantial storage savings. To our knowledge, this work provides the first low-rank SI–DSA framework for multidimensional steady-state RTE.
\end{abstract}

\section{Introduction\label{sec:introduction}}
Radiative transfer equation (RTE) is a kinetic equation modeling particle interactions with a background medium. It has critical applications in nuclear engineering, medical imaging, astrophysics, and remote sensing. This paper focuses on the steady-state RTE and aims to develop an efficient numerical solver to address the computational challenges inherent in RTE simulations.

Among classic iterative linear solvers for RTE, one effective and popular approach is the discrete ordinates ($S_N$) Source Iteration (SI) \cite{case1968linear, Adams2002FastIM} with Diffusion Synthetic Acceleration (DSA). In particular, SI decouples the linear system for different angular flux by freezing the scalar flux, while DSA introduces a correction to the scalar flux after each SI step by solving the diffusion limit of an error equation \cite{kopp1963synthetic,alcouffe1977dittusion,adams1992diffusion,Adams2002FastIM,larsen2009advances}. DSA can be seen as a crucial preconditioning step, particularly in scattering-dominated regimes, where the convergence of unaccelerated SI could be arbitrarily slow \cite{Adams2002FastIM,larsen2009advances}. Originally developed for the first-order RTE, SI-DSA has been extended to the second-order formulation \cite{palii2020convergent} and further accelerated by exploiting Krylov methods \cite{warsa2004krylov}.

Despite its popularity and effectiveness, SI-DSA, as with other classical RTE solvers, suffers from the curse of dimensionality, since the angular flux is posed in a high-dimensional phase space, requiring substantial memory and computational resources. To address this challenge, low-rank solvers leveraging advances in matrix and tensor decompositions have emerged as a powerful tool for efficient RTE simulations.  For time-dependent RTE, fixed-rank dynamical low-rank algorithms (DLRAs)  have been developed based on its first order formulation \cite{peng2020low,peng2023sweep}, the micro-macro decomposition \cite{einkemmer2021asymptotic,einkemmer2024asymptotic,ceruti2025galerkin}, and the high-order low-order (HOLO) framework  \cite{peng2021high}, with theoretical analysis in \cite{ding2021dynamical,yin2025towards}. Rank adaptive DLRAs for time-dependent RTE have been designed \cite{hauck2023predictor,baumann2024energy,frank2025asymptotic,ceruti2024parallel} with an extension to inverse problems \cite{baumann2025adaptive}. DLRAs have also been adapted to an energy dependent steady-state problem \cite{kusch2021robust} and an eigenvalue problem \cite{kusch2022low} by interpreting the energy discretization and power iteration as pseudo time-marching. 
Beyond DLRA, the step-and-truncation (SAT) method \cite{guo2022low,nakao2025reduced}, which allows for natural rank adaptivity, has been developed for time-dependent RTE based on micro-macro decomposition \cite{sands2024high}. For the second-order formulation of 1D1V steady-state RTE, a low-rank iterative method based on soft-thresholding \cite{bachmayr2017iterative} and exponential sums preconditioner \cite{ballani2013projection,hackbusch2015solution} has been proposed \cite{bachmayr2024low}. Proper generalized decomposition, utilizing CP tensor decomposition and alternating direction optimization, has been applied for RTE \cite{alberti2020reduced,dominesey2022reduced,dominesey2023reduced}. For a broader overview of low-rank methods for other kinetic equations, see the review paper \cite{einkemmer2025review}.

While significant progress has been made in developing low-rank methods for time-dependent RTE, research on low-rank iterative solvers and preconditioners for steady-state problems remains limited. The recent work \cite{bachmayr2024low} designs and rigorously analyzes a soft-thresholding based low-rank iterative solver for the 1D1V problem. In this work, we consider the multidimensional steady-state RTE in its second-order formulation and develop an inexact low-rank SI-DSA solver. To the best of our knowledge, this is the first such solver designed for multidimensional RTE.

Before outlining the proposed methodology, we briefly review low-rank iterative solvers for linear systems that are relevant to our development.
Such solvers have been designed based on alternating optimization \cite{holtz2012alternating,dolgov2014alternating}, Richardson iteration  \cite{matthies2012solving,bachmayr2017iterative}, low-rank Krylov method \cite{kressner2011low,ballani2013projection,dolgov2013tt} and Riemannian optimization \cite{kressner2016preconditioned}. Preconditioners for low-rank solvers have been developed by leveraging the exponential sums \cite{braess2005approximation,hackbusch2006low,beylkin2010approximation,hackbusch2015solution}, multigrid \cite{ballani2013projection,hackbusch2015solution}, averaged-coefficient
Sylvester preconditioner \cite{el2024sylvester} and basis-update-Galerkin approach \cite{meng2024preconditioning}. A preconditioned low-rank Anderson acceleration has been recently introduced for nonlinear problems \cite{appelo2025lraa}. Instead of an exhausting literature review, we refer to the recent review paper \cite{bachmayr2023low}.

The contributions of this paper are twofold:
\begin{enumerate}
    \item {\bf Low-rank SI-DSA framework}: We extend SI-DSA to the low-rank setting, representing the even parity in the second-order form of RTE and all intermediate high-dimensional quantities via low-rank matrix or hierarchical Tucker tensor (HTT) decompositions \cite{grasedyck2010hierarchical}, and solving each SI step with a preconditioned low-rank conjugate gradient (CG) method \cite{kressner2011low}.
    
    \item {\bf Adaptive inexact inner loop solves}: To further accelerate our method, we develop an adaptive strategy that dynamically adjusts the stopping criteria for the inner-loop low-rank CG solve, balancing accuracy and efficiency.
\end{enumerate}
We demonstrate the performance of our method through a series of 2D2V benchmark tests.  Even for challenging problems where the effective rank grows with grid resolution, our low-rank approach still achieves a substantial reduction in the required degrees of freedom (DOFs) compared to full-rank solvers.

In addition to considering multidimensional cases, our preconditioning strategy differs from the method for 1D1V problem in  \cite{bachmayr2017iterative}. We adopt the DSA framework, whereas \cite{bachmayr2024low} employs an exponential sums preconditioner \cite{ballani2013projection,hackbusch2015solution}, which in turn requires a Cholesky decomposition to transform the system into a Kronecker-sum form. However, in high dimensions, unless a tensor product grid is used and spatially separable material properties are assumed, both the memory requirement and the computational cost of Cholesky decomposition grow significantly due to the increased bandwidth of discrete spatial operators.

Although our focus is on DSA, the proposed low-rank approach can, in principle, be combined with other full-rank SA preconditioners, such as the quasi-diffusion method \cite{gol1964quasi,anistratov1993nonlinear,olivier2023family}, S2SA \cite{lorence1989s} or more recent reduced order model enhanced preconditioners \cite{mcclarren2022data,peng2024reduced,peng2025flexible}. Exploring such extensions is left for future work.

Our paper is organized as follows. In Sec. \ref{sec:background}, we introduce the model problem, our discretizations and briefly review the basic idea of SI-DSA. In Sec. \ref{sec:lr-method}, we develop our inexact low-rank SI-DSA. In Sec. \ref{sec:numerical}, we numerically demonstrate the performance of our method through a series of benchmark tests. At last, we draw our conclusions in Sec. \ref{sec:conclusion}.

\section{Background\label{sec:background}}
In this paper, we consider steady state, single group, linear RTE with isotropic scattering, source and vacuum inflow boundary conditions on the computational domain $\Gamma_{\bx}$:
\begin{subequations}
\label{eq:rte}
    \begin{align}
    &\BOmega \cdot \nabla \psi(\bx,\BOmega) + \sigma_t(\bx) \psi(\bx,\BOmega) = \sigma_s(\bx) \phi(\bx) + G(\bx), \quad\sigma_t(\bx)=\sigma_s(\bx)+\sigma_a(\bx),\\
    &\phi(\bx) = \frac{1}{4\pi}\int_{\mathbb{S}^2} \psi(\bx,\BOmega) d\BOmega, \quad \bx\in\Gamma_{\bx},\\
    &\psi(\bx,\BOmega) = 0, \quad \bx\in \partial\Gamma_{\bx},\;\BOmega\cdot \bn(\bx)<0.\label{eq:boundary-condition}
    \end{align}
\end{subequations}
Here, $\psi(\bx,\BOmega)$ is the angular flux (also known as radiation intensity or particle distribution depending on the underlying application) corresponding to angular direction on the unit sphere $\BOmega=(\BOmega_x,\BOmega_y,\BOmega_z)\in\mathbb{S}^2$ and spatial location $\bx=(x,y,z)\in\Gamma_{\bx}$,  $\phi(\bx)=\frac{1}{4\pi}\int_{\mathbb{S}^2} \psi(\bx,\BOmega) d\BOmega$ is the scalar flux (also known as the macroscopic density), $G(\bx)$ is an isotropic source, $\sigma_s(\bx)\geq 0$, $\sigma_a(\bx)\geq 0$ and $\sigma_t(\bx)=\sigma_s(\bx)+\sigma_a(\bx)\geq 0$ are the scattering, absorption and total cross sections, respectively. In the inflow boundary condition \eqref{eq:boundary-condition}, $\bn(\bx)$ is the outward normal direction of the computational domain $\Gamma_{\bx}$ at location $\bx\in\partial\Gamma_{\bx}$.

\subsection{Second-order formulation of RTE\label{sec:rte-2nd}}
We solve the RTE based on its second-order formulation and refer to the review paper \cite{lewis2009second} for an overview of numerical methods for second-order RTE. We start by defining the even-parity, $\psi^{+}$, and the odd-parity, $\psi^{-}$, as
\begin{equation}
    \psi^{\pm}(\bx,\BOmega)=\frac{1}{2}\left(\psi(\bx,\BOmega)\pm\psi(\bx,-\BOmega)\right).
\end{equation}
To derive equations for $\psi^{\pm}$, we first
decompose the angular space $\mathbb{S}^2$ as $\mathbb{S}^2=\mathbb{S}^+\oplus\mathbb{S}^-$, where $\mathbb{S}^\pm=\{\BOmega\in\mathbb{S}^2,\;\pm\BOmega_y\geq 0\}.$ Then, equation \eqref{eq:rte} can be rewritten as 
\begin{subequations}
    \begin{align}   \BOmega\cdot\nabla\psi(\bx,\BOmega)+\sigma_t\psi(\bx,\BOmega)&=\sigma_s\phi(\bx)+G(\bx),\quad\BOmega\in\mathbb{S}^+,\label{eq:rte_splus}\\
    -\BOmega\cdot\nabla\psi(\bx,-\BOmega)+\sigma_t\psi(\bx,-\BOmega)&=\sigma_s\phi(\bx)+G(\bx),\quad\BOmega\in\mathbb{S}^+.\label{eq:rte_sminus}
    \end{align}
\end{subequations}
By summing equations \eqref{eq:rte_splus} and \eqref{eq:rte_sminus}, and subtracting \eqref{eq:rte_sminus} from \eqref{eq:rte_splus}, we derive the following coupled system:
\begin{subequations}
    \begin{align}
        \BOmega\cdot\nabla\psi^-(\bx,\BOmega)+\sigma_t\psi^+(\bx,\BOmega)=\sigma_s\phi(\bx)+G(\bx),\quad\BOmega\in\mathbb{S}^+,
        \label{eq:even-odd1}\\
        \BOmega\cdot\nabla\psi^+(\bx,\BOmega)+\sigma_t\psi^{-}(\bx,\BOmega)=0,\quad\BOmega\in\mathbb{S}^+.
    \label{eq:even-odd2}
    \end{align}
\end{subequations}
When $\sigma_t>0$, equation \eqref{eq:even-odd2} leads to $\psi^{-}=-\sigma_t^{-1}\BOmega\cdot\nabla \psi^{+}$. Substituting this relation into equation \eqref{eq:even-odd1}, we can eliminate the odd-parity $\psi^-$ and obtain the following second-order equation for the even-parity $\psi^{+}$:
\begin{equation}
    \left(-(\BOmega\BOmega^T):\left(\nabla(\sigma_t^{-1}\nabla)^T\right)+\sigma_t\right)\psi^+(\bx,\BOmega)=\sigma_s\phi(\bx)+G(\bx),\quad\BOmega\in\mathbb{S}^+,
    \label{eq:even-parity}
\end{equation}
where $A:B=\textrm{tr}(AB^T)$ represents matrix inner product in the Frobenius sense. For simplicity of notation, we define the second-order operator $T$ as 
\begin{equation}
T = -(\BOmega\BOmega^T):\left(\nabla(\sigma_t^{-1}\nabla)^T\right).
    \label{eq:transport-operator}
\end{equation}
Then, equation \eqref{eq:even-parity} can be written as 
\begin{equation}
    (T+\sigma_t)\psi^{+}=\sigma_s\phi+G.
    \label{eq:even-parity-operator}
\end{equation}

In this paper, we focus on the 2D X-Y geometry, assuming that $\partial_z\psi^{+}=0$. Under this assumption, equation \eqref{eq:even-parity} reduces to the following form in the 2D2V setting.
\begin{equation}
    -\Big(\BOmega_x^2\partial_x(\sigma_t^{-1}\partial_x)+\BOmega_x\BOmega_y(\partial_x\left(\sigma_t^{-1}\partial_y)+\partial_y(\sigma_t^{-1}\partial_x)\right)+\BOmega_y^2\partial_y(\sigma_t^{-1}\partial_y)\Big)\psi^++\sigma_t\psi^+=\sigma_s\phi+G,\quad\BOmega\in\mathbb{S}^+.
    \label{eq:even-parity-2D}
\end{equation}
\subsection{Discretization\label{sec:discretization}}
We focus on the 2D X-Y geometry and apply the discrete ordinates method ($S_N$ method) in the angular space, combined with a finite difference method in the physical space.
\subsubsection{Angular discretization\label{sec:angular-discretization}} 
The discrete ordinates angular discretization \cite{pomraning2005equations} solves the RTE at a prescribed set of quadrature points in angular space. Specifically, we employ the Chebyshev–Legendre (CL) quadrature. The normalized CL quadrature rule, denoted by \(\mathrm{CL}(N_\theta, N_{\BOmega_z})\), is constructed as the tensor product of the normalized \(N_{\BOmega_z}\)-point Gauss–Legendre quadrature on [-1,1], $\{(\BOmega_{z,j},\omega_{z,j})\}_{j=1}^{N_{\BOmega_z}}$ with $\sum_{j=1}^{N_{\BOmega_z}}\omega_j=1$, and
the normalized $N_\theta$-point Chebyshev quadrature on the unit circle.
\begin{equation}
\left\{(\theta_j,\omega_{\theta,j}): \;\theta_j = \frac{2j\pi}{N_\theta}-\frac{\pi}{N_\theta}\,\text{and}\,\omega^\phi_j=\frac{1}{N_\theta},\; j=1,\dots,N_\theta\right\}.
\end{equation}
Then, with a slight abuse of notation, the quadrature points and the associated quadrature weights of the CL($N_{\theta}$,$N_{\BOmega_z}$) rule are defined as
\begin{equation}
\BOmega_j=\left(\BOmega_{j,x},\BOmega_{j,y},\BOmega_{j,z}\right)=\left(\cos(\theta_{j_1})\sqrt{1-\BOmega_{z,j_2}^2},\sin(\theta_{j_1})\sqrt{1-\BOmega_{z,j_2}^2},\BOmega_{z,j_2}\right)\quad\text{and}\quad\omega_j=\omega_{\theta,j_1}\omega_{\BOmega_z,j_2},
\end{equation}
with $j=(j_1,j_2)$, $j_1 = 1,\ldots, N_\theta$, $j_2 =1,\ldots,N_{\BOmega_z}$.
Note that if $N_\theta=2N$ and $N_{\BOmega_z}=N$,  then the quadrature rule \(\mathrm{CL}(N_\theta, N_{\BOmega_z})\) exactly integrates normalized polynomials of degree up to $2N$ on the unit sphere.  Furthermore, if $\BOmega_j$ is a CL quadrature point then $-\BOmega_j$ is also a quadrature point.

The operator form of the angular discretized second-order RTE is derived by approximating the scalar flux with numerical quadrature. We seek $\psi^+(\cdot,\BOmega_j)\approx\psi^+_j$ that solves
\begin{equation}
 (T_j+\sigma_t)\psi^+_j=\sigma_s\phi+G,\quad\text{for}\quad \BOmega_{j,y}>0,
\end{equation}
where $\phi=\sum_{\BOmega_{j,y}>0}\omega_j\psi_j^+$ and $T_j = -(\BOmega_j\BOmega_j^T):\left(\nabla(\sigma_t^{-1}\nabla)^T\right).$
\subsubsection{Finite difference spatial discretization\label{sec:spatial-discretization}}
In 2D X-Y geometry, we consider a rectangular computational domain $[x_L,x_R]\times[y_B,y_T]$.
Equation~\eqref{eq:even-parity-2D} is discretized using a central difference scheme on a uniform Cartesian grid with $N_x+1$ points in the $x$-direction and $N_y+1$ points in the $y$-direction.
Denote the grid size as $\dx=\frac{x_R-x_L}{N_x}$ and $\dy=\frac{y_T-y_B}{N_y}$. Define grid points as $x_{i_x}=x_L+i_x\Delta x$ and $y_{i_y}=y_B+i_y\Delta y$, where $i_x=0,\dots,N_x$ and $i_y=0,\dots,N_y$.
Let $\psi^+(x_{i_x},y_{i_y},\BOmega_j)\approx\psi^{+}_{i_x,i_y,j}$ and $\phi(x_{i_x},y_{i_y})\approx\phi_{i_x,i_y}$.
Then, the central difference $S_N$ approximation for the second-order system in the 2D X-Y geometry \eqref{eq:even-parity-2D} is defined as
\begin{equation}
\label{eq:even-parity-2D-FD}
-\Big(\BOmega_{j,x}^2 \BmcD_{xx}^{\sigma_t}+\BOmega_{j,x}\BOmega_{j,y}(\BmcD^{\sigma_t}_{xy}+\BmcD^{\sigma_t}_{yx})+\BOmega_{j,y}^2\BmcD^{\sigma_t}_{yy}\Big)\psi^+_{i_x,i_y,j}+\sigma_{t,i_x,i_y}\psi^+_{i_x,i_y,j}
=\sigma_{s,i_x,i_y}\phi_{i_x,i_y}+G_{i_x,i_y}
\end{equation}
with $\sigma_{t,i_x,i_y}=\sigma_t(x_{i_x},y_{i_y})$, $\sigma_{s,i_x,i_y}=\sigma_s(x_{i_x},y_{i_y})$ and $G_{i_x,i_y}=G(x_{i_x},y_{i_y})$.
Here, the finite difference operators are defined as
\begin{subequations}
\label{eq:fd-operators}
\begin{align}
\BmcD^{\sigma_t}_{xx}\psi^+_{i_x,i_y,j}&=\frac{1}{\dx^2}\left(\sigma^{-1}_{t,i_x-\half,i_y}\psi^{+}_{i_x-1,i_y,j}+(\sigma^{-1}_{t,i_x-\half,i_y}+\sigma^{-1}_{t,i_x+\half,i_y})\psi^{+}_{i_x,i_y,j}+\sigma^{-1}_{t,i_x+\half,i_y}\psi^{+}_{i_x+1,i_y,j}\right), \label{eq:dxx}\\
\BmcD^{\sigma_t}_{yy}\psi^+_{i_x,i_y,j}&=\frac{1}{\dy^2}\left(\sigma^{-1}_{t,i_x,i_y-\half}\psi^{+}_{i_x,i_y-1,j}+(\sigma^{-1}_{t,i_x,i_y-\half}+\sigma^{-1}_{t,i_x,i_y+\half})\psi^{+}_{i_x,i_y,j}+\sigma^{-1}_{t,i_x,i_y+\half}\psi^{+}_{i_x,i_y+1,j}\right), \label{eq:dyy}\\
\BmcD^{\sigma_t}_{xy}\psi^+_{i_x,i_y,j}&=\frac{1}{4\dx\dy}\left(\sigma^{-1}_{t,i_x+1,i_y}(\psi^{+}_{i_x+1,i_y+1,j}-\psi^{+}_{i_x+1,i_y-1,j})-\sigma^{-1}_{t,i_x-1,i_y}(\psi^{+}_{i_x-1,i_y+1,j}-\psi^{+}_{i_x-1,i_y-1,j})\right), \label{eq:dxy}\\
\BmcD^{\sigma_t}_{yx}\psi^+_{i_x,i_y,j}&=\frac{1}{4\dx\dy}\left(\sigma^{-1}_{t,i_x,i_y+1}(\psi^{+}_{i_x+1,i_y+1,j}-\psi^{+}_{i_x-1,i_y+1,j})-\sigma^{-1}_{t,i_x,i_y-1}(\psi^{+}_{i_x+1,i_y-1,j}-\psi^{+}_{i_x-1,i_y-1,j})\right), \label{eq:dyx}
\end{align}
\end{subequations}
where the grid function $\sigma^{-1}_{t,\alpha,\beta}$ is defined as $\sigma^{-1}_{t,\alpha,\beta}=\frac{1}{\sigma_t(x_{\alpha},y_{\beta})}$ with $\alpha=i_x\pm1,i_x\pm\half,i_x$; $\beta=i_y\pm1,i_y\pm\half,i_y$;  $x_\alpha=x_L+\alpha\dx$ and $y_\beta=y_B+\beta\dy$.

Define the total number of spatial DOFs as $N_{\bx}$.
Let $\bpsi^+=((\bpsi^+_1)^T,\dots,(\bpsi^+_{N_{\BOmega}/2})^T)\in\mathbb{R}^{N_{\bx}N_{\BOmega}/2}$, where $\bpsi_j^+$ stores DOFs for $\psi^+(x,y,\BOmega_j)$. Let $\bphi\in\mathbb{R}^{N_{\bx}}$ denote the vector storing DOFs for $\phi$. The matrix-vector formulation for the fully discretized system has a tensor product structure with respect to its spatial and angular DOFs:
\begin{align}
    \left(-(\BD_{\BOmega_x^2}\otimes\BD^{\sigma_t}_{xx}+\BD_{\BOmega_{x}\BOmega_y}\otimes\BD^{\sigma_t}_{\textrm{cross}}+\BD_{\BOmega_y^2}\otimes\BD^{\sigma_t}_{yy})+\BI\otimes\BSigma_t\right)\bpsi^+=(\BI\otimes\BSigma_s)\bphi+\boldsymbol{1}\otimes G,
    \label{eq:tensor-system}
\end{align}
where $\otimes$ denotes the Kronecker product, $\BD^{\sigma_t}_{xx}$, $\BD^{\sigma_t}_{yy}$  and $\BD^{\sigma_t}_{\textrm{cross}}$ represent matrices corresponding to the finite difference operators $\BmcD_{xx}^{\sigma_t}$, $\BmcD_{yy}^{\sigma_t}$ and $\BmcD_{xy}^{\sigma_t}+\BmcD_{yx}^{\sigma_t}$, $\BSigma_t,\BSigma_s\in\mathbb{R}^{N_{\bx}\times N_{\bx}}$ are diagonal matrices whose diagonal elements are given by the values of $\sigma_{t}$ and $\sigma_{s}$ evaluated at spatial grid points, $\BD_{\BOmega_x^2}, \BD_{\BOmega_x\BOmega_y}, \BD_{\BOmega_y^2}\in\mathbb{R}^{N_{\BOmega}/2\times N_{\BOmega}/2}$ are diagonal matrices whose diagonal elements are the values of $\BOmega_{x}^2$, $\BOmega_{y}^2$ and $\BOmega_{x}\BOmega_y$ evaluated at quadrature points in angular space,
$\BI\in\mathbb{R}^{N_{\BOmega}/2\times N_{\BOmega}/2}$ is the identity matrix, and $\boldsymbol{1}=(1,1,\dots,1)^T\in\mathbb{R}^{N_{\BOmega}/2}$. One can check that 
$ -(\BD_{\BOmega_x^2}\otimes\BD^{\sigma_t}_{xx}+\BD_{\BOmega_{x}\BOmega_y}\otimes\BD^{\sigma_t}_{\textrm{cross}}+\BD_{\BOmega_y^2}\otimes\BD^{\sigma_t}_{yy})+\BI\otimes\BSigma_t$ is a symmetric positive definite (SPD) matrix.

\subsection{Source Iteration with Diffusion Synthetic Acceleration (SI-DSA)\label{sec:si-dsa}}
We briefly review SI-DSA for solving RTE. For more details, we refer to review papers \cite{Adams2002FastIM,larsen2009advances} and DSA for the second-order formulation of RTE \cite{palii2020convergent}.

Let the initial guess of the scalar flux be $\bphi^{(0)}$. The $k$-th iteration of SI-DSA for the second-order equation parity \eqref{eq:even-parity-operator} consists of two main steps:
\begin{enumerate}
    \item \textbf{SI.} 
    For each angular direction, $\BOmega_j$, $1\leq j\leq N_{\BOmega}$, solve the angular decoupled equation
    \begin{equation}
    \left(-(\BOmega_{j,x}^2\BD^{\sigma_t}_{xx}+\BOmega_{j,x}\BOmega_{j,y}\BD^{\sigma_t}_{\textrm{cross}}+\BOmega_{j,y}^2\BD^{\sigma_t}_{yy})+\BSigma_t\right)\bpsi_j^{+,(k)}=\BSigma_s\bphi^{(k-1)}+G. \label{eq:SI}  
    \end{equation}
    \item \textbf{DSA.} Correct the scalar flux as $\bphi^{(k)}=\bphi^{(k,*)}+\delta\bphi^{(k)}$ by solving the diffusion limit of the error equation for the $k$-th iteration:
    \begin{equation}
        \big(-\frac{1}{3}(\BD^{\sigma_t}_{xx}+\BD^{\sigma_t}_{yy})+\BSigma_a\big)\delta\bphi^{(k)}=\BSigma_s(\bphi^{(k,*)}-\bphi^{(k)}),
    \label{eq:DSA}
    \end{equation} 
    where $\BSigma_a=\BSigma_t-\BSigma_s$ denotes the matrix for the absorption operator.
\end{enumerate}
Overall, SI-DSA can be interpreted as a preconditioned Richardson iteration, where DSA serves as the crucial preconditioning step. Without a preconditioner, the convergence speed of SI may become arbitrarily slow when $\sigma_s/\sigma_t \approx 1$ \cite{Adams2002FastIM}.

\section{Methodology\label{sec:lr-method}}
Two main challenges in solving steady-state RTEs are the substantial memory requirement arising from its high dimensionality and the design of effective preconditioners for multiscale problems. To address these issues, we extend the SI-DSA to a low-rank framework, leveraging recent advances in matrix and tensor decompositions to reduce storage costs. A key step is to take advantage of the tensor-product structure of the discrete linear system to develop an efficient low-rank solver.


In Sec.\ref{sec:lr-representation}, we present low-rank representations of the solution based on matrix and HTT decompositions \cite{grasedyck2010hierarchical} and describe how operators are applied within the low-rank framework. Sec.~\ref{sec:lr-SI-DSA} introduces our inexact low-rank SI–DSA algorithm, while Sec.~\ref{sec:lr-inner} provides details of the inner-loop solver used in each SI step.

\subsection{Low-rank representation\label{sec:lr-representation}}
\begin{figure}[]
  \begin{center} 
\includegraphics[width=0.6\textwidth]{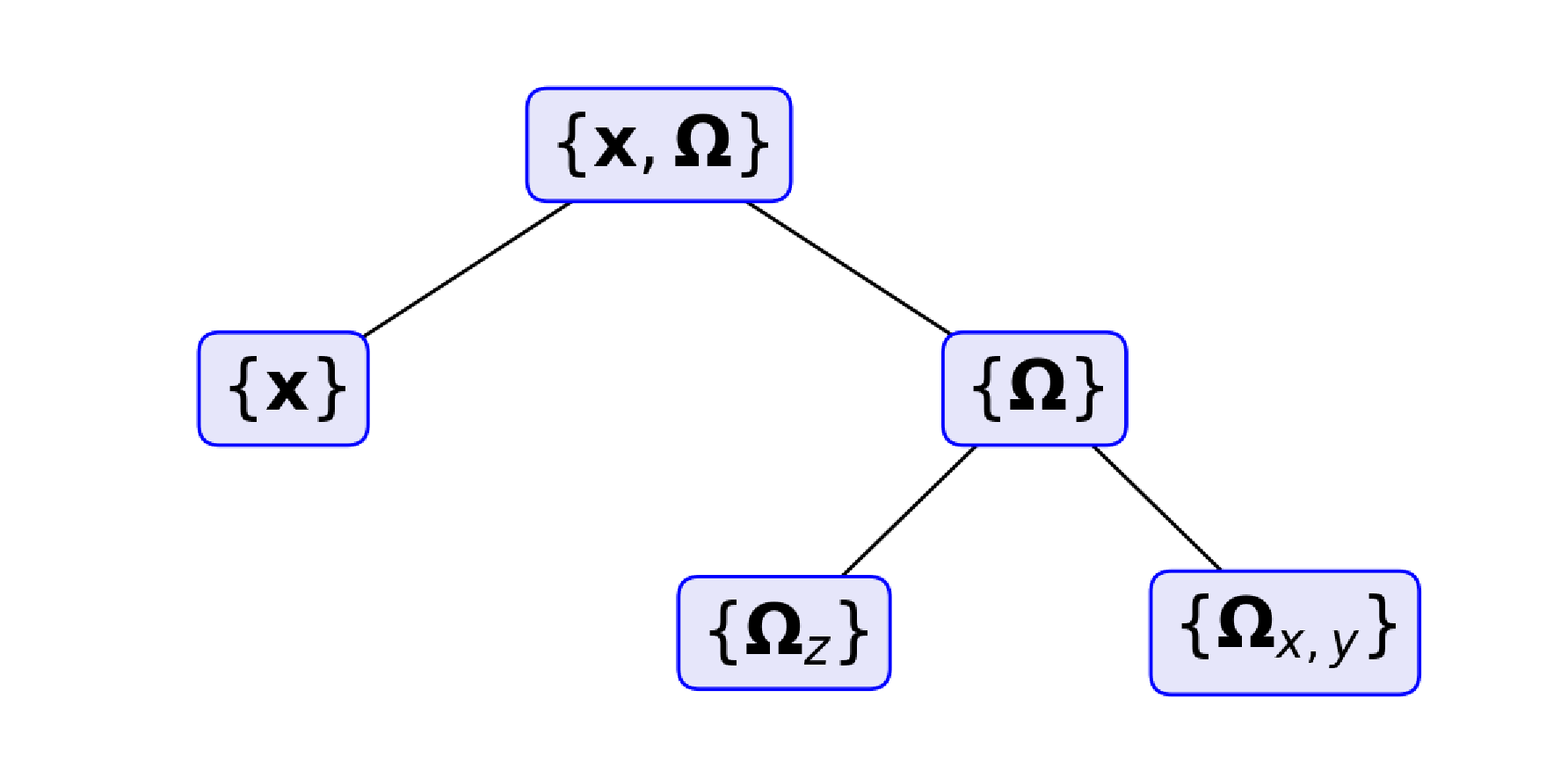}
  \caption{Dimension tree for the HTT decomposition  \eqref{eq:continuous-tensor}. \label{fig:dim_tree}}
  \end{center}
\end{figure}

We consider the following low-rank representations for the even-parity $\psi^+$. The first format decomposes the phase space as
\begin{equation}
\psi^+(\bx,\BOmega) \approx \sum_{l=1}^r \BB_{l}\,\BX_l(\bx)\BV_l(\BOmega),\label{eq:continuous-matrix},
\end{equation}
which is closely related to the Hilbert–Schmidt decomposition. Here, $r$ is the rank, $\{\BX_l\}_{l=1}^r$ and $\{\BV_l\}_{l=1}^r$ are the bases in physical and angular spaces, respectively, and $\{\BB_l\}_{l=1}^r$ are the coefficients. The second format is the HTT representation
\begin{subequations}
\label{eq:continuous-tensor}
\begin{align}
\psi^+(\bx,\BOmega) \approx
\sum_{l_1=1}^{r_1}\sum_{l_{23}=1}^{r_{23}} \BB^{123}_{l_1,l_{23}}\BX_{l_1}(\bx)\BV^{23}_{l_{23}}(\BOmega),\\
\BV^{23}_{l_{23}}(\BOmega)=\sum_{l_2=1}^{r_2}\sum_{l_3=1}^{r_3}\BB^{23}_{l_2,l_3,l_{23}}\BV^2_{l_2}(\BOmega_z)\BV^3_{l_3}(\BOmega_x,\BOmega_y),
\end{align}  
\end{subequations}
which further decomposes the angular variables. $r_1$, $r_{23}$, $r_2$, and $r_3$ form the hierarchical ranks of the HTT representation. $\BB^{123}$ and  $\BB^{23}$ are called transfer tensors. Note that the basis $\{\BV^{23}_{l_{23}}\}_{l_{23}=1}^{r_{23}}$ is defined through the lower-dimensional bases $\{\BV^2_{l_2}\}_{l_2=1}^{r_2}$ and $\{\BV^3_{l_3}\}_{l_3=1}^{r_3}$. In general, the HTT decomposition, also referred to as a tree tensor network, arranges the problem dimensions in a binary tree. The root node corresponds to the full set of dimensions, which is recursively divided into disjoint subsets represented by the non-root nodes. Figure~\ref{fig:dim_tree} illustrates the dimension tree associated with the HTT decomposition  \eqref{eq:continuous-tensor}.

At the discrete level, within the finite difference setting, we have the corresponding matrix decomposition 
\begin{equation}
\bpsi^+ \approx \BX\BB\BV^T=\sum_{l=1}^r\BB_l\,\BX_l\otimes\BV_l,\;\BX\in\mathbb{R}^{N_{\bx}\times r}, \BB\in\mathbb{R}^{r}, \BV\in\mathbb{R}^{N_{\BOmega}\times r},
\label{eq:lr-matrix}
\end{equation}
and HTT decomposition 
\begin{subequations}
\label{eq:lr-tensor}
\begin{align}
\bpsi^+\approx
\sum_{l_1=1}^{r_1}\sum_{l_{23}=1}^{r_{23}} \BB^{123}_{l_1,l_{23}}\BX_{l_1}\otimes\BV_{l_{23}}^{23},\; \BB_{l_1,l_{23}}^{123}\in\mathbb{R},\,\BX_{l_1}\in\mathbb{R}^{N_{\bx}},\,\BV_{l_{23}}^{23}\in\mathbb{R}^{N_{\BOmega}/2},\\
\BV^{23}_{l_{23}}=\sum_{l_2=1}^{r_2}\sum_{l_3=1}^{r_3}\BB^{23}_{l_2,l_3,l_{23}}\BV^{2}_{l_2}\otimes\BV^3_{l_3},\;\BB_{l_2,l_3,l_{23}}^{23}\in\mathbb{R},\,\BV_{l_2}^2\in\mathbb{R}^{N_{\BOmega_z}},\,\BV_{l_3}^3\in\mathbb{R}^{N_\theta}.
\end{align}  
\end{subequations}

\noindent\textbf{Memory compression.} 
In the full-rank format, the number of DOFs scales as $O(N_{\bx}N_{\BOmega})=O(N_{\bx}N_{\BOmega_z}N_{\theta})$. In contrast, the DOF count reduces to $O\left((N_{\bx} + N_{\BOmega}) r\right)$ in the low-rank matrix format, and to $O(N_{\bx}r_1+N_{\BOmega_z}r_2+N_\theta r_3+r_1r_{23}+r_2r_3r_{23})$ in the HTT format.


\noindent\textbf{Low-rank computation of operators.} For illustrative purposes, we focus on the matrix format, while also outlining how the approach can be extended to the HTT format. For example, 
matrix-vector multiplication in equation \ref{eq:tensor-system} $(\BD_{\BOmega_x^2}\otimes\BD_{xx}^{\sigma_t})\bpsi^+$ can be computed in the low-rank manner by applying the difference operators to the basis in the respective dimension, yielding
$$\sum_{l=1}^r\BB_l\,(\BD_{\BOmega_x^2}\BV_l)\otimes(\BD_{xx}^{\sigma_t}\BX_l),$$ which costs $O((N_{\bx}+N_{\BOmega})r)$.  
Then, the summation of low-rank outputs in \eqref{eq:tensor-system},  including $\sum_{l=1}^r\BB_l\,(\BD_{\BOmega_x^2}\BV_l)\otimes(\BD_{xx}^{\sigma_t}\BX_l)$, $\sum_{l=1}^r\BB_l\,(\BD_{\BOmega_y^2}\BV_l)\otimes(\BD_{yy}^{\sigma_t}\BX_l)$,
$\sum_{l=1}^r\BB_l\,(\BD_{\BOmega_x\BOmega_y}\BV_l)\otimes(\BD_{\textrm{cross}}^{\sigma_t}\BX_l)$,
can be approximated using a truncated sum of low-rank matrices. Note that truncated sum also controls the rank of the summation result by performing an SVD of the low-dimensional matrix cores. The overall cost for the truncation scales as $O((N_{\bx}+N_{\BOmega})r^2+r^3)$. See, e.g., Sec. 2.1 of \cite{guo2022low} or Alg. 2.3 of \cite{meng2024preconditioning} for more details. 

For the HTT format, discrete operators in the angular space can be further decomposed into the Kronecker product. For example, $\BD_{\BOmega_x}^2=\BD_{\cos^2(\theta)}\otimes\BD_{1-\BOmega_z^2}$ where $\BD_{\cos^2(\theta)}\in\mathbb{R}^{N_\theta\times N_\theta}$ and $\BD_{1-\BOmega_z^2}\in\mathbb{R}^{N_{\BOmega_z}\times N_{\BOmega_z}}$ are diagonal matrices whose diagonal entries are the point values of $\cos^2(\theta)$ and $1-\BOmega_z^2$, respectively, evaluated at the corresponding Chebyshev and Legendre points. The truncated sum of HTT controlling the rank of summation results is realized by leveraging hierarchical higher order SVD (HOSVD) \cite{grasedyck2010hierarchical}. We refer to \cite{bachmayr2023low,kressner2012htucker} for details of the truncated summation algorithm for HTT. In the HTT format, the computational cost of applying operators is $O(N_{\mathbf{x}} r_1 + N_{\boldsymbol{\Omega}z} r_2 + N_\theta r_3)$, while that of performing truncated tensor summation is $O(d N r^2 + d r^4)$, $r=\max(r_1,r_2,r_3,r_{23})$, $N=\max(N_{\bx},N_{\BOmega_z},N_\theta)$ \cite{grasedyck2010hierarchical}.

\subsection{Inexact low-rank SI-DSA\label{sec:lr-SI-DSA}}
We are now ready to formulate the main algorithm—an inexact low-rank variant of SI-DSA—for solving the second-order formulation of the RTE, with the goal of reducing computational cost.
The algorithm is summarized in Alg. \ref{alg:lr-si-dsa}.
\begin{algorithm}[ht]
\caption{Inexact low-rank SI-DSA.\label{alg:lr-si-dsa}}
\begin{algorithmic}[1]
\STATE{\textbf{Input}: threshold $\epsilon_{\textrm{Diff}}$, the maximum number of source iterations $N_{\textrm{SI-Iter}}$, initial tolerance for the relative residual of the inner loop low-rank CG solver $\epsilon_{\textrm{CG}}^{(1)}$, the maximum iteration number for inner loop low-rank CG solver $N_{\textrm{CG-Iter}}$, decay rate $0<\gamma<1$, and initial guess for the scalar flux $\bphi^{(0)}$.}
\STATE{\textbf{Initialization}: set $k=1$.}
\WHILE{$k\leq N_{\textrm{SI-Iter}}$}
\STATE{\textbf{SI step}: solve 
\begin{equation}
    \left(-(\BD_{\BOmega_x^2}\otimes\BD^{\sigma_t}_{xx}+\BD_{\BOmega_{x}\BOmega_y}\otimes\BD^{\sigma_t}_{\textrm{cross}}+\BD_{\BOmega_y^2}\otimes\BD^{\sigma_t}_{yy})+\BI\otimes\BSigma_t\right)\bpsi^{+,(k)}=(\BI\otimes\BSigma_s)\bphi^{(k-1)}+\boldsymbol{1}\otimes G, \label{eq:lr-si-step}
\end{equation}
with the low-rank CG algorithm in Alg. \ref{alg:lr-cg} until the relative residual is smaller than $\epsilon_{\textrm{CG}}^{(k)}$ or the inner loop iteration number reaches $N_{\textrm{CG-Iter}}$.
}
\STATE{\textbf{Update the scalar flux} $\bphi^{(k,*)}$ by numerically integrating $\bpsi^{+,(k)}$.}
\IF{$||\bphi^{(k,*)}-\bphi^{(k-1)}||_{\infty}\leq\epsilon_{\textrm{Diff}}$}
    \BREAK
\ENDIF
\STATE{\textbf{DSA step}: apply AMG to solve the diffusion correction equation
\begin{equation}
    \big(-\frac{1}{3}(\BD^{\sigma_t}_{xx}+\BD^{\sigma_t}_{yy})+\BSigma_a\big)\delta\bphi^{(k)}=\BSigma_s(\bphi^{(k,*)}-\bphi^{(k)}),
    \label{eq:lr-dsa}
\end{equation}
where $\BSigma_a=\BSigma_t-\BSigma_s$.
}
\STATE{Correct the scalar flux as $\bphi^{(k+1)}=\bphi^{(k,*)}+\delta\bphi^{(k)}$, and then set $k:=k+1$.}
\STATE{\textbf{Update the inner loop accuracy:} set $\epsilon_{\textrm{CG}}^{(k)}=\min(\epsilon_{\textrm{CG}}^{(1)},\gamma||\bphi^{(k,*)}-\bphi^{(k)}||)$.}
\ENDWHILE
\end{algorithmic}
\end{algorithm}

The key steps of the algorithm are outlined below:
\begin{enumerate}
    \item Since $-(\BD_{\BOmega_x^2}\otimes\BD^{\sigma_t}_{xx}+\BD_{\BOmega_{x}\BOmega_y}\otimes\BD^{\sigma_t}_{\textrm{cross}}+\BD_{\BOmega_y^2}\otimes\BD^{\sigma_t}_{yy})+\BI\otimes\BSigma_t$ is SPD, a preconditioned low-rank conjugate gradient (CG) \cite{kressner2011low} solver is applied to solve the inner loop equation \eqref{eq:lr-si-step}. Details of the inner-loop low-rank solver and the preconditioner for it will be discussed in Sec. \ref{sec:lr-inner}.
    \item The numerical integration of $\bphi^{(k,*)}=\sum_{\BOmega_j\in\mathbb{S}^+}\omega_j\bpsi_{j}^{(k)}$ is computed in the low-rank manner through tensor contraction with the cost $O((N_{\bx}+N_{\BOmega})r)$ for the matrix format and $O(r(N_{\bx}+N_{\Omega_z}+N_{\theta}))$ for the HTT format. 
    \item Following the DSA framework, we use the diffusion limit of the discrete error equation as our preconditioner.
    \item We adaptively control the required accuracy of the inner loop low-rank CG solver to further reduce the computational cost. The motivation to introduce this step will be detailed in Sec. \ref{sec:lr-inner}. 
\end{enumerate}
Overall, our algorithm can be seen as a preconditioned inexact fixed point iteration in the low-rank setting. As discussed in \cite{birken2015termination}, an inexact fixed point iteration will converge, if the inexactness in the inner loop solver approaches zero as iterations continue.

\noindent\textbf{Key difference between our method and the low-rank solver for second-order RTE in \cite{bachmayr2024low}.} 
A finite element discretization for the 1D1V problem is considered in \cite{bachmayr2024low}. To simplify the notation, we discuss the finite difference counterpart of their preconditioner. The discrete linear system for the 1D1V problem is
\begin{equation}
    (-\BD_{\BOmega_x^2}\otimes\BD_{xx}^{\sigma_t}+\BI\otimes\BSigma_t)\bpsi^+=(\BI\otimes\BSigma_s)\bphi+\boldsymbol{1}\otimes G.
    \label{eq:1d1v-system}
\end{equation}

In this work, we employ the DSA framework to accelerate the SI, which naturally generalizes to high-dimensional problems with low-rank formats. In contrast, \cite{bachmayr2024low} adopts a different framework based on the exponential sums preconditioner \cite{beylkin2010approximation,bachmayr2023low}.
To enable the use of this preconditioner, they first perform Cholesky decompositions of the discrete spatial operator $\BD_{xx}^{\sigma_t}+\BSigma_t=\BT_{x}\BT_{x}^T$, thereby transforming the system \eqref{eq:1d1v-system} into the form of Kronecker sum by solving $\hat{\bpsi}^+=(\BI\otimes\BT_{x})^T\bpsi^+$:
\begin{equation}
   (-\BD_{\BOmega_x^2}\otimes\BD_{xx}^{\sigma_t}+\BI\otimes\BSigma_t)
   = (\BI\otimes \BT_{x})(-\BD_{\BOmega_{x}^2}\otimes \BI+\BI\otimes\hat{\BSigma}_t)(\BI\otimes \BT_{x})^T.
\end{equation}

However, when extending to higher dimensions, a new transformation is required due to the extra cross derivatives. Additionally, the bandwidth of the spatial operator will increase, scaling as $O(N_{\bx}^{\frac{d-1}{d}})$ 
in dimension $d$. Let $n$ denote the number of grid points along each spatial direction.  Then, unless further tensor decomposition is enabled with spatially separable $\sigma_t,\sigma_s$  and a tensor product gird, the memory cost of standard Cholesky decomposition increases to $O(N_{\bx}^{2-\frac{1}{d}})=O(n^{2d-1})$. In particular, the memory cost scales as $O(N_{\bx}^{\frac{3}{2}})=O(n^3)$ in 2D and $O(N_{\bx}^{\frac{5}{3}})=O(n^5)$ in 3D.

\begin{rem}
Another closely related work is \cite{peng2023sweep}, which proposes an implicit DLRA for the first order time-dependent RTE based on SI and transport sweep. However, it does not incorporate rank adaptivity or DSA for preconditioning. The core idea in \cite{peng2023sweep} is to partition the angular domain into a small number of octant families, noting that all angles within the same octant share the same upwind ordering of spatial elements on a Cartesian mesh and thus the same sweeping pattern. This fact allows using transport sweep to solve the equation in each SI. Meanwhile, during each sweep, a dense block matrix of size $O(r)\times O(r)$ must be inverted for every spatial element. As shown in our numerical tests, for challenging problems the effective rank of the solution tends to grow with increasing resolution, leading to progressively larger dense matrices. 

\end{rem}
\subsection{Inner loop solver\label{sec:lr-inner}}
To solve the discrete second-order equation \eqref{eq:lr-si-step}, we apply a preconditioned low-rank CG method \cite{kressner2011low}. The algorithm detail is presented in Alg. \ref{alg:lr-cg}. 

\begin{algorithm}[ht]
\caption{Low-rank conjugate gradient (CG) algorithm to solve $\BA\bx=\bb$.\label{alg:lr-cg}}
\begin{algorithmic}[1]
\STATE{\textbf{Input}: tolerance for the relative residual $\epsilon_{\textrm{CG}}$, the maximum number of inner loop CG iterations $N_{\textrm{Iter}}$, low-rank preconditioner $\BM^{-1}$, low-rank representations for the initial guess  $\bx^{(0)}$ and the right hand side $\bb$.}
\STATE{\textbf{Initialization}: set $k=0$. If $\|\bx^{(0)}\|<10^{-15}$, $\br^{(0)}=\bb$. Otherwise, compute the residual $\br^{(0)}=\textrm{\tt Truncated sum}(\bb-\BA\bx^{(0)})$.}
\STATE{Apply the low-rank preconditioner $\bp^{(0)}=\BM^{-1}\br^{(0)}$.}
\STATE{Compute the low-rank inner product $\rho^{(0)}=\textrm{\tt Inner Product}(\br^{(0)},\bp^{(0)})$.}
\FOR{$k=1:N_{\textrm{Iter}}$}
    \STATE{Apply the operator $\BA$, $\bq^{(k-1)}=\BA\bp^{(k-1)}$.}
    \STATE{$\alpha=\rho/\textrm{\tt Inner Product}(\bp^{(k-1)},\bq^{(k-1)})$.}
    \STATE{$\bx^{(k)}=\textrm{\tt Truncated Sum}(\br^{(k-1)}-\alpha\bq^{(k-1)})$.}
    \IF{$||\br^{(k)}||\leq \epsilon_{\textrm{CG}}$}
        \BREAK
    \ENDIF
    \STATE{Apply the preconditioner $\bz^{(k)}=\BM^{-1}\br^{(k)}$.}
    \STATE{Compute $\rho^{(k)}=\textrm{\tt Inner Product}(\bz^{(k)},\br^{(k)})$ and $\beta=\rho^{(k)}/\rho^{(k-1)}$.}
    \STATE{Update the searching direction $\bp^{(k)}=\textrm{\tt Truncated sum}(\bz^{(k)},\beta\bp^{(k)})$.}
\ENDFOR
\end{algorithmic}
\end{algorithm}

To accelerate the convergence, it is crucial to design an effective preconditioner. Here, we apply a preconditioning strategy similar to \cite{kressner2011low}. For parametric elliptic problems, \cite{kressner2011low} uses a uniform preconditioner for all parameters by solving the problem corresponding to a ``weighted average" parameter based on the spectrum of the parametric problem. Inspired by this strategy, we develop a \textit{two-level} Diffusion Preconditioner (DPC). Taking the matrix format as an example,  the preconditioner is defined as
\begin{equation}
    \BM=\BI\otimes\Big(-\frac{1}{3}(\BD_{xx}^{\sigma_t}+\BD_{yy}^{\sigma_t})+\BSigma_a\Big).
\end{equation}
Therefore,  $\BM^{-1}(\sum_{l=1}^r \BB_l\,\BX_l\otimes\BV_l)$ can be computed as
\begin{equation}
   \sum_{l=1}^r \BB_l\,\BV_l\otimes \left( \big(-\frac{1}{3}(\BD_{xx}^{\sigma_t}+\BD_{yy}^{\sigma_t})+\BSigma_a\big)^{-1}\BX_{l}\right).
\end{equation}
Then, we incorporate an additional level of approximation for $\left(-\frac{1}{3}(\BD_{xx}^{\sigma_t}+\BD_{yy}^{\sigma_t})+\BSigma_a\right)^{-1}$ through AMG. This two-level preconditioner is simple to apply and compatible with the low-rank formats used. 

At the conventional full-rank level, an efficient preconditioner can be constructed by applying directly a standard AMG to approximate $\left(-(\BOmega_{j,x}^2\BD^{\sigma_t}_{xx}+\BOmega_{j,x}\BOmega_{j,y}\BD^{\sigma_t}_{\textrm{cross}}+\BOmega_{j,y}^2\BD^{\sigma_t}_{yy})+\BSigma_t\right)^{-1}$ for each $\BOmega_j$. Nevertheless, this preconditioner does not possess the Kronecker product structure, and therefore cannot be directly extended to the low-rank case. As detailed in App.~\ref{apx:inner-pc}, the full-rank counterpart of our two-level DPC is less efficient than the standard AMG preconditioner, because DPC requires additional inner iterations at its second level to approximate $\left(-\tfrac{1}{3}(\BD_{xx}^{\sigma_t}+\BD_{yy}^{\sigma_t})+\BSigma_a\right)^{-1}$.
Hence, to offset the slower convergence caused by DPC, we employ an inexact inner loop solver with adaptively controlled error tolerance, thereby reducing computational cost.

Though low-rank geometric multigrid methods have been developed in \cite{ballani2013projection,hackbusch2015solution}, their extensions to angular discretizations is challenging, since the quadrature points in the angular space with differing accuracy orders lack a hierarchical structure. Meanwhile, angular and space-angle multigrid methods have been designed for the first order RTE in \cite{morel1991angular,gao2009fast,lee2012multigrid}, their implementation and extensions to the low-rank setting are non-trivial. 
In addition, for the reasons discussed in the comparison between our method and \cite{bachmayr2024low} (see Sec. \ref{sec:lr-SI-DSA}), we do not pursue an exponential sums preconditioner. Instead, this paper focuses on DPC, leaving the development of more efficient low-rank space-angle multigrid preconditioners to future work. 

\begin{rem}
A key observation made in App. \ref{apx:inner-pc} is that both our two-level DPC and the full-rank AMG preconditioner are quasi-optimal in the diffusion regime in the sense that the inner loop iterations required for convergence almost does not grow as spatial resolution grows. However, for multiscale or transport dominant problems, both of them requires increasing iterations to solve the inner loop problem \eqref{eq:lr-si-step}, as resolution improves. On the other hand, transport sweep for the first order system always converges with one iteration per angle. However, designing an efficient low-rank extension of transport sweep compatible with Kronecker product structure is non-trivial, when the effective rank of the underlying solution is not low or the mesh is not Cartesian.
\end{rem}

\begin{rem}\textbf{Inexact full-rank SI-DSA.} To our surprise, we numerically observe that our inexact low-rank SI-DSA algorithm is faster and more robust than its full-rank counterpart. We refer readers to  Sec. \ref{sec:diffusion} for more details. 
\end{rem}

\section{Numerical examples\label{sec:numerical}}
In this section, we demonstrate the performance of the proposed method through a series of benchmark tests in the 2D X-Y geometry, with particular emphasis on comparing our low-rank SI–DSA method against its full-rank counterpart.

In the simulations, we set the stopping criteria for both low-rank and full-rank SI-DSA methods as $\|\phi^{(k)}-\phi^{(k-1)}\|\leq 10^{-5}$. In addition, for our low-rank solver,  truncation tolerance in SVD and hierarchical HOSVD is $\epsilon_{\textrm{truncate}}=10^{-6}$. For the inner loop low-rank CG solver, the maximum number of iterations is capped at $200$, the initial relative residual is $\epsilon_{\textrm{CG}}^{(1)}=10^{-2}$, and the decay rate that controls the inner-loop accuracy in the inexact solver is $\gamma=0.1$. For the inner loop solver of full-rank SI-DSA, we adopt the following two approaches: (1) a full-rank CG solver preconditioned by AMG  with relative residual tolerance $10^{-12}$ and (2) a full-rank CG solver preconditioend by the full-rank counterpart of the two-level DPC, see Sec. \ref{sec:lr-inner} with relative residual tolerance $10^{-8}$. We denote the full-rank results obtained by AMG preconditioner as ``FR–AMG" and those obtained by two-level DPC as ``FR–DPC". In both the full-rank and low-rank solvers, the diffusion equation in the DSA step is solved using CG solver precondionted by AMG with a relative residual tolerance of $10^{-12}$. In all tests, we employ a mesh with $(N_x+1)\times (N_y+1)$ grid points in the physical space and the CL($N_{\theta},N_{\BOmega_z}$) quadrature rule for angular discretization. 

Our code is implemented in {\tt{Matlab}}. Specifically, we leverage {\tt{Hierarchical Tucker Toolbox}} \cite{kressner2012htucker} to implement our low-rank representations
and {\tt{iFEM}} package \cite{chen2009integrated} for AMG. All tests are performed on a Macbook pro laptop with M1 chip.

\noindent\textbf{Notation.} Before delving into numerical results, we introduce the following shorthand notation. 
To measure the accuracy of our low-rank methods, we use results given by FR-AMG as reference solutions, since it nearly solves the inner-loop problem exactly with an AMG solver. The error in the even-parity variable is defined as
\begin{equation}
    \BmcE_{\psi^+}^{\textrm{FR}} = \|\bpsi_{\textrm{LR}}^+ -\bpsi_{\textrm{reference}}^+\|_{l_2},
\end{equation}
where $\bpsi_{\textrm{LR}}^+$ denotes the solution on the full grid reshaped from the  low-rank solution. Similarly, we can define the error for the scalar flux $\BmcE_{\phi}^{\textrm{FR}}$. We consider two compression ratios for the low-rank method to assess its performance. The first is the compression ratio of the solution, defined as
\begin{equation}
\textrm{Solution compression ratio} = \frac{\textrm{DOFs of the low-rank solution}}{\textrm{DOFs of the full-rank solution}}.
\end{equation}
The second is the overall compression ratio for the $k$-th iteration of the iterative solver, which accounts for the memory reduction in both the Krylov vectors and residual:
\begin{equation}
    \parbox{4cm}{\centering The compression ratio \\ at the $k$-th iteration}  = \frac{\textrm{DOF}(\bx^{(k)}_{\textrm{LR}})+\textrm{DOF}(\bp^{(k)}_{\textrm{LR}})+\textrm{DOF}(\bq^{(k)}_{\textrm{LR}})+\textrm{DOF}(\bz^{(k)}_{\textrm{LR}})+\textrm{DOF}(\br^{(k)}_{\textrm{LR}})}
{\textrm{DOF}(\bx^{(k)}_{\textrm{LR}})+\textrm{DOF}(\bp^{(k)}_{\textrm{FR}})+\textrm{DOF}(\bq^{(k)}_{\textrm{FR}})+\textrm{DOF}(\bz^{(k)}_{\textrm{FR}})+\textrm{DOF}(\br^{(k)}_{\textrm{FR}})},    
\end{equation}
where $\bc^{(k)}_{\textrm{LR}}$ and $\bc^{(k)}_{\textrm{FR}}$ denote the low-rank and full-rank representations, respectively, of a generic vector  $\bc\in\{\bx,\bp,\bq,\br,\bz\}$. The definitions of $\bx,\bp,\bq,\br,\bz$  are provided in the low-rank CG algorithm (Alg. \ref{alg:lr-cg}).

\subsection{Diffusion dominant problem \label{sec:diffusion}}
\begin{figure}[]
   \begin{tikzpicture}
    \node[anchor=south west, inner sep=0] (image) at (0,0) {
     \includegraphics[trim=10mm 9mm 10mm 10.6mm, clip,width=0.32\textwidth]{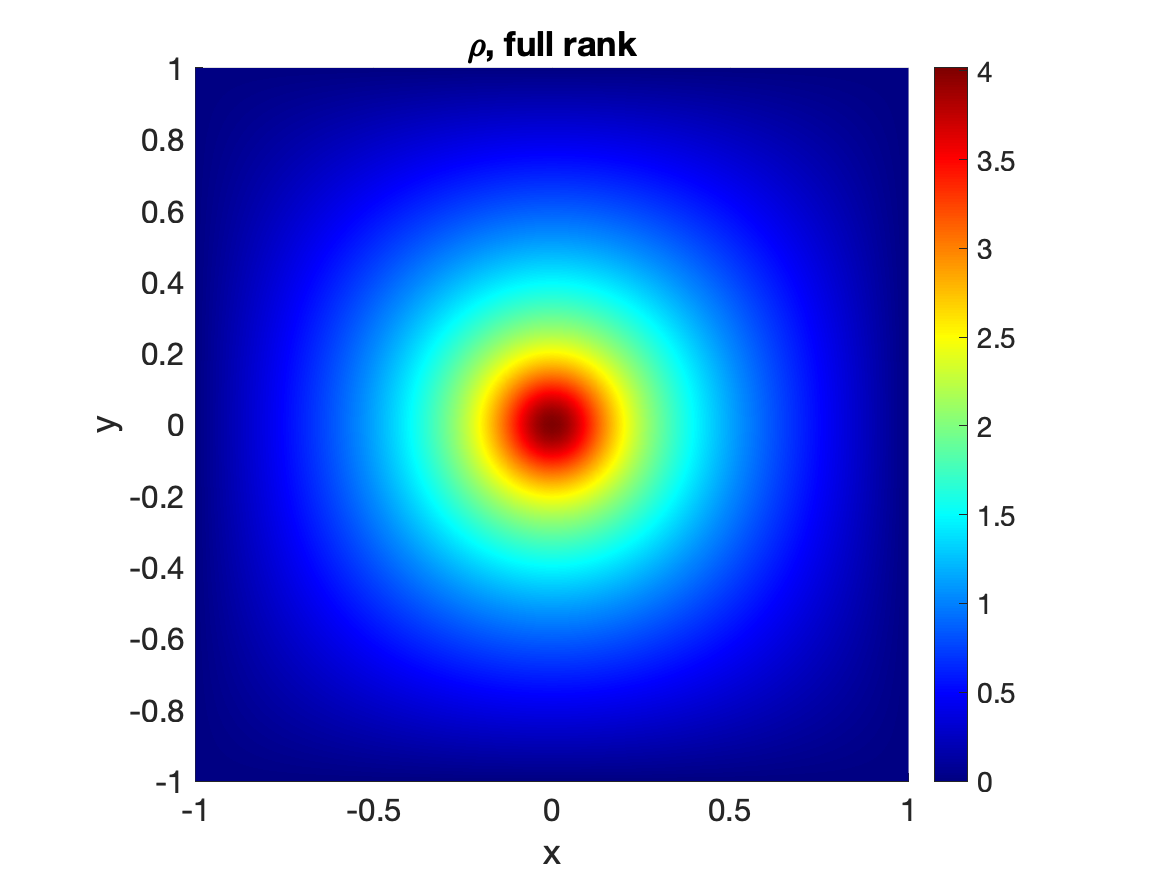}
    };
    \node[font=\bfseries, align=center, above=-1.2mm of image.north] {
      \tiny $\phi$, full-rank
    };
  \end{tikzpicture}
  \hspace{-2mm}
  \begin{tikzpicture}
    \node[anchor=south west, inner sep=0] (image) at (0,0) {
     \includegraphics[trim=10mm 9mm 10mm 10.6mm, clip,width=0.32\textwidth]{pictures/rho_diffusion_fr.eps}
    };
    \node[font=\bfseries, align=center, above=-1.2mm of image.north] {
      \tiny $\phi$, low-rank, tensor
    };
  \end{tikzpicture}
\includegraphics[width=0.32\textwidth]{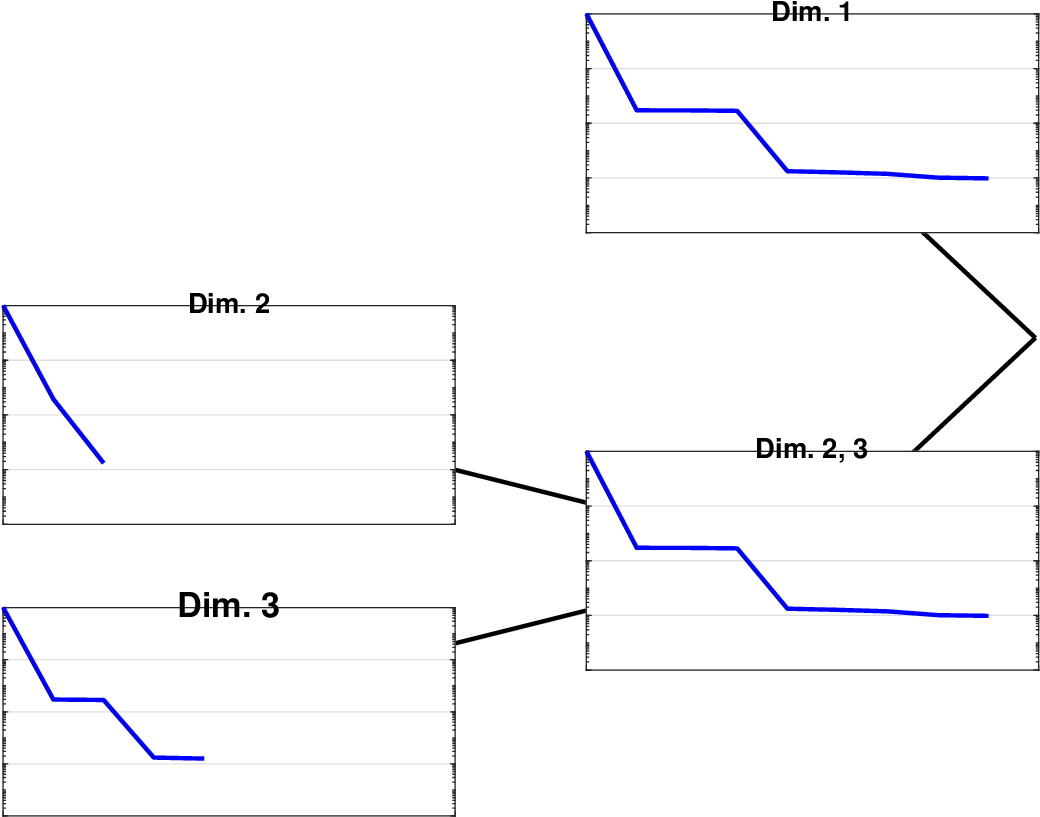}
  \caption{Full-rank solution, low-rank solution and the singular values for the matricization of the low-rank tensor for the diffusion dominant problem in Sec. \ref{sec:diffusion} with $(N_x,N_y,N_\theta,N_{\BOmega_z})=(96,96,60,30)$. \label{fig:diffusion-solution}}
\end{figure}

We first conduct a convergence study for a diffusion dominant problem on the computational domain $[-1,1]^2$ with vacuum boundary conditions, $\sigma_s(x,y)=100$, $\sigma_a(x,y)=0$ and a Gaussian source $G(x,y)=\exp(-100(x^2+y^2))$. At the $L$-th refinement level, we set $N_x=N_y=16 L$, $N_\theta=10L$ and $N_{\BOmega_z}=5L$, with $L=2,3,\dots,6$. Fig.~\ref{fig:diffusion-solution} presents a reference full-rank solution, along with the low-rank solutions in the HTT format and the hierarchical singular values from their matricizations.
As shown in Tab. \ref{tab:diffusion-results}, the $l_2$ difference in $\bphi$ and $\bpsi^+$ between the low-rank and the full-rank results are on the order of $O(10^{-4})$, validating the correctness of our low-rank solver.

\begin{figure}[]
  \begin{center} 
\includegraphics[width=0.45\textwidth]{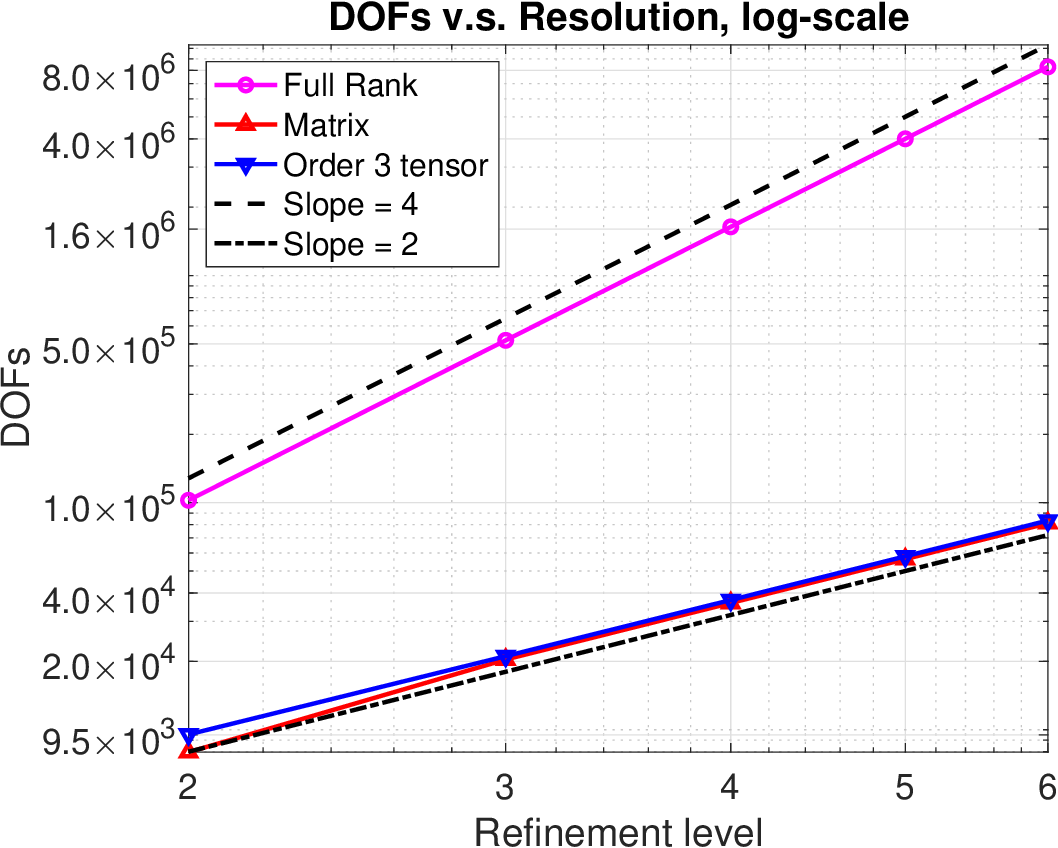}
\includegraphics[width=0.45\textwidth]{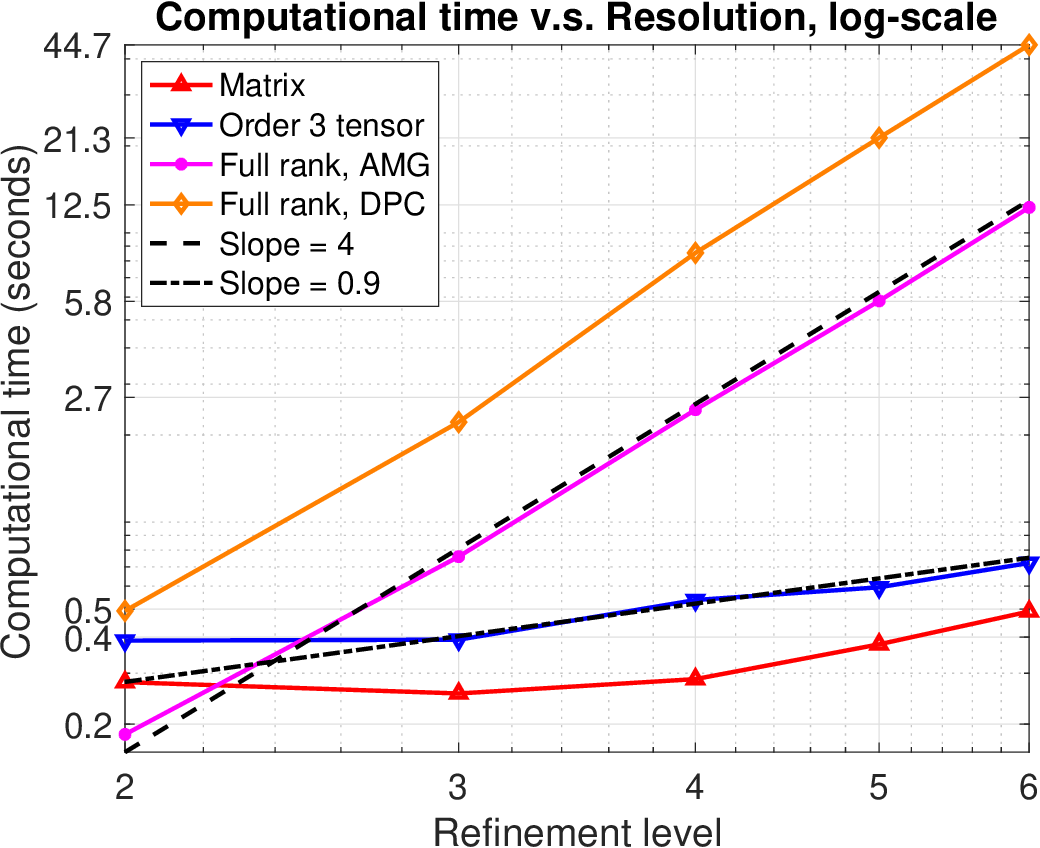}
  \caption{DOFs and wall clock computational time for the diffusion dominant problem in Sec. \ref{sec:diffusion} with different refinement level. \label{fig:diffusion-refinement-study}}
  \end{center}
\end{figure}

\noindent\textbf{Memory compression.} In the diffusion dominant regime, the angular flux exhibits a low-rank structure in the angular domain. As shown in Tab. \ref{tab:diffusion-results}, the ranks of the solutions obtained by our low-rank method using either matrix or HT decompositions remain the same as the mesh resolution grows, demonstrating that the proposed method can effectively capture the inherent low-rank solution structure. Furthermore,  in the left panel of Fig. \ref{fig:diffusion-refinement-study}, it is observed that the DOFs of the low-rank solutions scale linearly with DOFs in the physical space, i.e. $N_{\bx}=O(N_xN_y)\propto  O(L^2)$, 
whereas the DOFs of the full-rank counterpart scale as $O(N_xN_yN_\theta N_{\Omega_z})\propto O(L^4)$, indicating that the effective ranks of the low-rank solutions are largely independent of the angular resolution. Furthermore, for the case $(N_x,N_y,N_\theta,N_{\BOmega_z})=(96,96,60,30)$, our low-rank solver only requires approximately $1\%$ DOFs compared to the full-rank counterpart.

\noindent\textbf{Computational time.} As shown in the right panel of Fig. \ref{fig:diffusion-refinement-study}, though the inner-loop two-level DPC is less efficient than a full-rank AMG preconditioner, the proposed low-rank solver still achieves significant computational savings by effectively capturing the inherent low-rank structure. In fact, when $(N_x,N_y,N_\theta,N_{\BOmega_z})=(96,96,60,30)$, our low-rank solver with either the matrix or HTT decompositions achieves speedups of approximately $24.97\times$ and $16.97\times$ over FR-AMG, and $91.13\times$ and $61.94\times$ over FR-DPC, respectively. For this example, the matrix format is about $1.5\times$ faster than the HTT format, which we attribute to the additional overhead of the HTT representation.

\begin{figure}[]
  \begin{center} 
\includegraphics[width=0.45\textwidth]{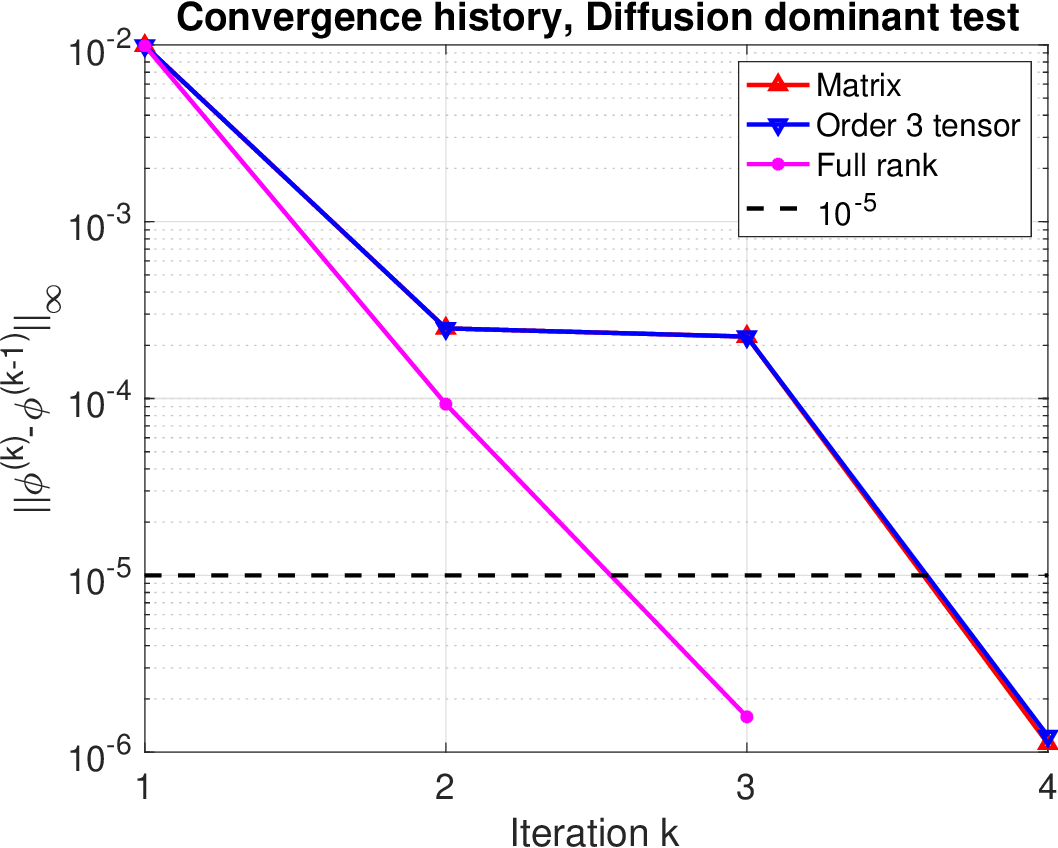}
\includegraphics[width=0.45\textwidth]{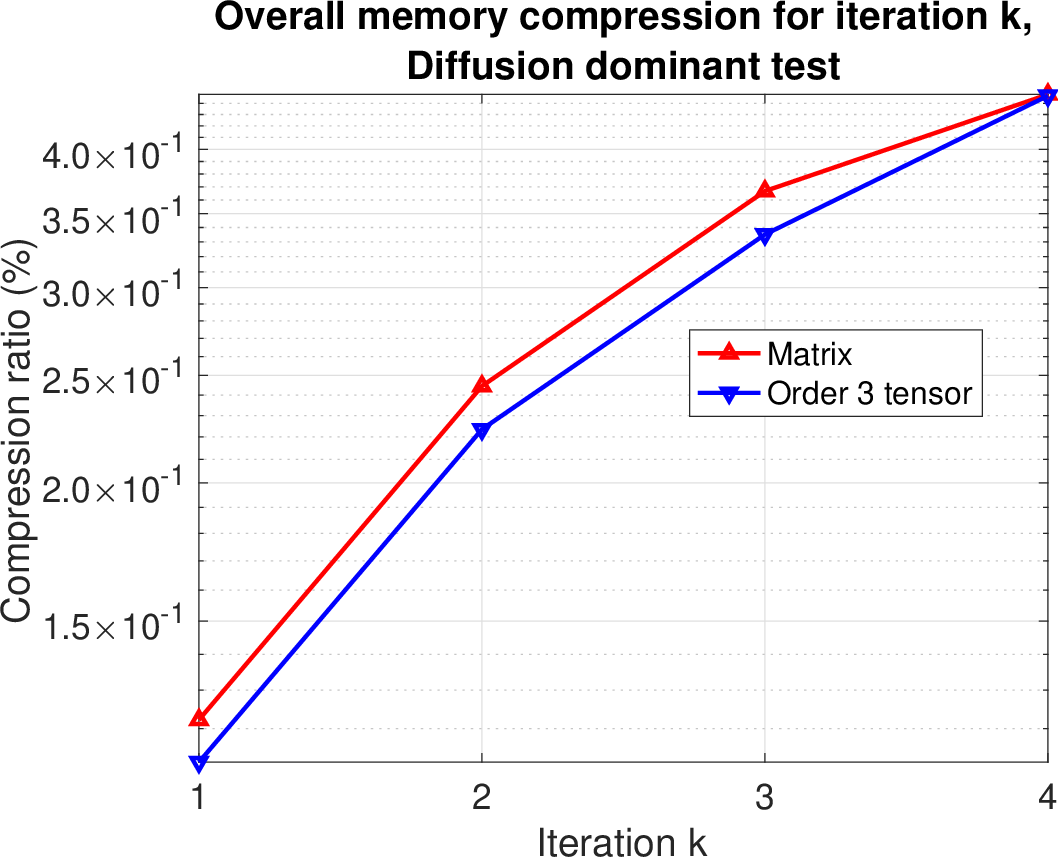}
  \caption{The history of $||\phi^{(k)}-\phi^{(k-1)}||$ and the memory compression ratio for the diffusion dominant problem in Sec. \ref{sec:diffusion} with $(N_x,N_y,N_\theta,N_{\BOmega_z})=(96,96,60,30)$. \label{fig:diffusion-convergence-history}}
  \end{center}
\end{figure}

\noindent\textbf{Convergence history.} In Fig. \ref{fig:diffusion-convergence-history}, we present the evolution of $\|\bphi^{(k)}-\bphi^{(k-1)}\|_{\infty}$ as the SI proceeds, along with the overall compression ratio at each step. In the final iteration, the ratio is observed to be slightly lower than that of the solution $\bpsi^{+}$, indicating that the Krylov tensors and the residual exhibit lower ranks than the solution itself.

\noindent\textbf{Comparison between the inexact low-rank and full-rank SI-DSA.} We further compare the inexact low-rank solver against its full-rank counterpart using the two-level DPC as the inner-loop preconditioner. For accuracy assessment, we benchmark the inexact SI results against the full-rank SI-DSA solution, obtained using the AMG inner-loop solver with a relative residual tolerance set to $10^{-12}$. As shown in Tab. \ref{tab:diffusion-comparison}, the inexact low-rank solver with $\epsilon_{CG}^{(1)}=10^{-2}$ and $\gamma=0.1$ yields smaller errors, ranging from $9\times10^{-5}$ to $1.43\times10^{-4}$, across all tested resolutions. Surprisingly, the errors for the inexact full-rank counterpart are much larger, ranging from $2.40\times10^{-1}$ to $2.22$ for $(\epsilon_{\mathrm{CG}}^{(1)},\gamma)=(10^{-2},0.1)$, and from $4.18\times10^{-4}$ to $6.48\times10^{-3}$ for $(\epsilon_{\mathrm{CG}}^{(1)},\gamma)=(10^{-5},0.01)$. Furthermore, the inexact low-rank solver is significantly faster than the full-rank version, underscoring both the efficiency and accuracy of the proposed low-rank approach for diffusion dominant problems.

In addition, similar observations were made across all numerical tests. For brevity, we omit the detailed results, highlighting here only the variable scattering test in Sec. \ref{sec:variable-scattering}. Even with a tight tolerance and fast decay rate, $(\epsilon_{\textrm{CG}}^{(1)},\gamma)=(10^{-5},0.01)$, the inexact full-rank SI-DSA yields an error of approximately $1.09\times10^{-2}$ for $(N_{\bx},N_{\by},N_{\theta},N_{\BOmega_z})=(96,96,60,30)$. In contrast, with $(\epsilon_{\textrm{CG}}^{(1)},\gamma)=(10^{-2},0.1)$,  inexact low-rank SI-DSA attains an error on the order of $10^{-5}$ on the same mesh while running about $2.91\times$ faster in matrix format and $2.91\times$ faster in the HTT format.

Surprisingly, the inexact low-rank SI-DSA is more robust and efficient than its full-rank counterpart. Understanding the theoretical basis for this interesting observation remains an open question and a promising direction for future research.

\begin{table}[htbp]
  \centering
 \medskip
 \begin{subtable}{\textwidth}
   \centering
    \begin{tabular}
    {|l|c|c|c|c|c|c|c|c|c|c|c|}
    \hline
    $N_x\times N_y\times N_{\theta}\times N_{\BOmega_z}$ &  $\BmcD_{\bpsi^+}$ & $\BmcD_{\bphi}$ &Solution compression & Rank of $\psi^+$ \\ \hline
$32\times32\times20\times 10$& $1.22\times10^{-4}$& $1.72\times 10^{-4}$ & $7.78\%$ & $7$\\
$48\times48\times30\times 15$& $1.59\times10^{-4}$&$2.25\times10^{-4}$&$3.93\%$&$8$\\
$64\times64\times40\times 20$&$1.81\times10^{-4}$ &$2.56\times10^{-4}$&$2.21\%$ & $8$ \\
$80\times80\times50\times 25$& $1.93\times10^{-4}$&$2.73\times10^{-4}$&$1.41\%$&$8$\\
$96\times96\times60\times 30$& $1.99\times10^{-4}$&$2.81\times10^{-4}$&$0.97\%$&$8$ \\ \hline
    \end{tabular}
      \caption{Matrix format. \label{tab:diffusion-matrix}}
\end{subtable}
\medskip
 \begin{subtable}{\textwidth}
   \centering
    \begin{tabular}
    {|l|c|c|c|c|c|c|c|c|c|c|c|}
    \hline
    $N_x\times N_y\times N_{\theta}\times N_{\BOmega_z}$ &  $\BmcD_{\bpsi^+}$ & $\BmcD_{\bphi}$ &Solution compression & Rank of $\psi^+$ \\ \hline
$32\times32\times20\times 10$& $9.47\times10^{-5}$& $1.34\times 10^{-4}$ & $9.31\%$ & $[9,3,5]$\\
$48\times48\times30\times 15$& $1.18\times10^{-4}$&$1.67\times10^{-4}$&$4.07\%$&$[9,3,5]$\\
$64\times64\times40\times 20$&$1.32\times10^{-4}$ &$1.87\times10^{-4}$&$2.27\%$ & $[9,3,5]$ \\
$80\times80\times50\times 25$& $1.40\times10^{-4}$&$1.98\times10^{-4}$&$1.45\%$&$[9,3,5]$\\
$96\times96\times60\times 30$& $1.43\times10^{-4}$&$2.02\times10^{-4}$&$1.01\%$&$[9,3,5]$ \\ \hline
    \end{tabular}
      \caption{Tensor format. \label{tab:diffusion-matrix}}
\end{subtable}
\caption{The difference between the low-rank solver  and the full-rank reference, the compression ratio and the rank for the diffusion dominant problem in Sec. \ref{sec:diffusion}.\label{tab:diffusion-results} }
\end{table}
\begin{table}[htbp]
\centering
\begin{subtable}{\textwidth}
\centering
\begin{tabular}
    {|l|c|c|c|c|c|c|c|c|c|c|c|}
    \hline
    $N_x\times N_y\times N_{\theta}\times N_{\BOmega_z}$ &  $||\bpsi_{\textrm{ILR}}^{+,\textrm{MAT}}-\bpsi_{\textrm{FR}}^{+}||_2$ & $||\bpsi_{\textrm{ILR}}^{+,\textrm{HTT}}-\bpsi_{\textrm{FR}}^{+}||_2$ &
    $||\bpsi_{\textrm{IFR},1}^{+}-\bpsi_{\textrm{FR}}^{+}||_2$ & $||\bpsi_{\textrm{IFR},2}^{+}-\bpsi_{\textrm{FR}}^{+}||_2$\\ \hline
$32\times32\times20\times 10$& $1.22\times10^{-4}$& $9.47\times 10^{-5}$ & $2.40\times10^{-1}$ & $4.18\times10^{-4}$\\
$48\times48\times30\times 15$& $1.59\times10^{-4}$&$1.18\times10^{-4}$&$5.55\times10^{-1}$&$1.11\times10^{-3}$\\
$64\times64\times40\times 20$&$1.32\times10^{-4}$ &$1.32\times10^{-4}$&$9.99\times10^{-1}$ & $2.11\times10^{-3}$ \\
$80\times80\times50\times 25$& $1.93\times10^{-4}$&$1.40\times10^{-4}$&$1.57\times10^{0}$&$3.58\times10^{-3}$\\
$96\times96\times60\times 30$& $1.99\times10^{-4}$&$1.43\times10^{-4}$&$2.22\times10^{0}$&$6.48\times10^{-3}$ \\ \hline
\end{tabular}
\caption{Accuracy.}
\end{subtable}
\medskip
\begin{subtable}{\textwidth}
\centering
\begin{tabular}
    {|l|c|c|c|c|c|c|c|c|c|c|c|}
    \hline
    $N_x\times N_y\times N_{\theta}\times N_{\BOmega_z}$ &  $T_{\textrm{ILR}}^{\textrm{MAT}}$ & $T_{\textrm{ILR}}^{\textrm{HTT}}$ &
    $T_{\textrm{IFR},1}$ & $T_{\textrm{IFR},2}$ & $T_{\textrm{FR}}$\\ \hline
$32\times32\times20\times 10$& $0.28$& $0.39$ & $0.11$ & $1.02$ & $0.18$\\
$48\times48\times30\times 15$& $0.26$&$0.39$&$0.44$&$1.65$ & $0.76$\\
$64\times64\times40\times 20$&$0.29$ &$0.54$&$1.52$ & $5.91$ & $2.45$\\
$80\times80\times50\times 25$& $0.38$&$0.59$&$3.72$&$15.08$ & $5.82$\\
$96\times96\times60\times 30$& $0.49$&$0.72$&$7.65$&$30.86$ & $12.25$ \\ \hline
\end{tabular}
\caption{Wall clock computational time (seconds).}
\end{subtable}
\caption{Comparison between the inexact low-rank and full-rank SI-DSA for the diffusion dominant problem in Sec. \ref{sec:diffusion}. ${\textrm{FR}}$: full rank results with AMG as inner loop solver and tolerance for the relative residual as $10^{-12}$.  $_{\textrm{ILR}}^{\textrm{MAT}}/ _{\textrm{ILR}}^{\textrm{HTT}}$:results for the inexact low-rank SI-DSA with $\epsilon_{\textrm{CG}}^{(1)}=10^{-2}$ and $\gamma=0.1$ in the matrix/HTT format. $\textrm{IFR},1$: inexact full-rank results with $\epsilon_{\textrm{CG}}^{(1)}=10^{-2}$ and $\gamma=0.1$. $\textrm{IFR},2$: inexact full-rank with results $\epsilon_{\textrm{CG}}^{(2)}=10^{-5}$ and $\gamma=0.01$. \label{tab:diffusion-comparison} }
\end{table}
\subsection{Transport dominant problem \label{sec:transport}}
\begin{figure}[]
   \begin{tikzpicture}
    \node[anchor=south west, inner sep=0] (image) at (0,0) {
     \includegraphics[trim=10mm 9mm 10mm 10.6mm, clip,width=0.32\textwidth]{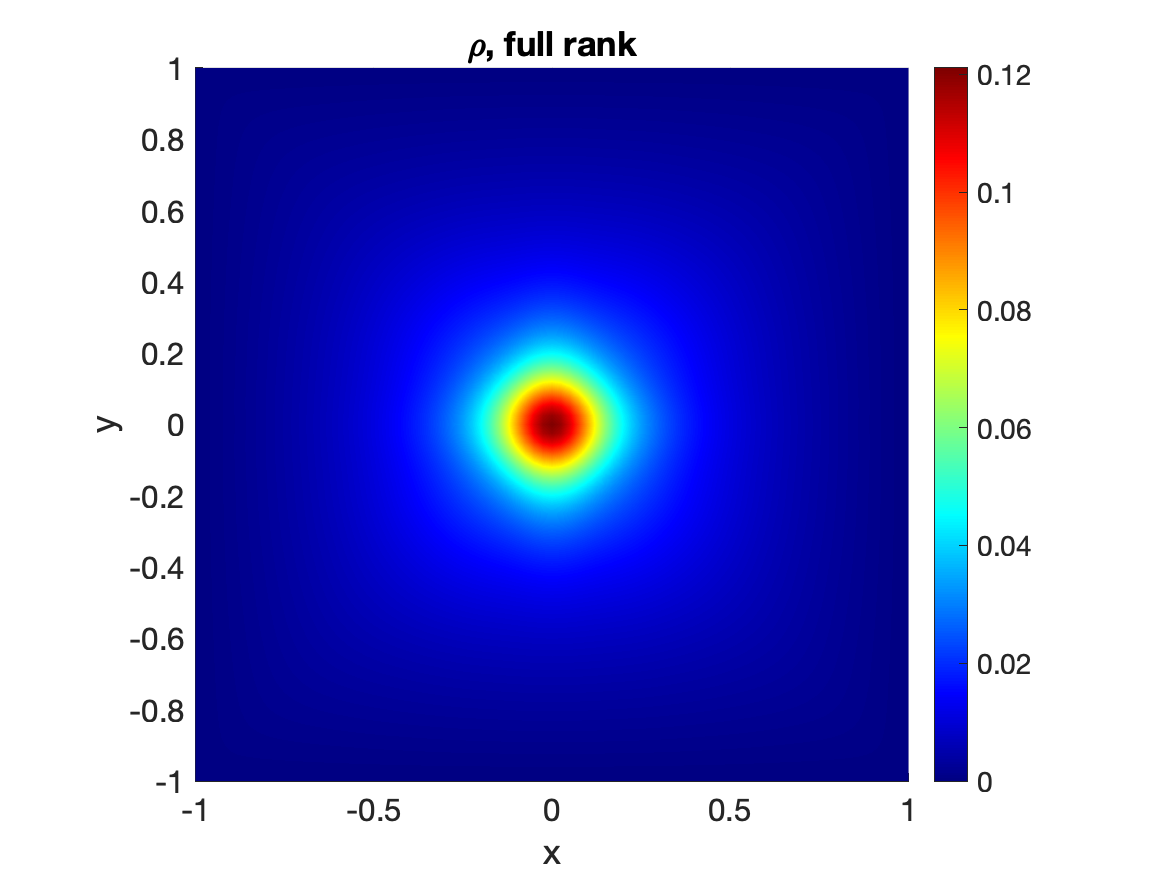}
    };
    \node[font=\bfseries, align=center, above=-1.2mm of image.north] {
      \tiny $\phi$, full-rank
    };
  \end{tikzpicture}
  \hspace{-2mm}
  \begin{tikzpicture}
    \node[anchor=south west, inner sep=0] (image) at (0,0) {
     \includegraphics[trim=10mm 9mm 10mm 10.6mm, clip,width=0.32\textwidth]{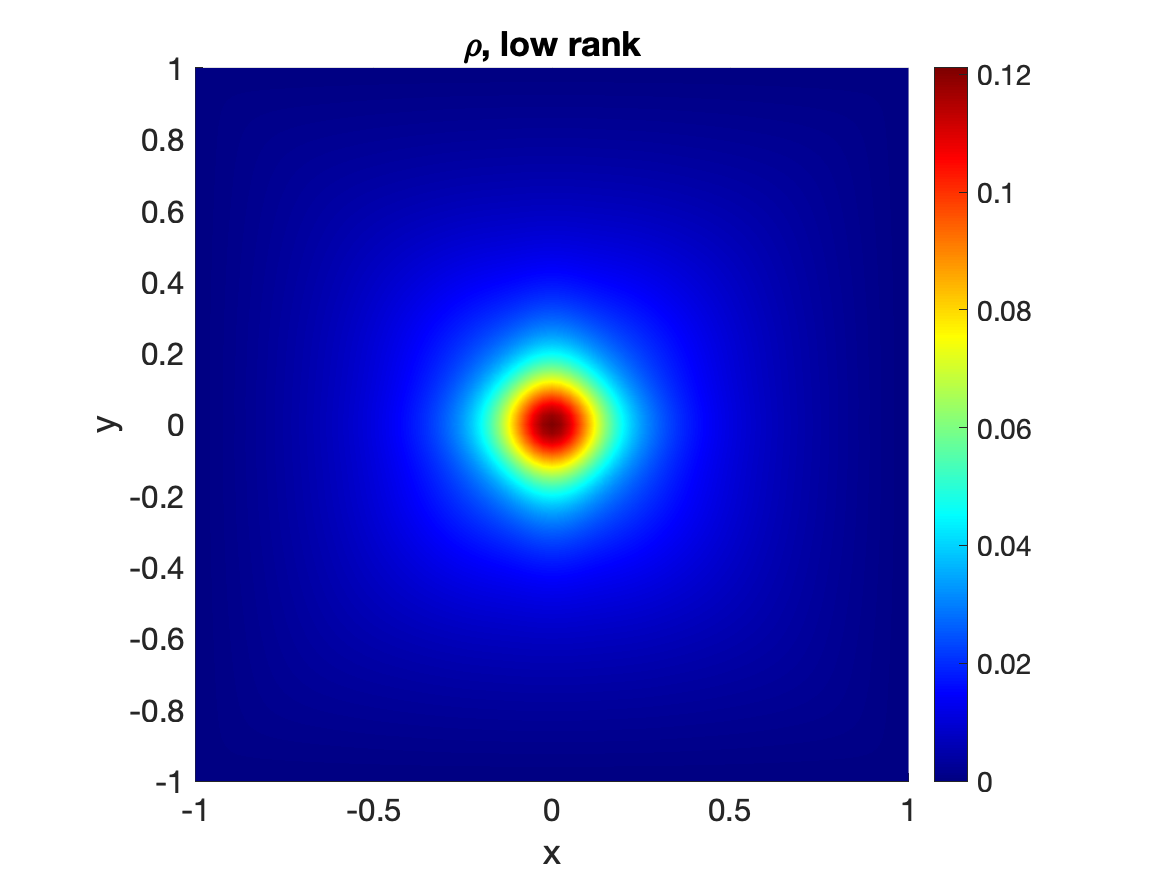}
    };
    \node[font=\bfseries, align=center, above=-1.2mm of image.north] {
      \tiny $\phi$, low-rank, tensor
    };
  \end{tikzpicture}
\includegraphics[width=0.32\textwidth]{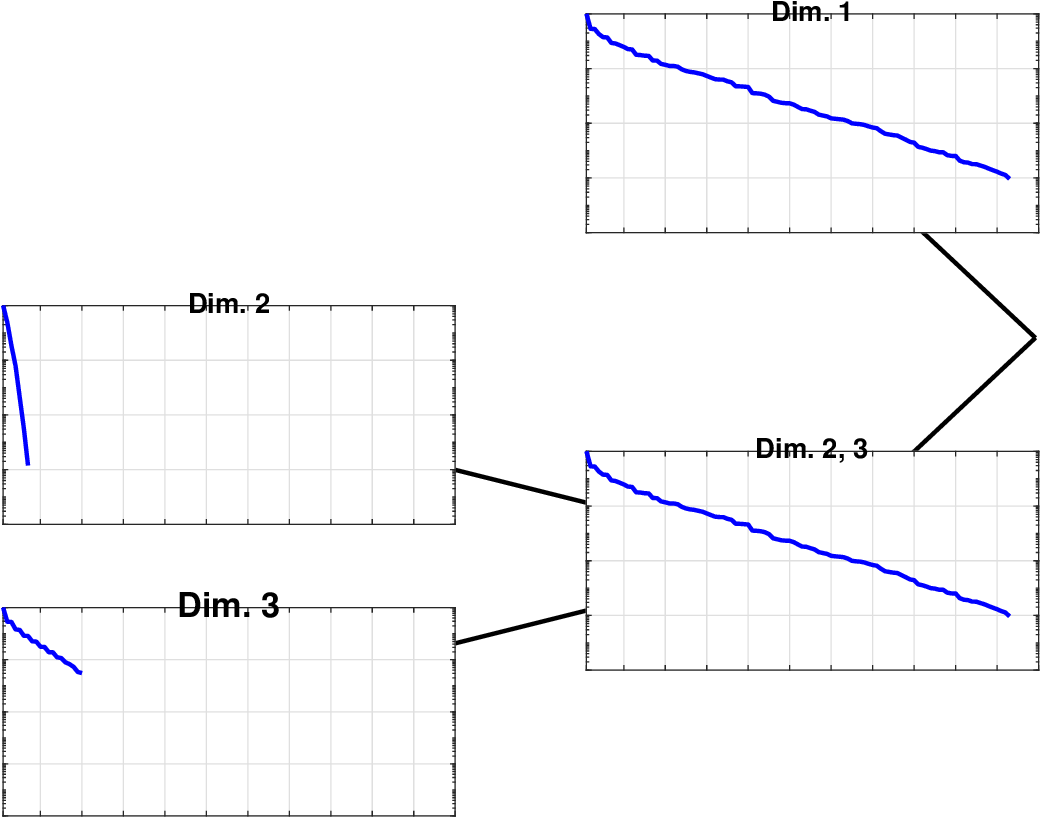}
  \caption{Full-rank solution, low-rank solution and the singular values for the matricization of the low-rank tensor for the transport dominant problem in Sec. \ref{sec:diffusion} with $(N_x,N_y,N_\theta,N_{\BOmega_z})=(80,80,40,20)$. \label{fig:transport-solution}}
\end{figure}

We consider a transport dominant problem on the computational domain $[-1,1]^2$ with vacuum boundary conditions, $\sigma_s(x,y)=1$, $\sigma_a(x,y)=0$ and a Gaussian source $G(x,y)=\exp(-100(x^2+y^2))$. In Fig. \ref{fig:transport-solution}, we plot the full-rank solution and the low-rank solution in the HTT format for the mesh $(N_x,N_y,N_\theta,N_{\BOmega_z})=(80,80,40,20)$ for comparison.  The singular values for the matricization of the low-rank solution
tensor are also shown.

We set $(N_x,N_y,N_\theta,N_{\BOmega_z}) = (16L,16L,8L,4L)$ with $L=2,3,4,5$. As shown in Tab. \ref{tab:transport-results}, on all meshes the difference between the low-rank solutions and the full-rank solutions with nearly exact inner-loop solves is approximately $10^{-8}$. 

\noindent\textbf{Rank of the solution and memory compression.} It is well known that the Kolmogorov $n$-width, i.e., the optimal error of approximating the solution manifold with a $n$-dimensional linear space, may decay slowly for transport-dominated problems \cite{ohlberger2015reduced,greif2019decay}. When the $n$-width exhibits such slow decay, the rank of the solution may become relatively high, thereby reducing the effectiveness of low-rank methods \cite{guo2022low}. 

For this transport-dominated problem, we observe that the solution rank is nearly proportional to $N_\theta$, in contrast to the previous diffusion dominant case where the rank remains almost constant. Nevertheless, even for this challenging problem, the low-rank method achieves a substantial reduction in memory usage.

\begin{figure}[]
  \begin{center} 
\includegraphics[width=0.45\textwidth]{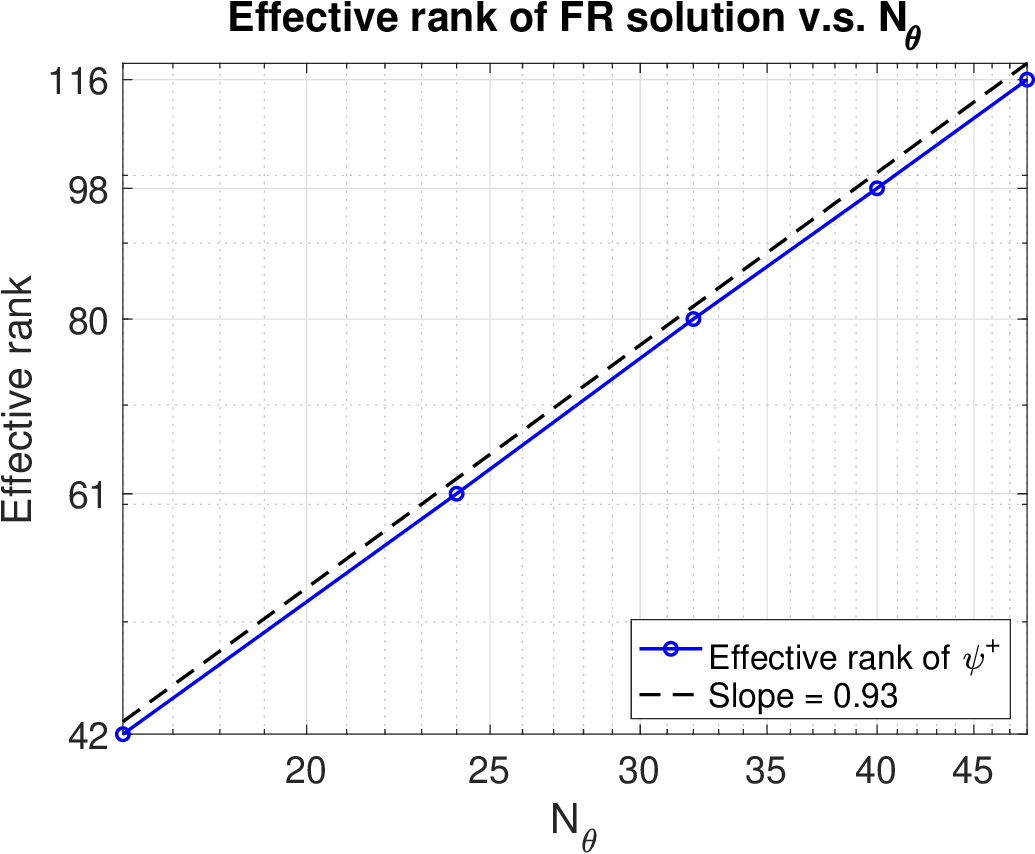}
\includegraphics[width=0.45\textwidth]{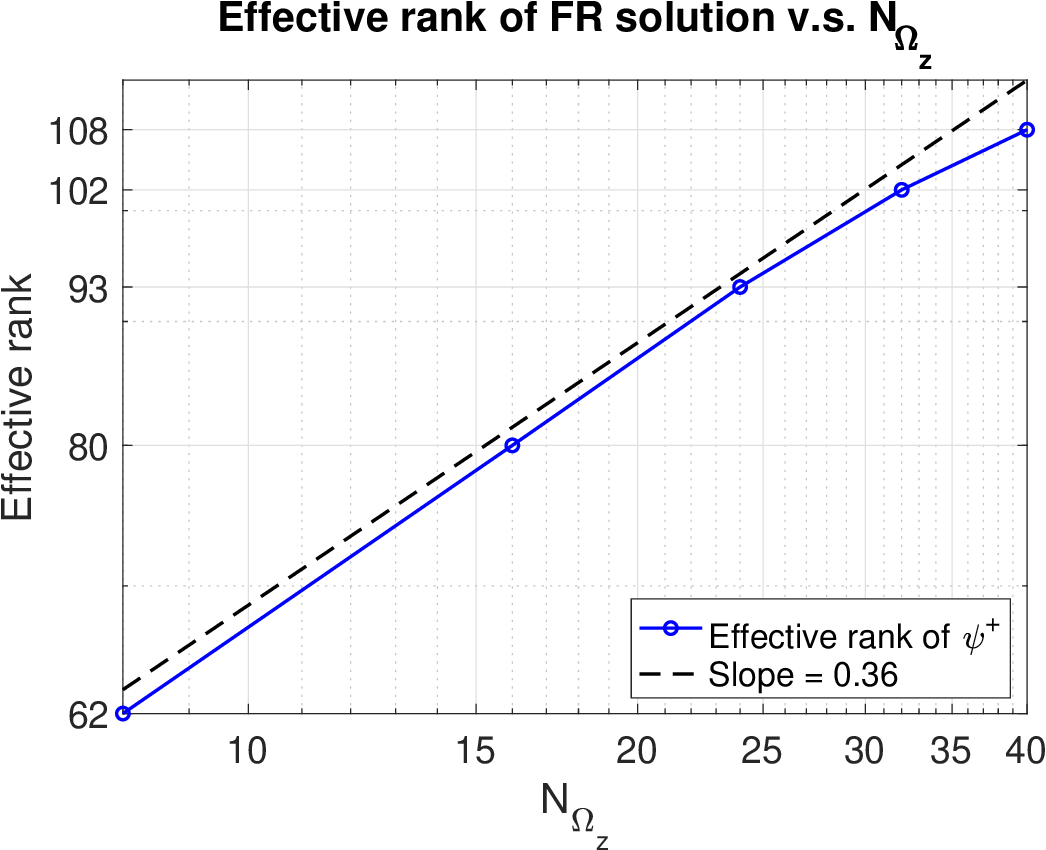}
  \caption{Effective rank of the full-rank solution in the matrix format for the transport dominant problem in Sec. \ref{sec:transport} with a fixed spatial resolution ($N_x=N_y=80$). FR: full-rank. Left: refinement study w.r.t $N_\theta$, $N_{\BOmega_z}=16$. Right: refinement study w.r.t $N_{\Omega_z}$,  $N_\theta=32$. \label{fig:transport-rank}}
  \end{center}
\end{figure}

\begin{figure}[]
  \begin{center} 
\includegraphics[width=0.32\textwidth]{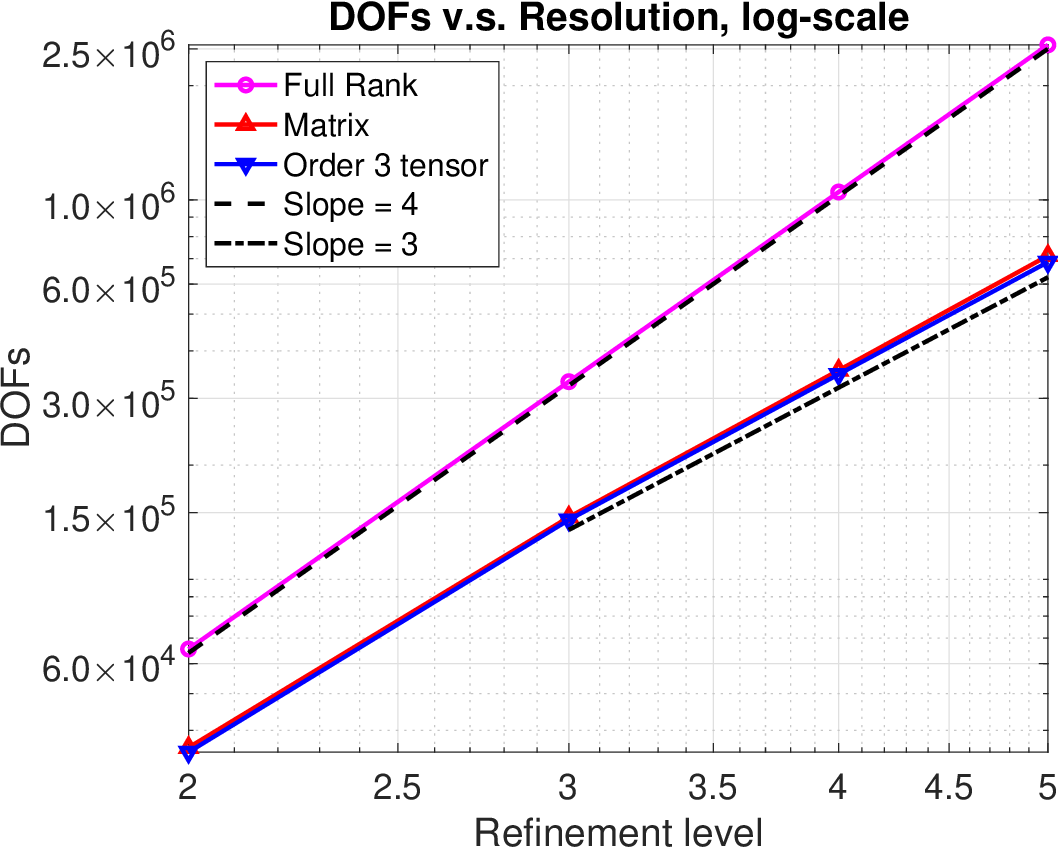}
\includegraphics[width=0.32\textwidth]{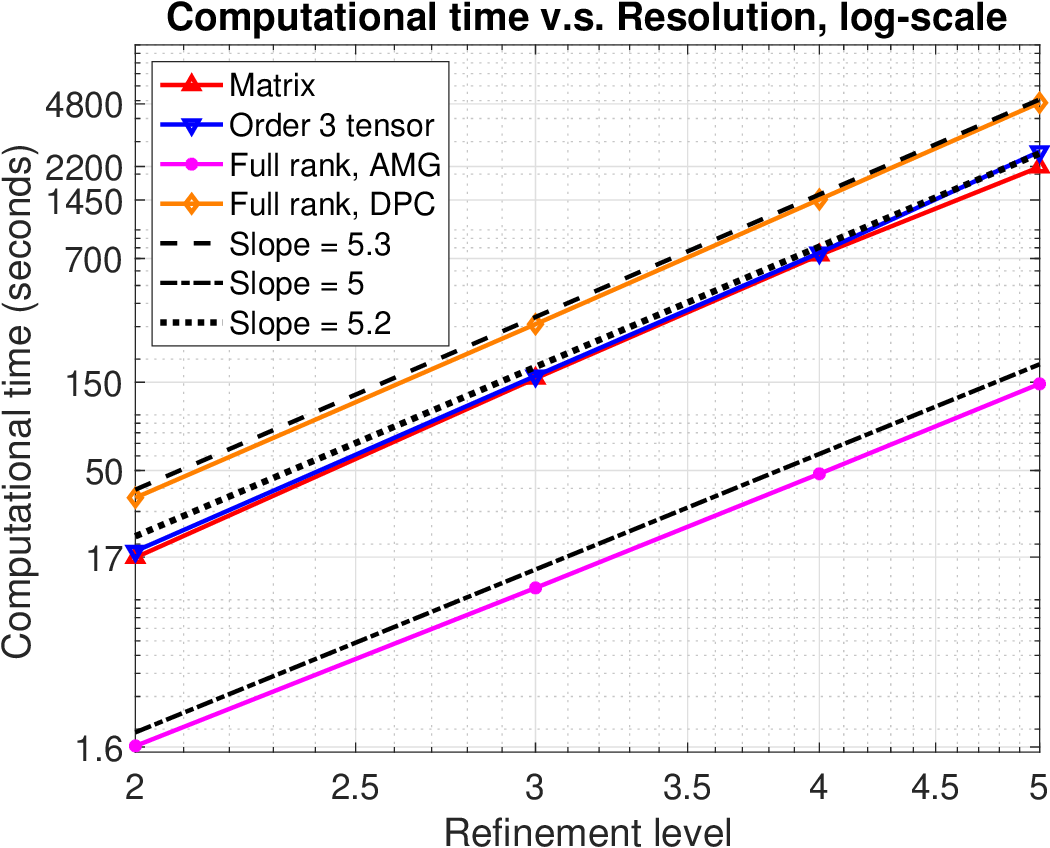}
\includegraphics[width=0.32\textwidth]{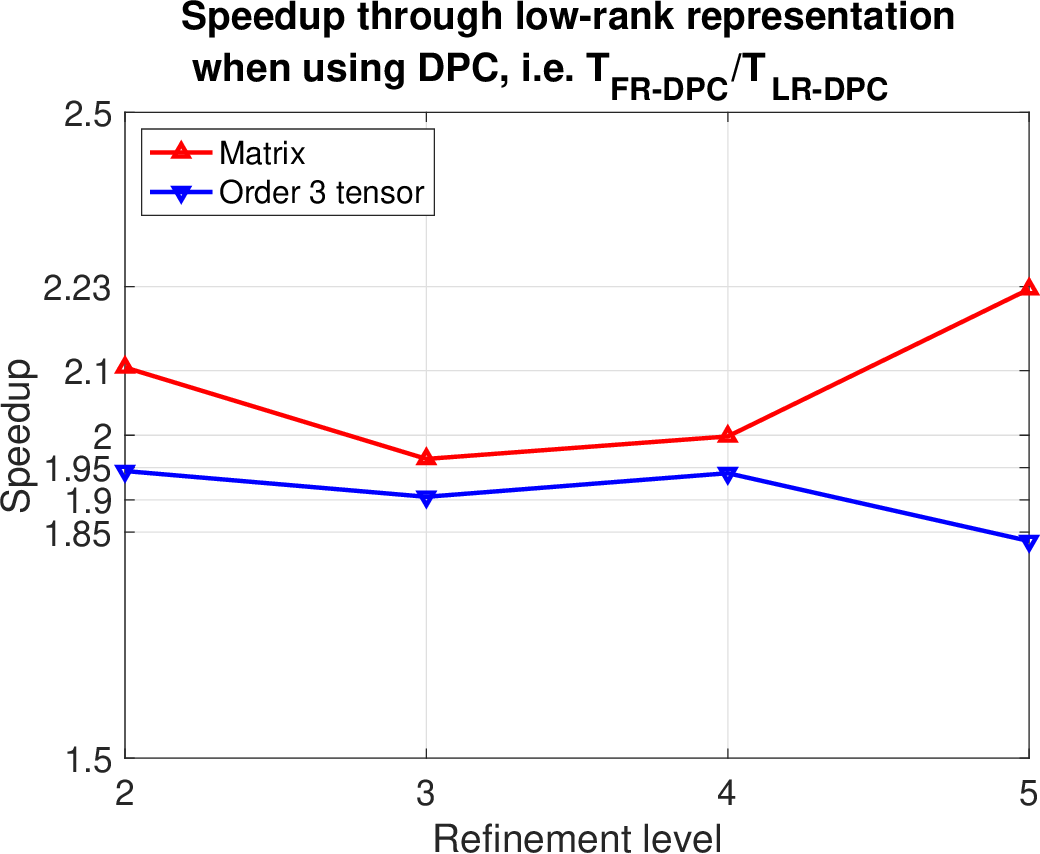}
  \caption{Resutls for the transport dominant problem in Sec. \ref{sec:transport} with different refinement level. Left: DOFs v.s. resolution. Middle: computational time v.s. resolution. Right: speedup achieved by the low-rank method when using the two-level DPC as the inner loop preconditioner. \label{fig:transport-refinement-study}}
  \end{center}
\end{figure}

\begin{enumerate}
\item \textbf{Effective rank of the full-rank solution.} We represent the full-rank solution in the matrix format, and then compute its effective rank through SVD with the same truncation strategy as the low-rank method. Specifically, we investigate the relation between the effective rank of the numerical solution and the spatial and angular resolution.


\begin{enumerate}
    \item \textbf{Effective rank versus spatial resolution.} Fixing the angular resolution at $(N_\theta,N_{\BOmega_z})=(32,16)$ and increasing only the spatial mesh resolution, we observe that the rank remains nearly constant with respect to spatial resolution (see Tab. \ref{tab:transport-rank}). In other words, the spatial resolution has negligible influence of the solution rank at the discrete level.
    
    \item \textbf{Effective rank versus the angular resolution.} 
We fix the spatial resolution at $N_x=N_y=80$.
With $N_{\BOmega_z}=16$ fixed, the effective rank of the solution, as shown in the left panel of Fig. \ref{fig:transport-rank}, exhibits sub-linear growth with respect to $N_\theta$, scaling approximately as $O(N_\theta^{0.93})$.
Conversely, with $N_\theta=32$ fixed, the effective rank, shown in the right panel of Fig. \ref{fig:transport-rank}, scales approximately as $O(N_{\BOmega_z}^{0.36})$.

\end{enumerate}
Hence, for this transport dominant problem, once the spatial resolution is sufficiently refined, its influence on the effective rank of the solution becomes negligible, whereas the angular resolution has a significant impact.

\item \textbf{Memory compression.} In the left panel of Fig. \ref{fig:transport-refinement-study}, the DOFs scale as $O(L^4)$ for the full-rank method and as $O(L^3)$ for the low-rank method. Consequently, as shown in Tab. \ref{tab:transport-results}, higher resolutions lead to greater compression for the low-rank methods. Specifically, for $(N_x,N_y,N_\theta,N_{\BOmega_z})=(80,80,40,20)$, the low-rank methods in the matrix and HTT formats reduce the DOFs by factors of approximately $3.59$ and $3.74$, respectively.

\item \textbf{The rank of the low-rank method versus the effective rank.} A comparison of the fourth columns in Tab. \ref{tab:transport-rank} and Tab. \ref{tab:transport-results} shows that the low-rank solutions generated by our low-rank solver qualitatively match the effective rank of the full-rank solution.
\end{enumerate}

\noindent\textbf{Computational time.} The wall-clock times for this example are shown in the right panel of Fig. \ref{fig:transport-refinement-study}. We make the following observations.
\begin{enumerate}
   \item When the same inner loop preconditioner (DPC) is applied in the low-rank and the full-rank methods, the low-rank method always achieves approximately a twofold speedup.
\item Compared with the full-rank method using AMG, the low-rank method with the two-level DPC is slower. We anticipate that the low-rank solver could match or exceed the performance of FR-AMG once a more efficient low-rank preconditioner, comparable to AMG but formulated in low-rank format, becomes available. This will be left for future work.
\end{enumerate}

\begin{figure}[]
  \begin{center} 
\includegraphics[width=0.45\textwidth]{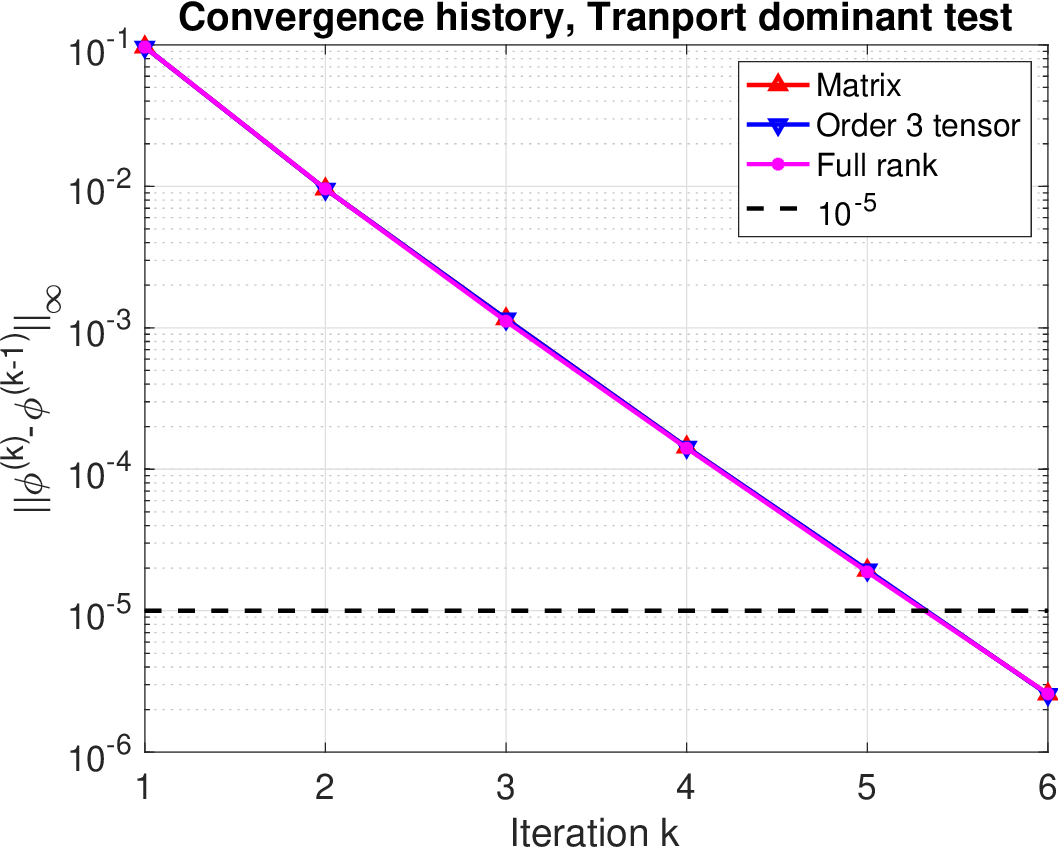}
\includegraphics[width=0.45\textwidth]{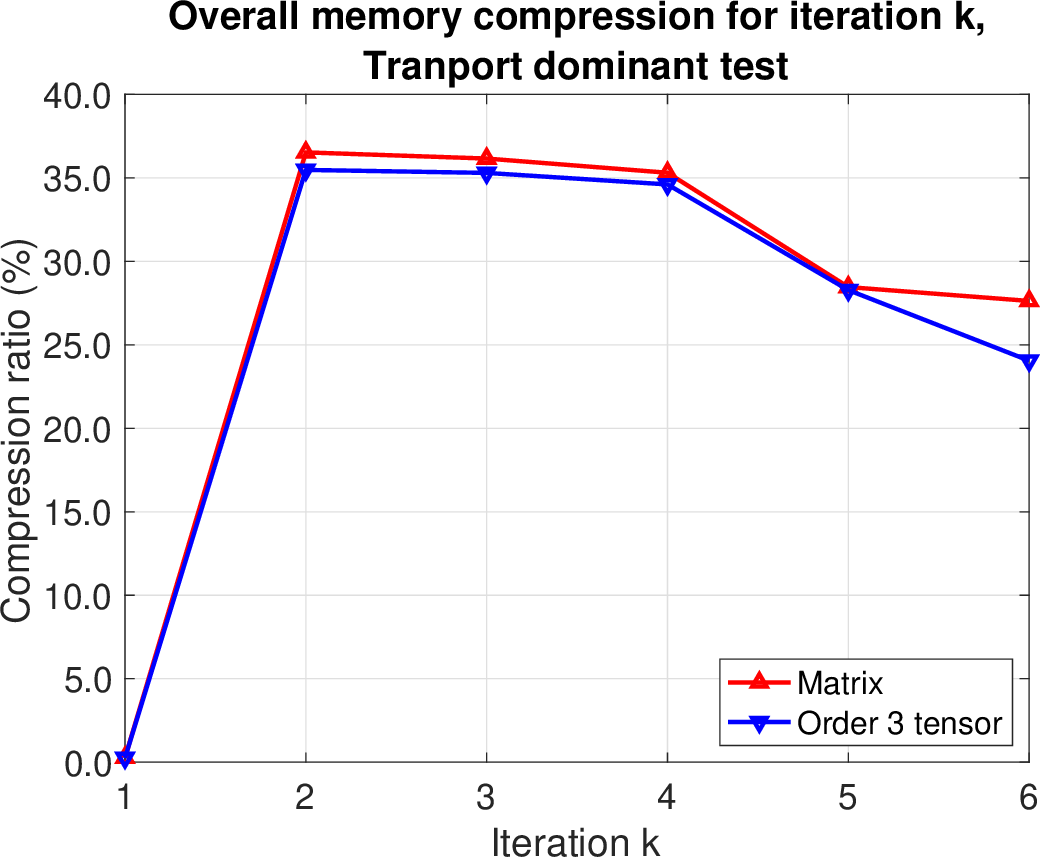}
  \caption{The history of $||\phi^{(k)}-\phi^{(k-1)}||$ and the memory compression ratio for the transport problem in Sec. \ref{sec:transport} with $(N_x,N_y,N_\theta,N_{\BOmega_z})=(80,80,40,20)$. \label{fig:transport-convergence-history}}
  \end{center}
\end{figure}

\noindent\textbf{Convergence history.} As shown in Fig. \ref{fig:transport-convergence-history}, the evolution of $||\bphi^{(k)}-\bphi^{(k-1)}||_{\infty}$ for the two low-rank methods closely matches that of the full-rank method.  Interestingly, the DOFs in the low-rank methods increase sharply after the first iteration as the Krylov update rapidly enriches the solution space, and then gradually decrease in subsequent iterations as the solver primarily refines existing components rather than introducing new directions.
We observe a similar convergence pattern for the following two problems.

\noindent\textbf{Comparison between the two low-rank formats.}  The matrix format is slightly faster owing to lower overhead in matrix computations, whereas the HTT format provides slightly greater compression. However, the overall differences in both computational time and required DOFs are minimal.

\textbf{In summary}, although the transport dominant problem exhibits much lower compressibility than the diffusion-dominated case, the proposed low-rank methods still deliver substantial memory savings. With the same preconditioner, they outperform the full-rank method in speed, and with further advances in low-rank preconditioning, their performance has the potential to match or even surpass that of the full-rank method equipped with an efficient AMG inner-loop preconditioner.

\begin{table}[htbp]
  \centering
 \medskip
    \begin{tabular}
    {|l|c|c|c|c|c|c|c|c|c|c|c|}
    \hline
               & $N_x=N_y=20$ & $N_x=N_y=40$ &  $N_x=N_y=60$ & $N_x=N_y=80$ & $N_x=N_y=81$\\ \hline
Effective rank & $77$ & $79$ & $80$ & $80$ & $81$ \\ \hline
    \end{tabular}
\caption{Effective rank of the full-rank solution in the matrix format for the transport dominant problem in Sec. \ref{sec:transport}  fixed angular resolution ($N_\theta=32$ and $N_{\BOmega_z}=16$).\label{tab:transport-rank} }
\end{table}
\begin{table}[htbp]
  \centering
 \medskip
 \begin{subtable}{\textwidth}
   \centering
    \begin{tabular}
    {|l|c|c|c|c|c|c|c|c|c|c|c|}
    \hline
    $N_x\times N_y\times N_{\theta}\times N_{\BOmega_z}$ &  $\BmcD_{\bpsi^+}$ & $\BmcD_{\bphi}$ &Solution compression & Rank of $\psi^+$ \\ \hline
$32\times32\times16\times 8$& $2.36\times10^{-8}$& $7.87\times 10^{-9}$ & $54.99\%$ & $32$\\
$48\times48\times24\times 12$& $4.02\times10^{-8}$&$2.41\times10^{-9}$&$43.96\%$&$58$\\
$64\times64\times32\times 16$&$4.18\times10^{-8}$ &$3.48\times10^{-9}$&$33.90\%$ & $80$ \\
$80\times80\times40\times 20$& $1.17\times10^{-7}$&$6.42\times10^{-9}$&$27.82\%$&$103$\\
\hline
    \end{tabular}
      \caption{Matrix format. \label{tab:diffusion-matrix}}
\end{subtable}
\medskip
 \begin{subtable}{\textwidth}
   \centering
        \begin{tabular}
    {|l|c|c|c|c|c|c|c|c|c|c|c|}
    \hline
    $N_x\times N_y\times N_{\theta}\times N_{\BOmega_z}$ &  $\BmcD_{\bpsi^+}$ & $\BmcD_{\bphi}$ &Solution compression & Rank of $\psi^+$ \\ \hline
$32\times32\times16\times 8$& $2.30\times10^{-8}$& $7.72\times 10^{-9}$ & $53.49\%$ & $[32,4,8]$\\
$48\times48\times24\times 12$& $3.45\times10^{-8}$&$2.15\times10^{-9}$&$43.24\%$&$[59,5,12]$\\
$64\times64\times32\times 16$&$3.27\times10^{-8}$ &$3.41\times10^{-9}$&$33.09\%$ & $[81,6,16]$ \\
$80\times80\times40\times 20$& $1.17\times10^{-7}$&$1.06\times10^{-8}$&$26.77\%$&$[103,7,20]$\\ \hline
    \end{tabular}
      \caption{HTT format. \label{tab:diffusion-matrix}}
\end{subtable}
\caption{The difference between the low-rank solver  and the full-rank solver, the compression ratio and the rank for the transport dominant problem in Sec. \ref{sec:transport}.\label{tab:transport-results} }
\end{table}
\subsection{Multiscale variable scattering \label{sec:variable-scattering}}
\begin{figure}[]
    \centering
    \includegraphics[trim=10mm 9mm 10mm 0mm, clip,width=0.45\textwidth]{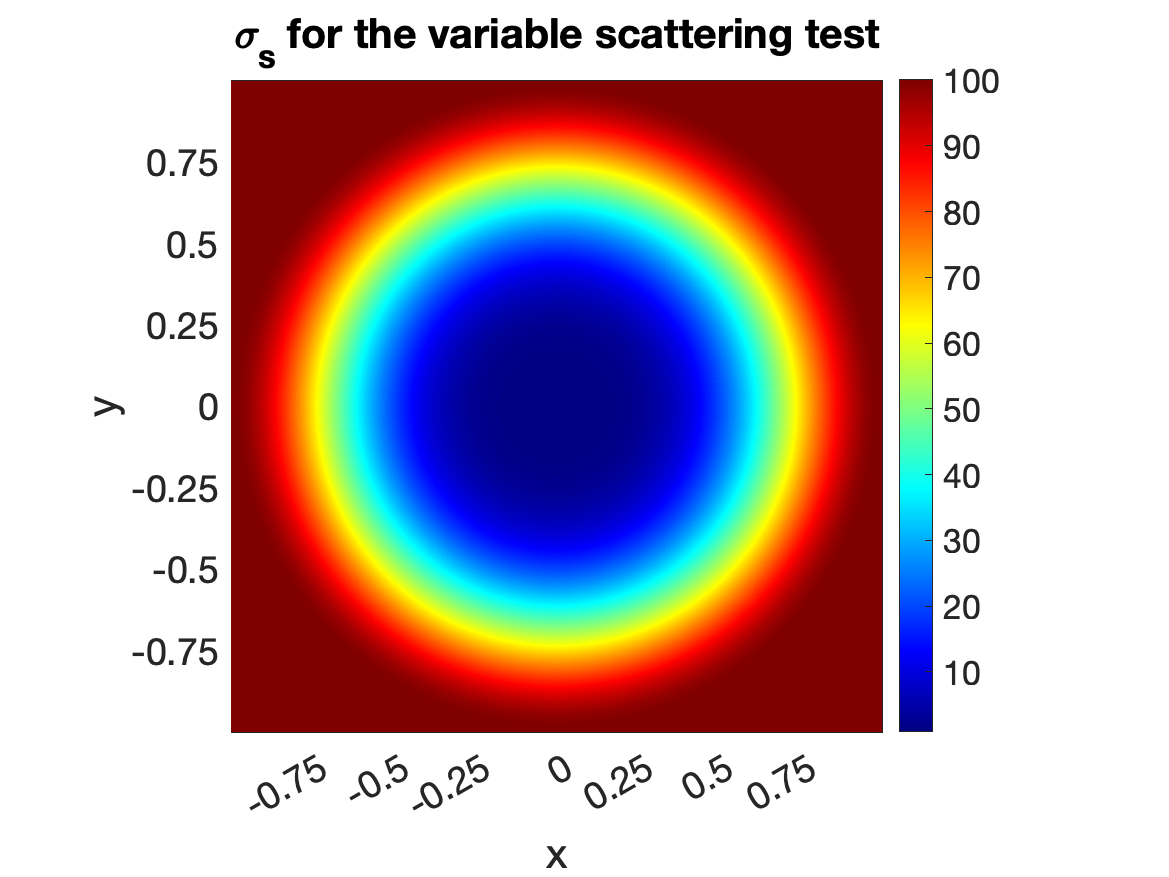}
    \hspace{-0.5cm}
    \includegraphics[width=0.45\textwidth]{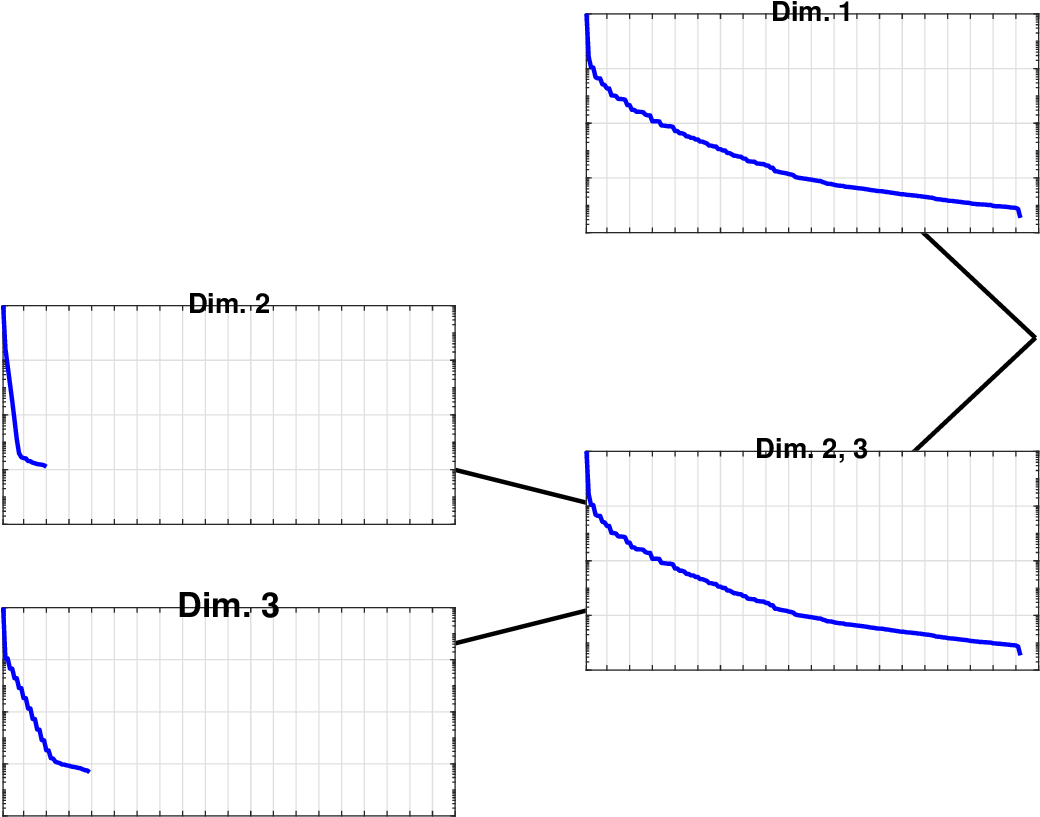}
   \begin{tikzpicture}
    \node[anchor=south west, inner sep=0] (image) at (0,0) {
     \includegraphics[trim=10mm 9mm 10mm 10mm, clip,width=0.45\textwidth]{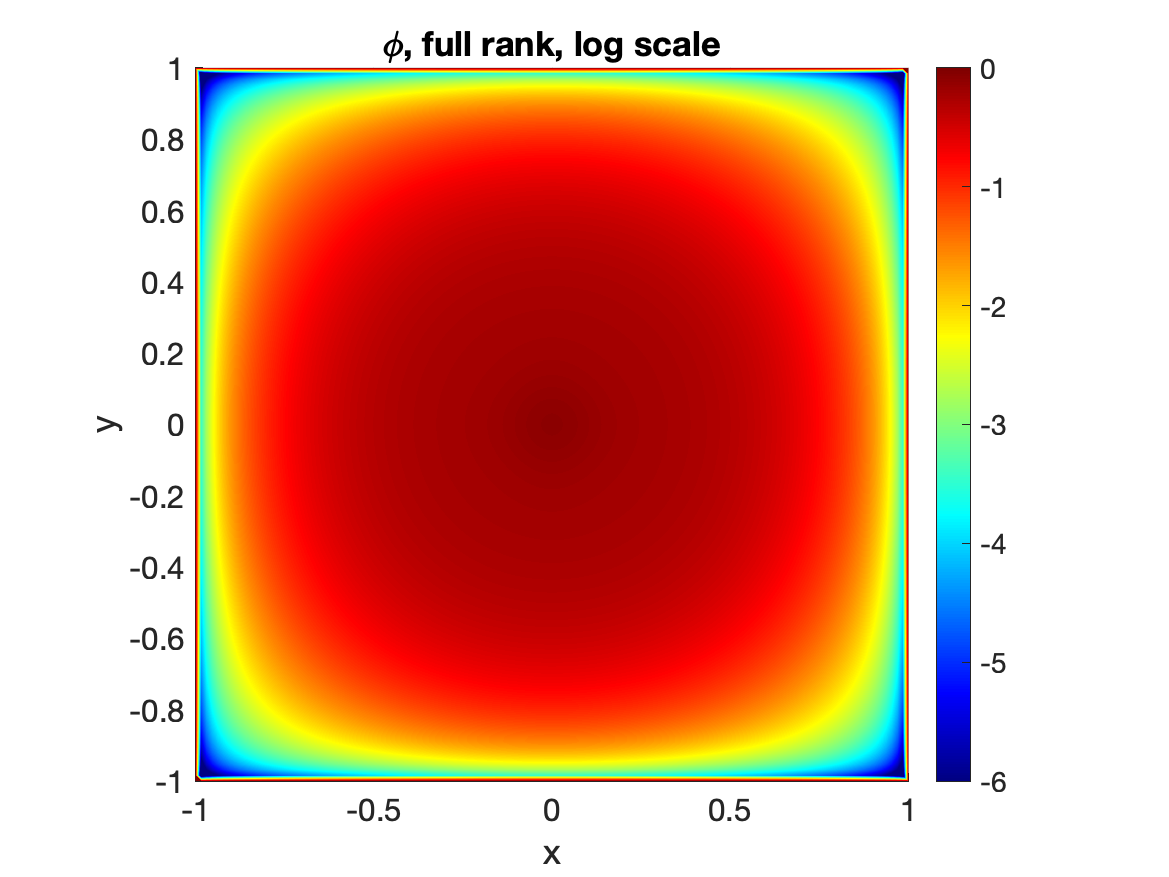}
    };
    \node[font=\bfseries, align=center, above=-1.2mm of image.north] {
      \tiny $\phi$, full-rank, log-scale
    };
  \end{tikzpicture}
  \hspace{-2mm}
  \begin{tikzpicture}
    \node[anchor=south west, inner sep=0] (image) at (0,0) {
     \includegraphics[trim=10mm 9mm 10mm 10mm, clip,width=0.45\textwidth]{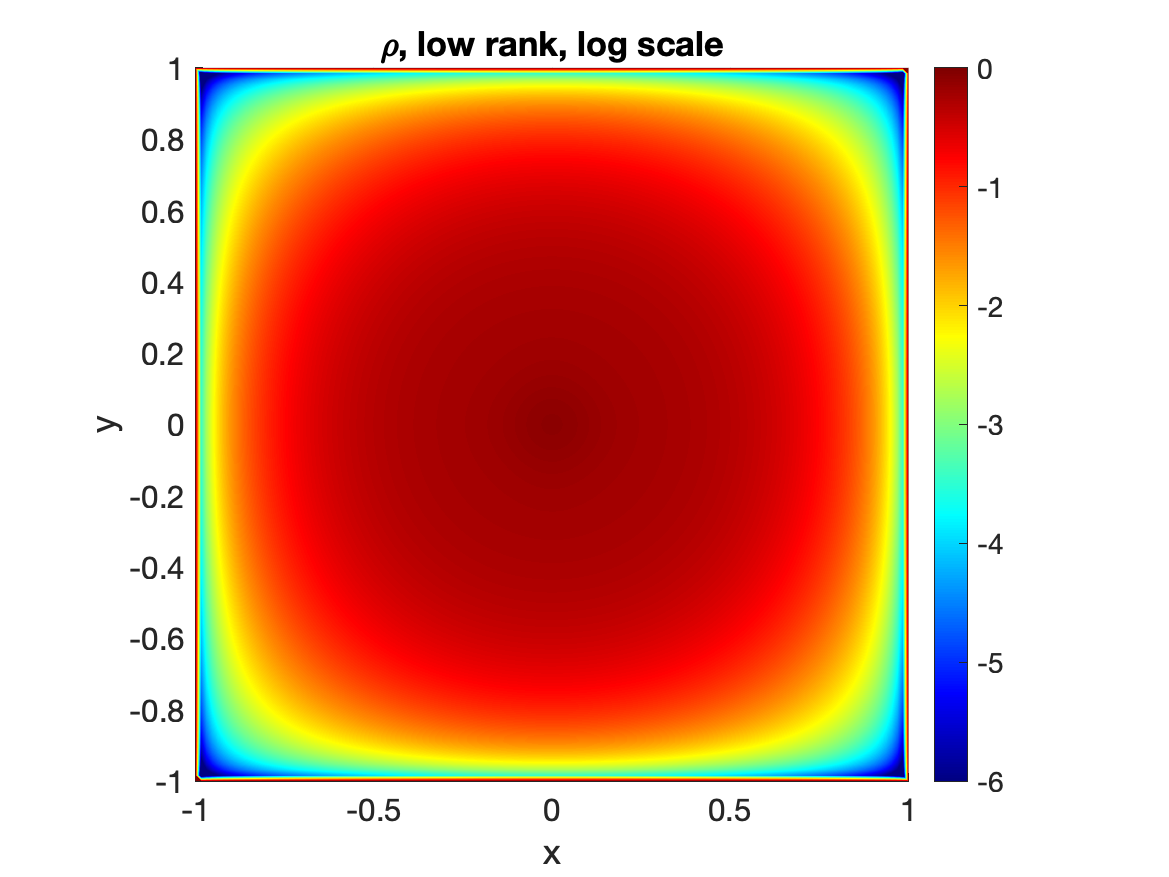}
    };
    \node[font=\bfseries, align=center, above=-1.2mm of image.north] {
      \tiny $\phi$, low-rank, tensor, log-scale
    };
  \end{tikzpicture}
  \caption{Variable scattering problem in Sec. \ref{sec:variable-scattering}, $(N_x,N_y,N_\theta,N_{\BOmega_z})=(80,80,40,20)$. Top left: set-up for the scattering cross section. Top right: singular values for the matricization of the low-rank tensor. Bottom left: full-rank solution. Bottom right: low-rank solution in the order 3 tensor format. \label{fig:vs-solution}}
\end{figure}

We consider a problem on the computational domain $[-1,1]^2$ with vacuum boundary conditions, no absorption, a Gaussian source $G(x,y)=\exp(-100(x^2+y^2))$ and a variable scattering defined as 
\begin{equation}
    \sigma_s(x,y) = \begin{cases}
                    99r^4(r^2-2)^2+1;, \quad\text{if}\; r=\sqrt{x^2+y^2}\leq 1,\\
                    100,\quad \text{otherwise}.
                    \end{cases}
\end{equation}
The configuration of $\sigma_s(x,y)$ can be found in Fig. \ref{fig:vs-solution}.
From the center of the computational domain to the boundary, the scattering strength smoothly grows from $1$ to $100$. In other words, the problem smoothly transitions from transport dominant to diffusion dominant, leading to strong multiscale effects. In the simulation, we set $N_x=N_y=16L$, $N_\theta=10L$ and $N_{\BOmega_z}=5$ with  
 $L=2,3,\dots,8$.
 
In Fig. \ref{fig:vs-solution}, we present the full-rank and low-rank solutions for $L=8$. Across all meshes, the difference between the two is on the order of $10^{-5}$ in the $l_2$ norm (see Tab. \ref{tab:vs-results}).

\begin{figure}[]
  \begin{center} 
\includegraphics[width=0.32\textwidth]{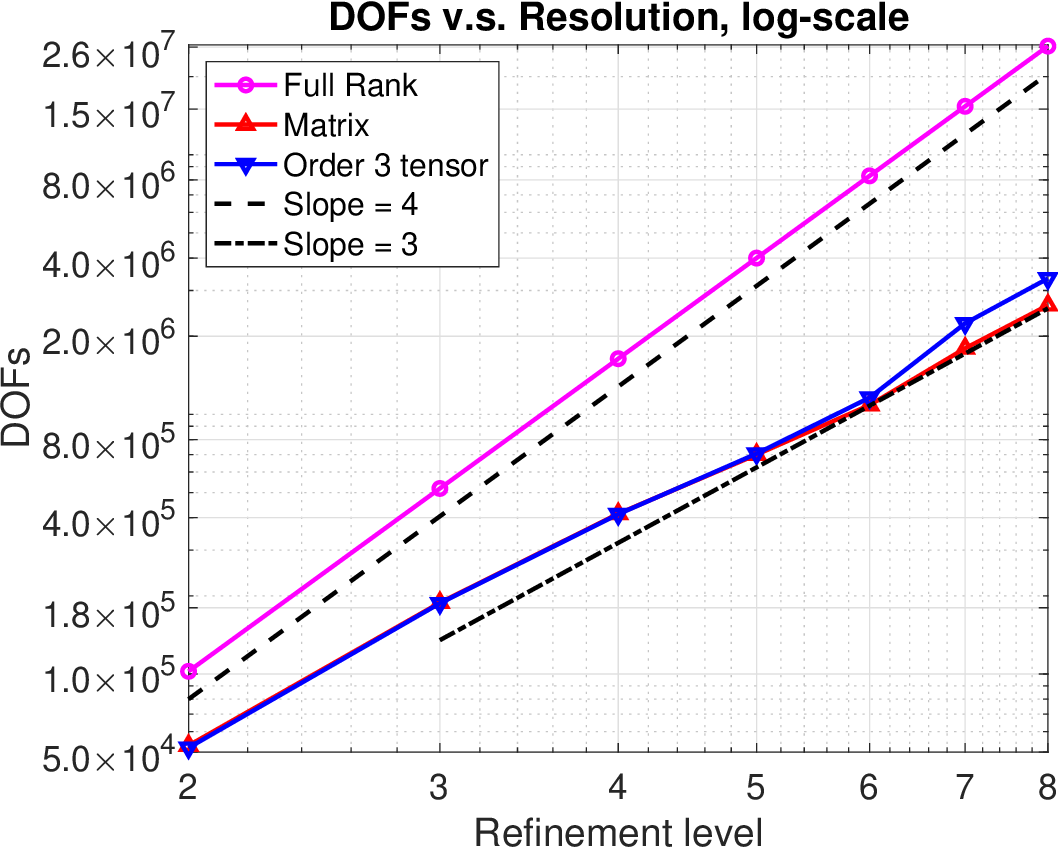}
\includegraphics[width=0.32\textwidth]{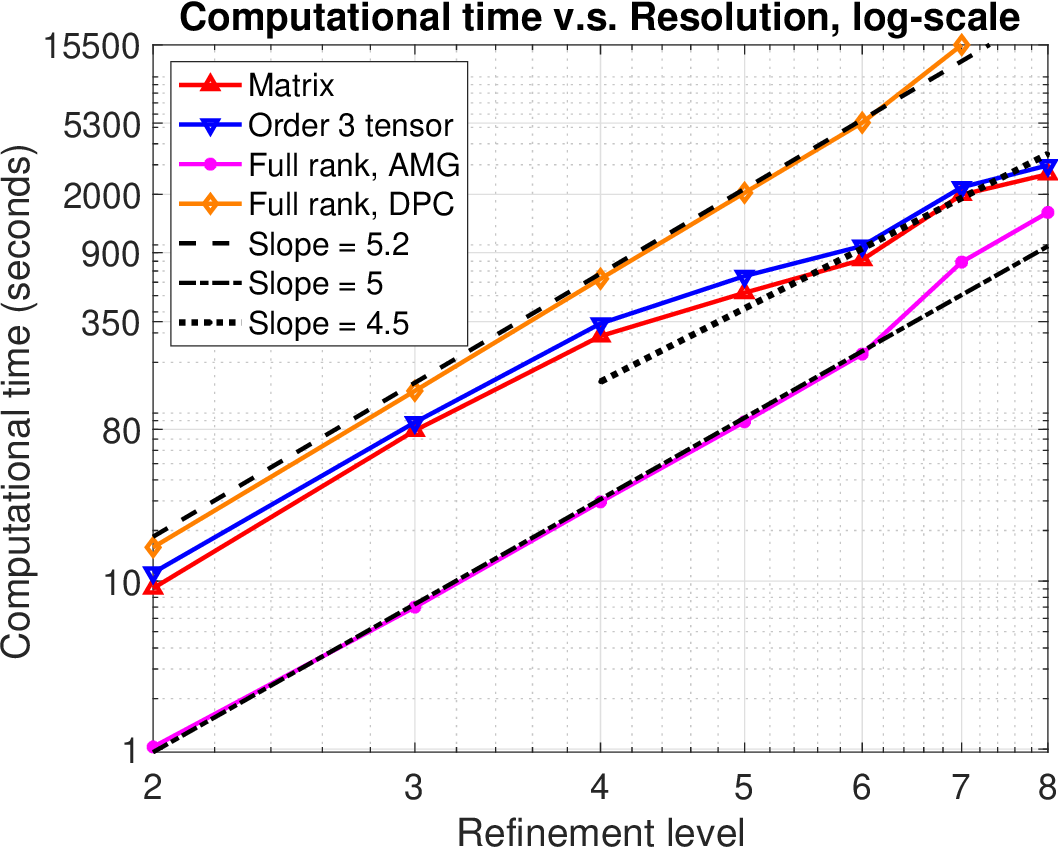}
\includegraphics[width=0.32\textwidth]{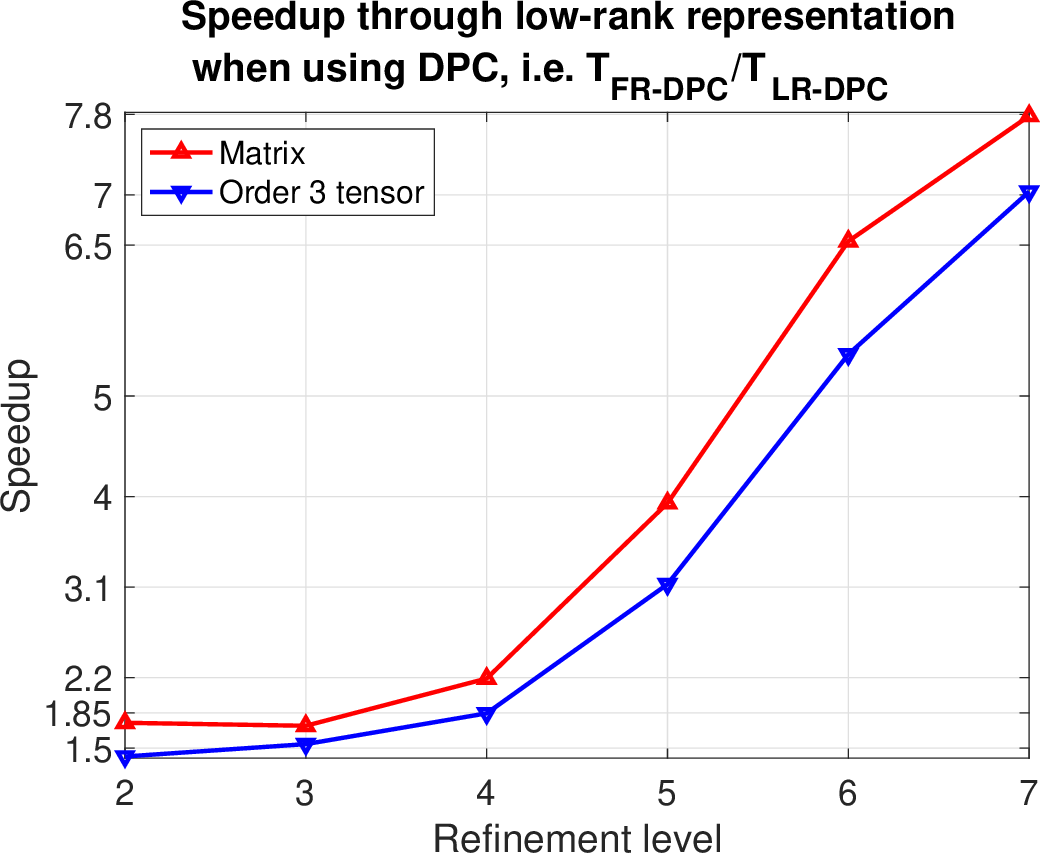}
  \caption{Results for the variable scattering problem in Sec. \ref{sec:variable-scattering} with different refinement level. Left: DOFs v.s. resolution. Middle: wall clock time v.s. resolution. Right: speedup obtained by the low-rank method when using two-level DPC.  \label{fig:vs-refinement-study}}
  \end{center}
\end{figure}

\noindent\textbf{Memory compression.} This variable scattering problem  exhibits multiscale structures and contains a large transport-dominated region. As demonstrated in Tab. \ref{tab:vs-results}, the rank of the low-rank solution grows as the resolution grows. As shown in right plot of Fig. \ref{fig:vs-refinement-study}, DOFs of the low-rank methods scale approximately as $O(L^3)$, while DOFs of the full-rank method scales as $O(L^4)$. As the resolution grows, more significant computational saving is achieved by the low-rank method. Specifically, when $(N_x,N_y,N_\theta,N_{\BOmega_z})=(128,128,80,40)$, low-rank methods reduces the DOFs by a factor of approximately $9.96$ with the matrix format and $7.86$ with the tensor format.

\noindent\textbf{Computational time.} As presented in the middle panel of Fig. \ref{fig:vs-refinement-study}, the computational time roughly scales as $O(L^5)$ for the full-rank methods and $O(L^{4.5})$ for the low-rank method, where $L$ is the refinement level. As discussed in App.~\ref{apx:inner-pc}, the standalone AMG preconditioner is much more efficient than the two-level DPC. As a result, the full-rank method with AMG is faster than the low-rank method with DPC. Nevertheless, because the computational time of the low-rank method grows more slowly with resolution, we anticipate that, as the resolution increases, the low-rank method with DPC will eventually outperform the full-rank method with AMG. Furthermore, when the two-level DPC is employed as the inner-loop preconditioner, the low-rank solver already achieves a speedup of $1.5\times$ to $7.8\times$ over the full-rank method, with larger efficiency gains observed at higher resolutions (see the right panel of Fig.\ref{fig:vs-refinement-study}).

\begin{figure}[]
  \begin{center} 
\includegraphics[width=0.45\textwidth]{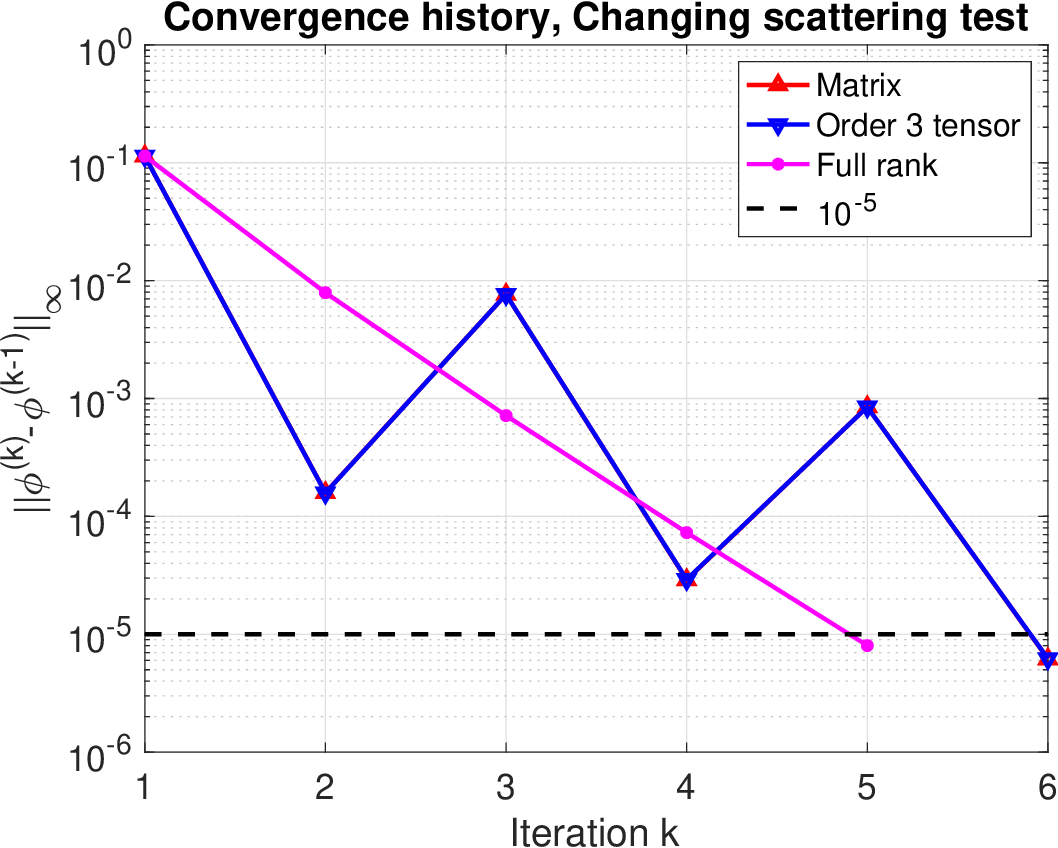}
\includegraphics[width=0.45\textwidth]{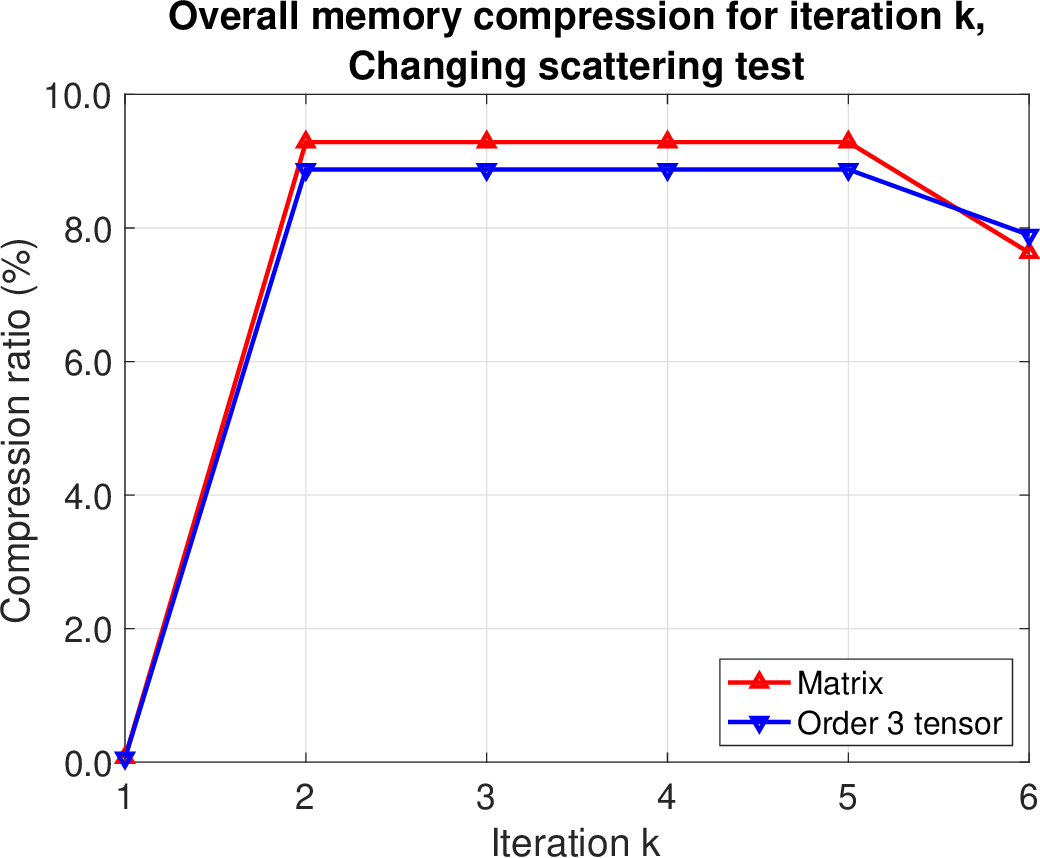}
  \caption{The history of $||\phi^{(k)}-\phi^{(k-1)}||$ and the memory compression ratio for the variable scattering problem in Sec. \ref{sec:variable-scattering} with $(N_x,N_y,N_\theta,N_{\BOmega_z})=(128,128,80,40)$. \label{fig:vs-convergence-history}}
  \end{center}
\end{figure}

\noindent\textbf{Convergence history.} As shown in Fig. \ref{fig:vs-convergence-history}, the inexact low-rank methods exhibit oscillatory behavior in $||\bphi^{(k)}-\bphi^{(k-1)}||$. The DOFs increase sharply after the first iteration, then reach a plateau as the iterations proceed, and decrease slightly in the final iteration.
\begin{table}[htbp]
  \centering
 \medskip
 \begin{subtable}{\textwidth}
   \centering
    \begin{tabular}
    {|l|c|c|c|c|c|c|c|c|c|c|c|}
    \hline
    $N_x\times N_y\times N_{\theta}\times N_{\BOmega_z}$ &  $\BmcD_{\bpsi^+}$ & $\BmcD_{\bphi}$ &Solution compression & Rank of $\psi^+$ \\ \hline
$32\times32\times20\times10$& $4.54\times10^{-5}$& $3.06\times 10^{-5}$ & $51.79\%$ & $45$\\
$48\times48\times30\times15$& $4.76\times10^{-5}$&$3.28\times10^{-5}$&$49.88\%$&$72$\\
$64\times64\times40\times20$&$4.86\times10^{-5}$ &$3.44\times10^{-5}$&$25.28\%$ & $90$ \\
$80\times80\times50\times25$& $5.87\times10^{-5}$&$5.74\times10^{-5}$&$17.50\%$&$98$\\
$96\times96\times60\times30$& $4.85\times10^{-5}$&$3.21\times10^{-5}$&$13.09\%$&$106$\\
$112\times112\times70\times35$& $4.93\times10^{-5}$&$3.30\times10^{-5}$&$11.69\%$&$129$\\
$128\times128\times80\times40$& $5.82\times10^{-5}$&$5.37\times10^{-5}$&$10.04\%$&$145$\\
\hline
    \end{tabular}
      \caption{Matrix format. \label{tab:diffusion-matrix}}
\end{subtable}
\medskip
 \begin{subtable}{\textwidth}
   \centering
        \begin{tabular}
    {|l|c|c|c|c|c|c|c|c|c|c|c|}
    \hline
    $N_x\times N_y\times N_{\theta}\times N_{\BOmega_z}$ &  $\BmcD_{\bpsi^+}$ & $\BmcD_{\bphi}$ &Solution compression & Rank of $\psi^+$ \\ \hline
$32\times32\times20\times10$& $4.54\times10^{-5}$& $3.06\times 10^{-5}$ & $50.75\%$ & $[46,5,10]$\\
$48\times48\times30\times15$& $4.75\times10^{-5}$&$3.27\times10^{-5}$&$36.12\%$&$[75,5,15]$\\
$64\times64\times40\times20$&$4.86\times10^{-5}$ &$3.45\times10^{-5}$&$25.18\%$ & $[94,7,20]$ \\
$80\times80\times50\times25$& $5.84\times10^{-5}$&$5.69\times10^{-5}$&$17.65\%$&$[104,9,23]$\\
$96\times96\times60\times30$& $4.86\times10^{-5}$&$3.26\times10^{-5}$&$13.98\%$&$[118,14,28]$\\
$112\times112\times70\times35$& $4.99\times10^{-5}$&$3.44\times10^{-5}$&$14.59\%$&$[167,18,34]$\\
$128\times128\times80\times40$& $5.89\times10^{-5}$&$5.53\times10^{-5}$&$12.73\%$&$[192,20,39]$\\
\hline
    \end{tabular}
      \caption{Tensor format. \label{tab:diffusion-matrix}}
\end{subtable}
\caption{The difference between the low-rank solver  and the full-rank solver, the compression ratio and the rank for the variable scattering problem in Sec. \ref{sec:variable-scattering}.\label{tab:vs-results} }
\end{table}
\subsection{Pin-cell problem \label{sec:pin-cell}}
\begin{figure}[]
\centering
    \includegraphics[trim=10mm 9mm 10mm 0mm, clip,width=0.45\textwidth]{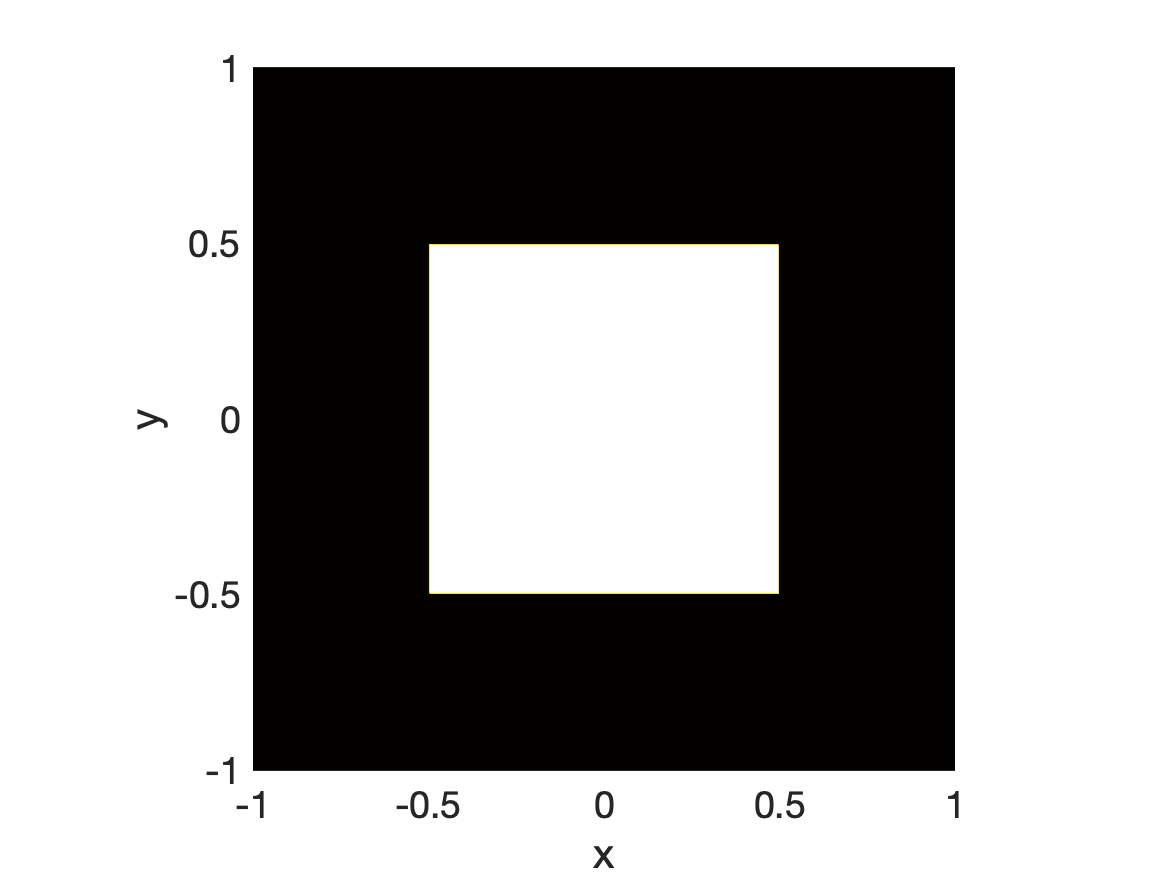}
    \hspace{-0.5cm}
    \includegraphics[width=0.45\textwidth]{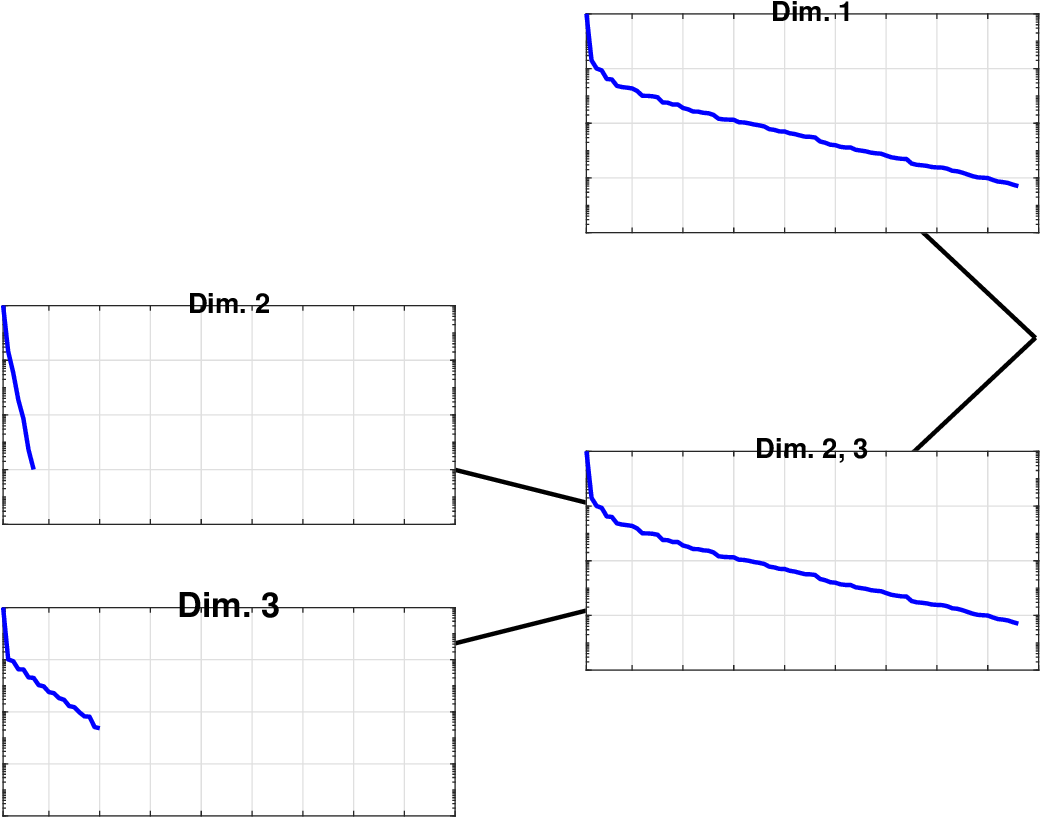}
   \begin{tikzpicture}
    \node[anchor=south west, inner sep=0] (image) at (0,0) {
     \includegraphics[trim=10mm 9mm 10mm 10mm, clip,width=0.45\textwidth]{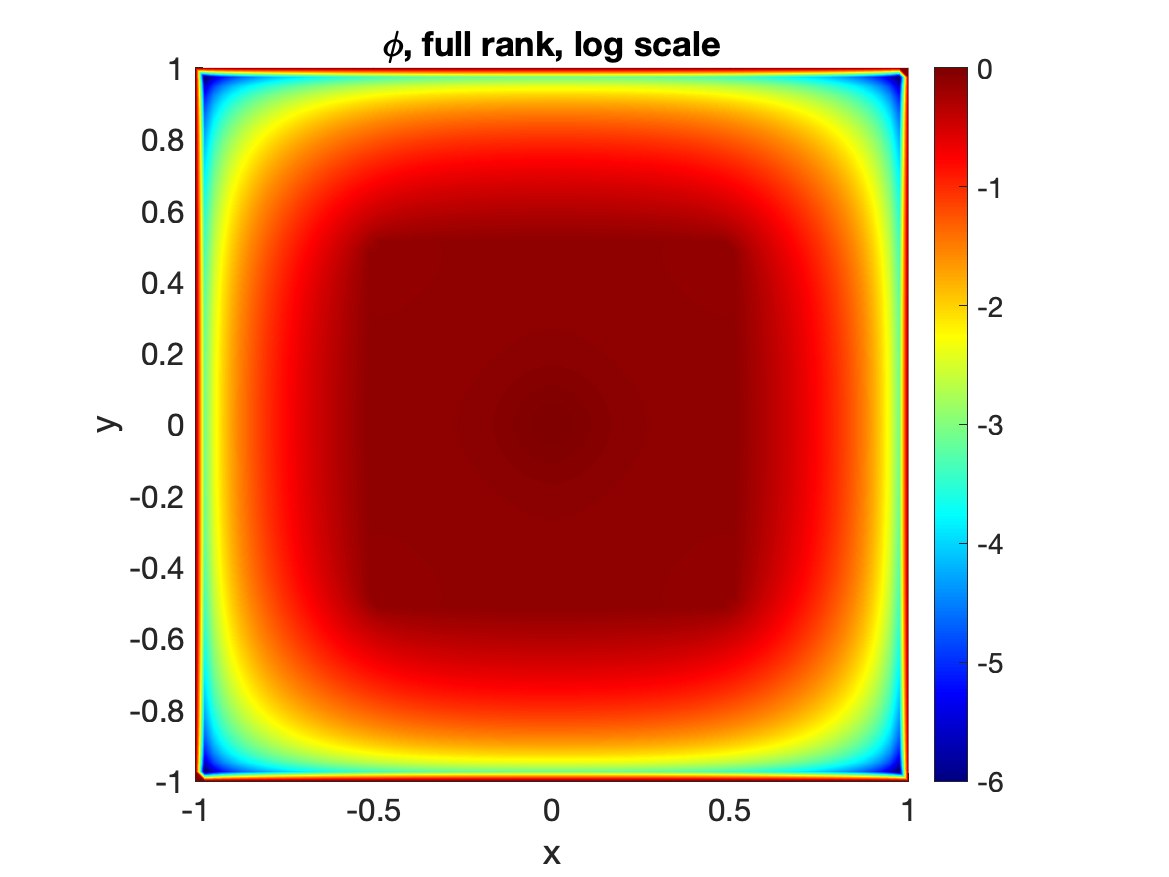}
    };
    \node[font=\bfseries, align=center, above=-1.2mm of image.north] {
      \tiny $\phi$, full-rank, log-scale
    };
  \end{tikzpicture}
  \hspace{-2mm}
  \begin{tikzpicture}
    \node[anchor=south west, inner sep=0] (image) at (0,0) {
     \includegraphics[trim=10mm 9mm 10mm 10mm, clip,width=0.45\textwidth]{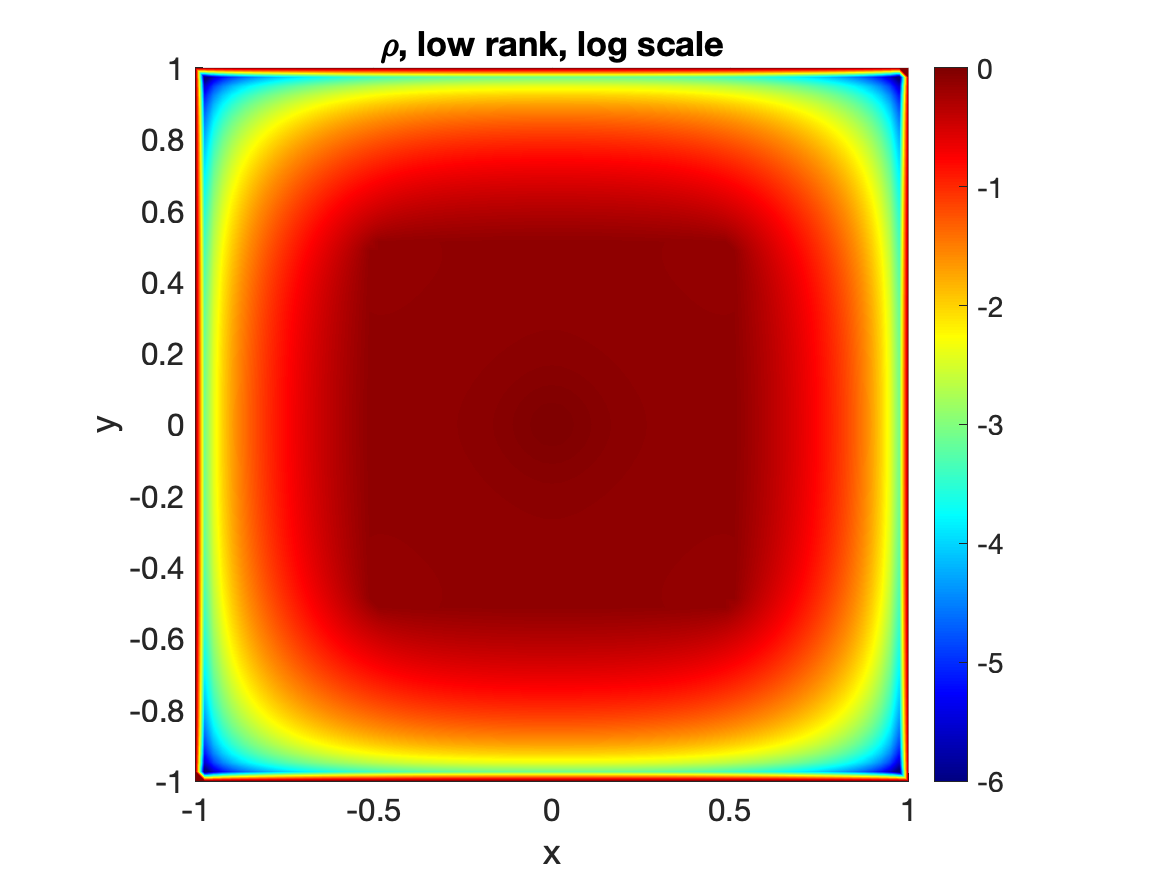}
    };
    \node[font=\bfseries, align=center, above=-1.2mm of image.north] {
      \tiny $\phi$, low-rank, tensor, log-scale
    };
  \end{tikzpicture}
  \caption{Pin-cell problem in Sec. \ref{sec:pin-cell}, $(N_x,N_y,N_\theta,N_{\BOmega_z})=(78,78,40,20)$. Top left: set-up for the scattering cross section. Top right: singular values for the matricization of the low-rank tensor. Bottom left: full-rank solution. Bottom right: low-rank solution in the order 3 tensor format.  \label{fig:pin-cell-solution}}
\end{figure}

We consider a pin-cell problem on the computational domain $[-1,1]^2$ with vacuum boundary conditions, no absorption and  a Gaussian source $G(x,y)=\exp(-100(x^2+y^2))$. The scattering cross section is 
\begin{equation}
    \sigma_s(x,y) = \begin{cases}
                    1,\quad\text{if}\;|x|\leq0.5\;\text{and}\;|y|\leq0.5,\\
                    100,\quad\text{otherwise}.
                    \end{cases}
\end{equation}
The visualization for the set-up of the scattering cross section is presented in the top left plot of Fig. \ref{fig:pin-cell-solution}. The outer black region is diffusion dominant, while the inner white region is transport dominant. We set $(N_x,N_y,N_\theta,N_{\BOmega_z})=(78,78,40,20)$. 

\begin{figure}[]
  \begin{center} 
\includegraphics[width=0.45\textwidth]{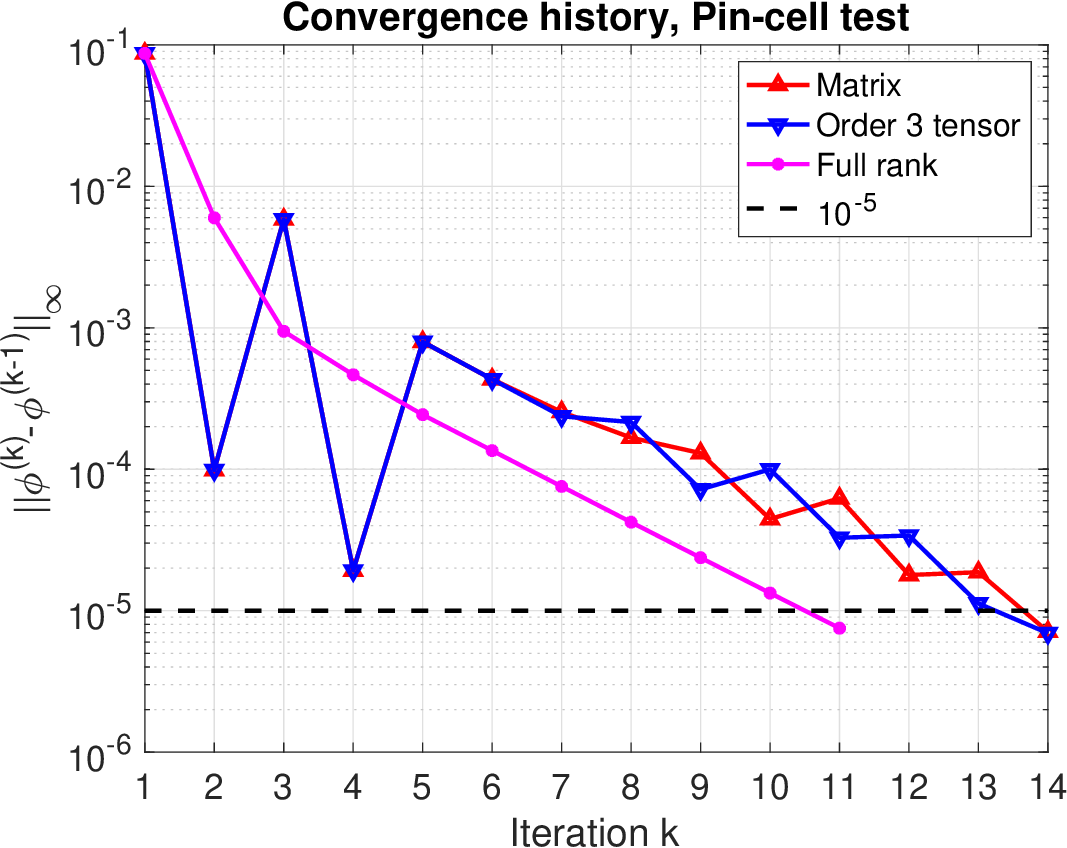}
\includegraphics[width=0.45\textwidth]{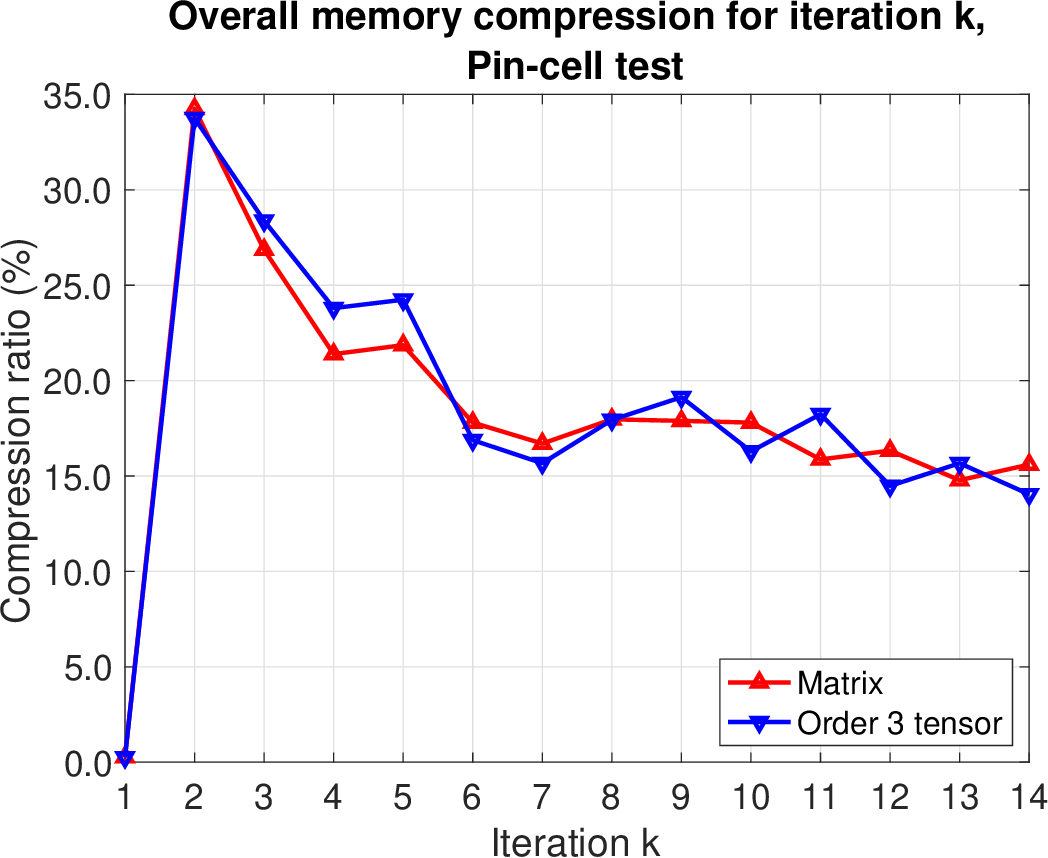}
  \caption{The history of $||\phi^{(k)}-\phi^{(k-1)}||$ and the memory compression ratio for the pin-cell problem in Sec. \ref{sec:pin-cell} with $(N_x,N_y,N_\theta,N_{\BOmega_z})=(78,78,40,20)$. \label{fig:pin-cell-convergence-history}}
  \end{center}
\end{figure}

As shown in Tab. \ref{tab:pin-cell-results}, the difference between the low-rank and full-rank results is approximately $5\times10^{-5}$ in the $l_2$ norm, while using only about $21.89\%$ of the DOFs in the matrix format and $22.34\%$ in the tensor format relative to the full-rank solver.
As presented in Fig. \ref{fig:pin-cell-convergence-history}, the low-rank solvers require slightly more iterations. The DOFs still increase sharply after the first iteration and then exhibit a gradual decreasing trend. 
\begin{table}[htbp]
  \centering
 \medskip
 \begin{subtable}{\textwidth}
   \centering
    \begin{tabular}
    {|l|c|c|c|c|c|c|c|c|c|c|c|}
    \hline
    $N_x\times N_y\times N_{\theta}\times N_{\BOmega_z}$ &  $\BmcD_{\bpsi^+}$ & $\BmcD_{\bphi}$ &Solution compression & Rank of $\psi^+$ \\ \hline
$78\times78\times40\times20$& $4.69\times10^{-5}$& $6.55\times 10^{-5}$ & $21.89\%$ & $81$\\
\hline
    \end{tabular}
      \caption{Matrix format. \label{tab:pin-cell-matrix}}
\end{subtable}
\medskip
 \begin{subtable}{\textwidth}
   \centering
        \begin{tabular}
    {|l|c|c|c|c|c|c|c|c|c|c|c|}
    \hline
    $N_x\times N_y\times N_{\theta}\times N_{\BOmega_z}$ &  $\BmcD_{\bpsi^+}$ & $\BmcD_{\bphi}$ &Solution compression & Rank of $\psi^+$ \\ \hline
$78\times78\times40\times20$& $4.97\times10^{-5}$& $6.89\times 10^{-5}$ & $22.34\%$ & $[86,7,20]$\\
\hline
    \end{tabular}
      \caption{HTT format. \label{tab:pin-cell-tensor}}
\end{subtable}
\caption{The difference between the low-rank solver  and the full-rank solver, the compression ratio and the rank for the pin-cell problem in Sec. \ref{sec:pin-cell}.\label{tab:pin-cell-results} }
\end{table}
\subsection{Lattice problem\label{sec:lattice}}
\begin{figure}[]
\centering
    \includegraphics[trim=10mm 9mm 10mm 0mm, clip,width=0.45\textwidth]{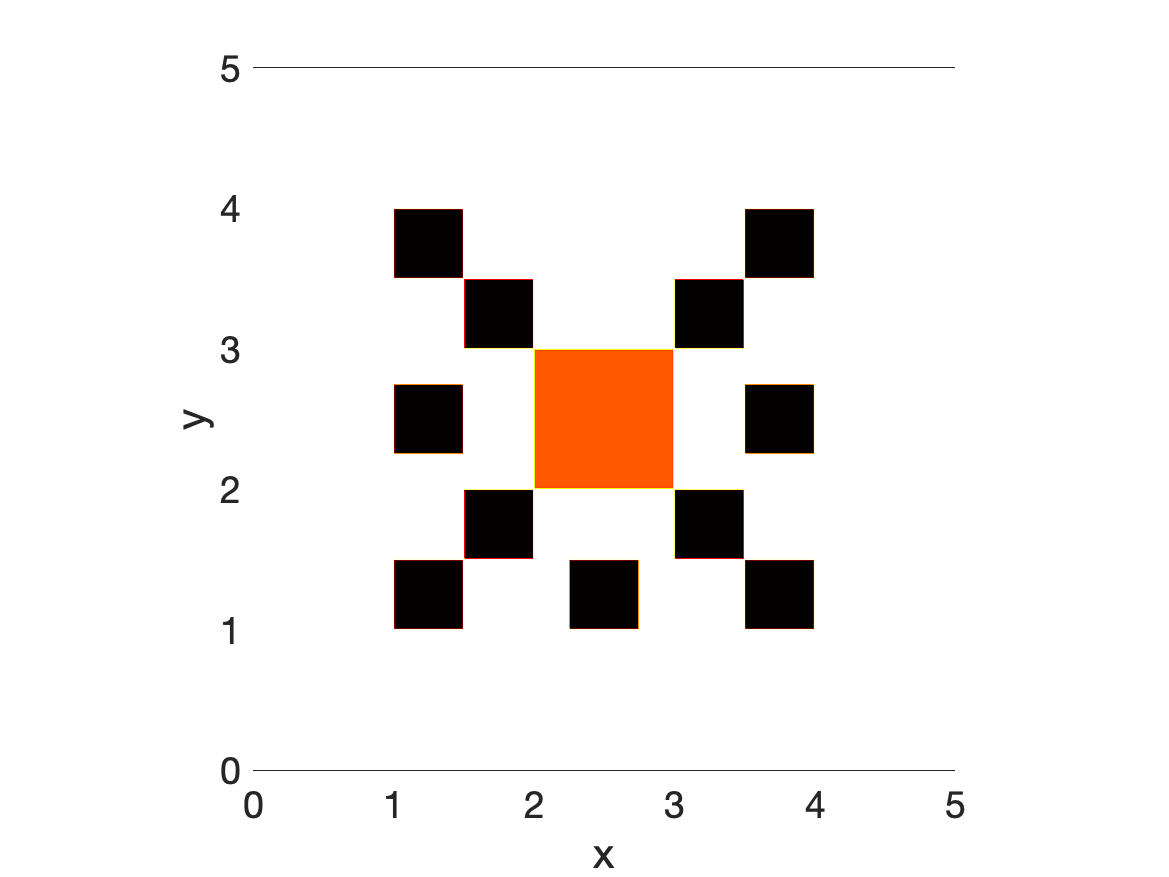}
    \hspace{-0.5cm}
    \includegraphics[width=0.45\textwidth]{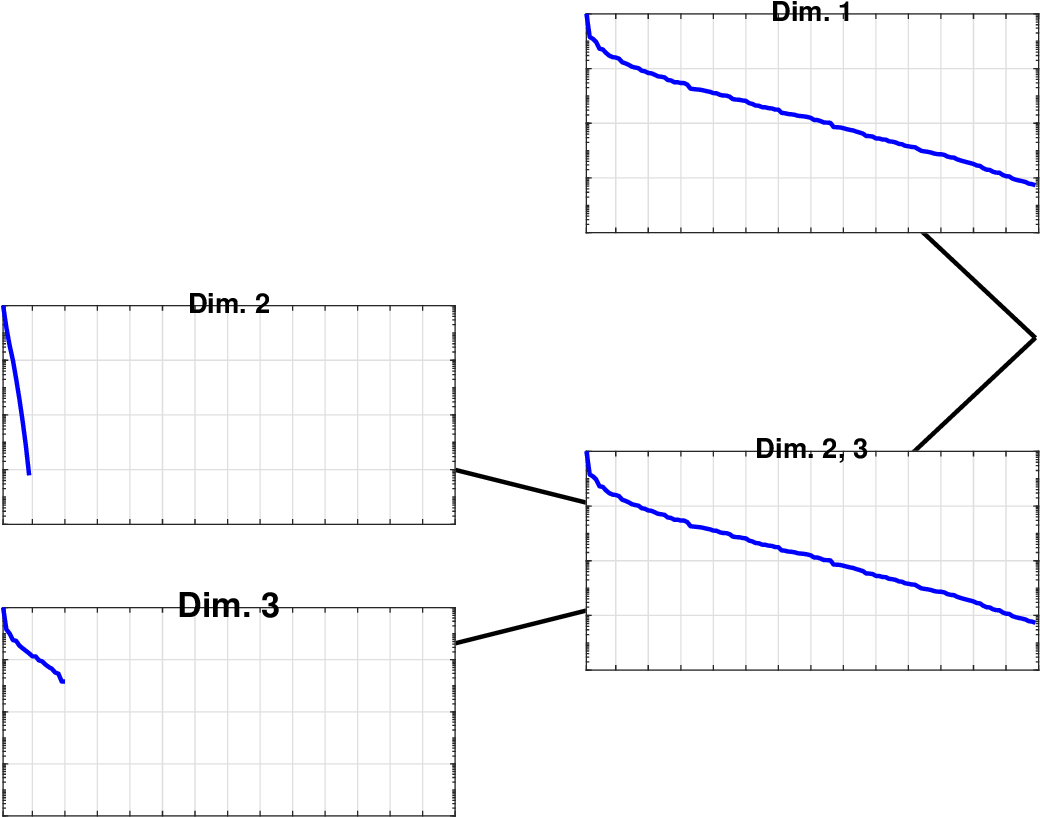}
   \begin{tikzpicture}
    \node[anchor=south west, inner sep=0] (image) at (0,0) {
     \includegraphics[trim=10mm 9mm 10mm 10mm, clip,width=0.45\textwidth]{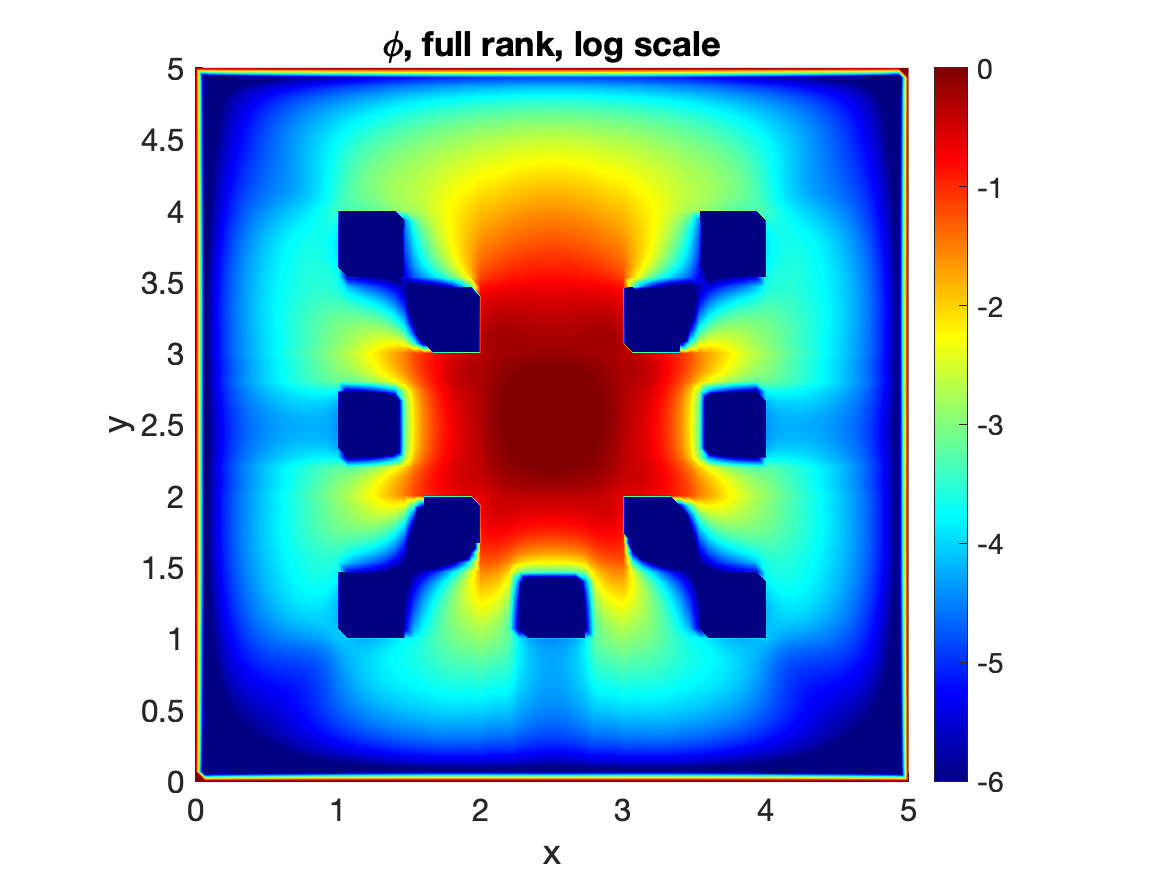}
    };
    \node[font=\bfseries, align=center, above=-1.2mm of image.north] {
      \tiny $\phi$, full-rank, log-scale
    };
  \end{tikzpicture}
  \hspace{-2mm}
  \begin{tikzpicture}
    \node[anchor=south west, inner sep=0] (image) at (0,0) {
     \includegraphics[trim=10mm 9mm 10mm 10mm, clip,width=0.45\textwidth]{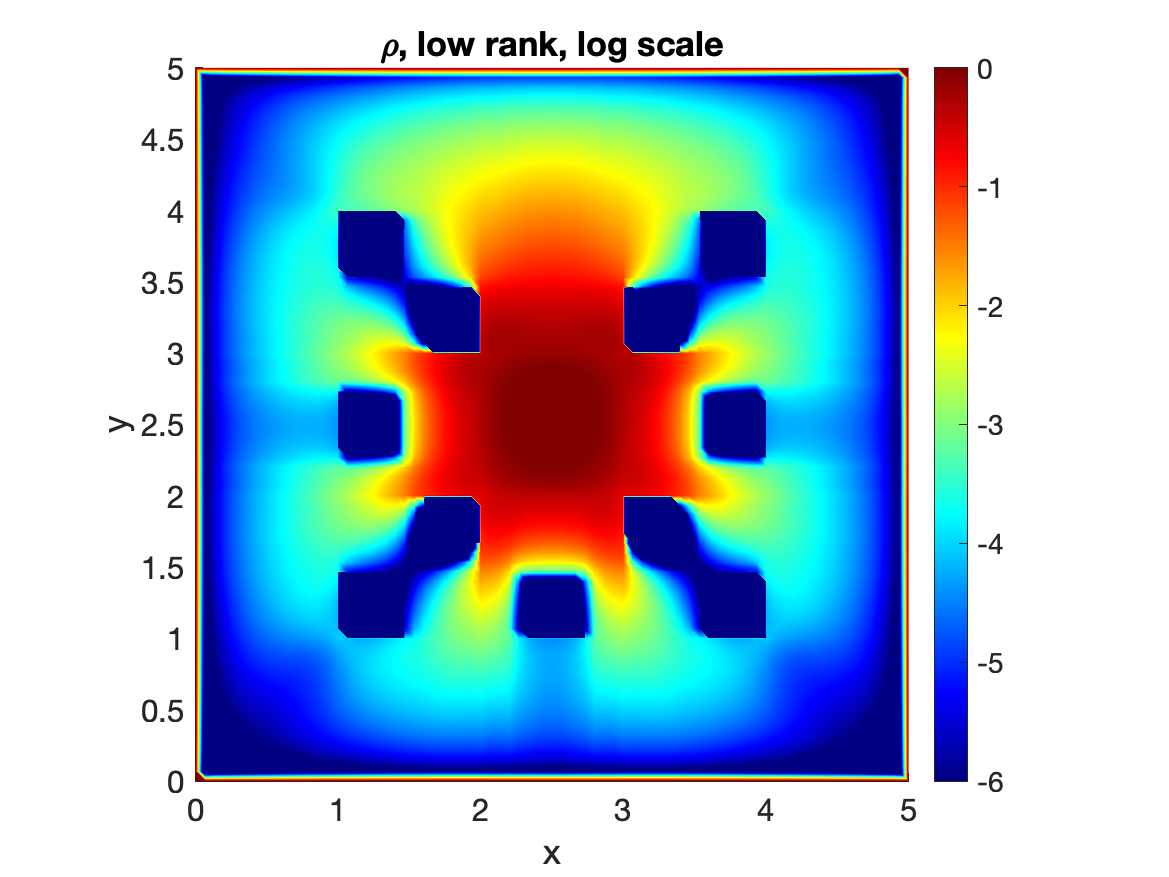}
    };
    \node[font=\bfseries, align=center, above=-1.2mm of image.north] {
      \tiny $\phi$, low-rank, tensor, log-scale
    };
  \end{tikzpicture}
  \caption{Lattice problem in Sec. \ref{sec:lattice}, $(N_x,N_y,N_\theta,N_{\BOmega_z})=(75,75,40,20)$. Top left: set-up for the material property. Top right: singular values for the matricization of the low-rank tensor. Bottom left: full-rank solution. Bottom right: low-rank solution in the order 3 tensor format.  \label{fig:lattice-solution}}
\end{figure}

In the last example, we consider a lattice problem with vaccum boundary conditions on the computational domain $[0,5]^2$. The material properties are shown in the top-left plot of Fig. \ref{fig:lattice-solution}. The black region represents pure absorption, with $\sigma_s(x,y)=0$ and $\sigma_a(x,y)=100$. The white and orange regions represent pure scattering, with $\sigma_s(x,y)=1$ and $\sigma_a(x,y)=0$. A constant source is imposed in the orange region.

\begin{figure}[]
  \begin{center} 
\includegraphics[width=0.45\textwidth]{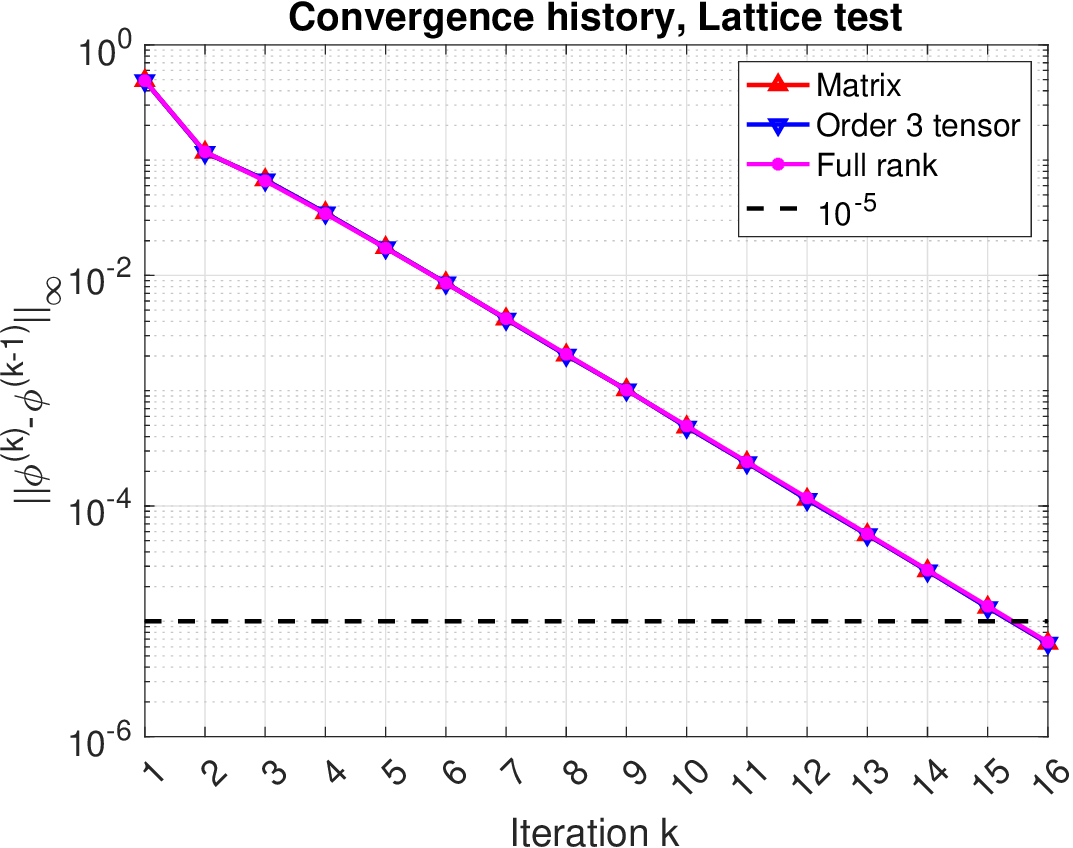}
\includegraphics[width=0.45\textwidth]{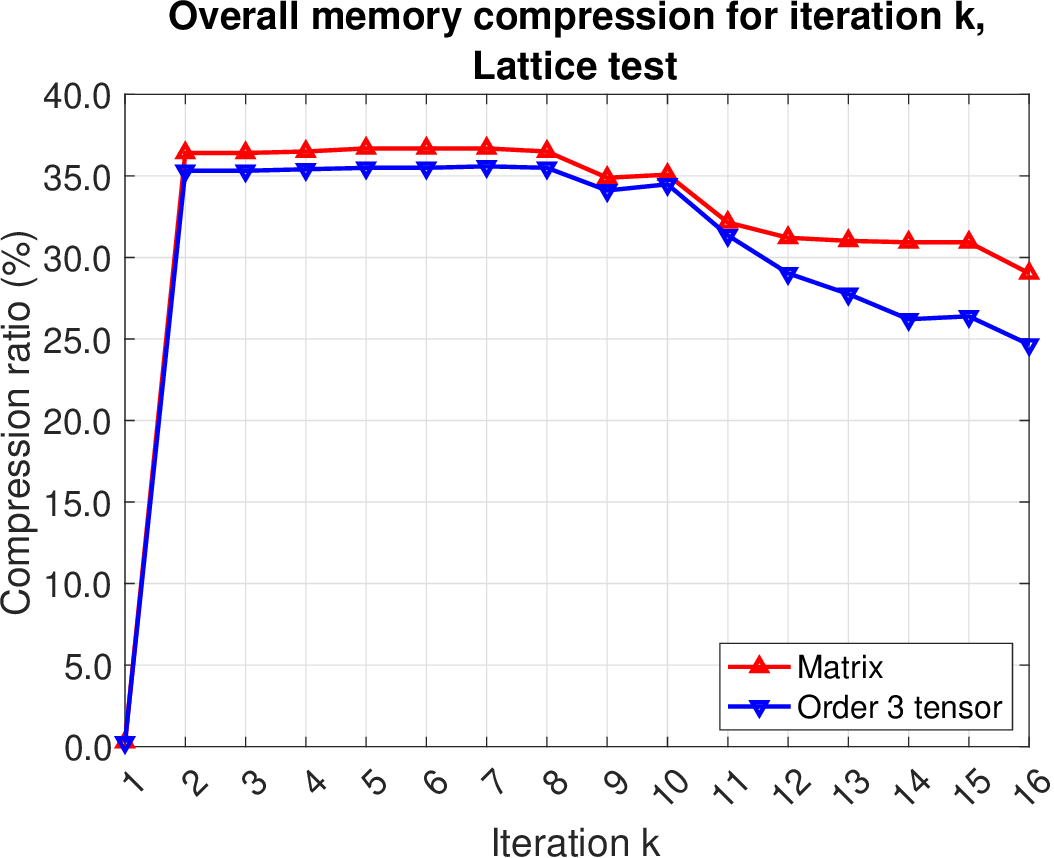}
  \caption{The history of $||\phi^{(k)}-\phi^{(k-1)}||$ and the memory compression ratio for the lattice problem in Sec. \ref{sec:lattice} with $(N_x,N_y,N_\theta,N_{\BOmega_z})=(75,75,40,20)$. \label{fig:lattice-convergence-history}}
  \end{center}
\end{figure}

As shown in Tab. \ref{tab:lattice-results}, the $l_2$ difference between the low-rank and full-rank results is smaller than $2\times 10^{-6}$, while the low-rank method reduces the DOFs by a factor of approximately $2.70$ in the matrix format and $2.72$ in the tensor format. As presented in Fig. \ref{fig:lattice-convergence-history}, the convergence histories of the full-rank method and the two low-rank methods qualitatively match. In the low-rank methods, the DOFs increase sharply after the first iteration, then stabilize at a plateau, and gradually decrease in the final iterations.
\begin{table}[htbp]
  \centering
 \medskip
 \begin{subtable}{\textwidth}
   \centering
    \begin{tabular}
    {|l|c|c|c|c|c|c|c|c|c|c|c|}
    \hline
    $N_x\times N_y\times N_{\theta}\times N_{\BOmega_z}$ &  $\BmcD_{\bpsi^+}$ & $\BmcD_{\bphi}$ &Solution compression & Rank of $\psi^+$ \\ \hline
$75\times75\times40\times20$& $1.74\times10^{-6}$& $4.12\times 10^{-7}$ & $37.04\%$ & $135$\\
\hline
    \end{tabular}
      \caption{Matrix format. \label{tab:lattice-matrix}}
\end{subtable}
\medskip
 \begin{subtable}{\textwidth}
   \centering
        \begin{tabular}
    {|l|c|c|c|c|c|c|c|c|c|c|c|}
    \hline
    $N_x\times N_y\times N_{\theta}\times N_{\BOmega_z}$ &  $\BmcD_{\bpsi^+}$ & $\BmcD_{\bphi}$ &Solution compression & Rank of $\psi^+$ \\ \hline
$75\times75\times40\times20$& $1.27\times10^{-6}$& $5.15\times 10^{-7}$ & $36.80\%$ & $[139,9,20]$\\
\hline
    \end{tabular}
      \caption{Tensor format. \label{tab:lattice-tensor}}
\end{subtable}
\caption{The difference between the low-rank solver  and the full-rank solver, the compression ratio and the rank for the lattice problem in Sec. \ref{sec:lattice}.\label{tab:lattice-results} }
\end{table}
\section{Conclusions\label{sec:conclusion}}
In this paper, to reduce memory cost, we develop an inexact low-rank SI-DSA solver for multi-D mulit-V RTE based on its second-order formulation. To the best of our knowledge, apart from the recent work on 1D1V problems \cite{bachmayr2024low}, this paper is the first to design preconditioners and low-rank iterative solvers for steady-state RTEs.  

We numerically demonstrate that, even for challenging problems where the  effective rank of the solution grows as resolution increases, our low-rank method can still effectively compress the solution and reduce the required DOFs in the numerical simulation. The performance of our method using the matrix and HTT format are comparable. Hence, we recommend the matrix format due to its simpler implementation and smaller computational overhead. The development of more efficient HTT decompositions for RTE simulations to further reduce overhead is left for future work.
 
When the problem admits a low-rank structure, the proposed low-rank method can achieve significant acceleration over full-rank methods. However, for more challenging problems where the rank of the low-rank solution grows with mesh resolution, performance strongly depends on the efficiency of the inner-loop preconditioner. With the same preconditioner used, the low-rank method is consistently faster than the full-rank method. Yet, because the proposed two-level DPC is less efficient than standalone AMG, the low-rank method with DPC may become slower than the full-rank method with AMG for these challenging problems.  We expect that once a low-rank preconditioner with efficiency comparable to standalone AMG is developed, the low-rank approach is capable of delivering a clear speedup across all scenarios. 

The most urgent next step is to design a more efficient low-rank inner-loop preconditioner. A promising direction is to extend the space–angle multigrid method to the low-rank setting. In addition, since the second-order formulation requires $\sigma_t>0$, a low-rank solver based on the first-order formulation is necessary for problems where $\sigma_t=0$ in certain subregions. Finally, extensions to multigroup problems, anisotropic scattering, and nonlinear thermal radiation remain highly desirable.

\subsection*{Acknowledgments}
W. Guo was partially supported by the Air Force Office of Scientific Research FA9550-18-1-0257. Z. Peng was partially supported by the Hong Kong Research Grants Council grants Early Career Scheme 26302724 and General Research Fund 16306825. 
\appendix
\section{Performance of the inner loop diffusion and AMG preconditioner\label{apx:inner-pc}}
Here, at the full-rank level, we compare the performance of the preconditioners for the inner-loop CG solver in each source iteration, namely the standalone AMG and the two-level DPC discussed in Sec.~\ref{sec:lr-inner}. We use CL($32,16$) quadrature rule in the angular space and perform a spatial refinement study for the diffusion dominant problem in Sec. \ref{sec:diffusion}, the transport dominant problem in Sec. \ref{sec:transport} and the variable scattering problem in Sec. \ref{sec:variable-scattering}.

To facilitate the performance comparison of the two preconditioners, we introduce the following shorthand notations.
\begin{itemize}
    \item $N_{\textrm{SI}}$: Total number of outer source iterations
    \item $\bar{N}_{\textrm{Inner-CG}}$: Average  iterations of inner loop CG solver per angular direction per outer source iteration
    \item $\bar{N}_{\textrm{IPC}}$: Average iterations at the second level of DPC to solve the system determine by the diffusion operator per angular direction per inner iteration
\end{itemize}
We set the stopping criteria of the outer loop source iteration as $||\bphi^{(k)}-\bphi^{(k-1)}||\leq 10^{-5}$ and the inner loop solver as the relative residual smaller than $10^{-9}$.

The comparison results for the AMG inner loop preconditioner and DPC are summarized in Tab. \ref{tab:pc-compare}. Main observations are as follows.
\begin{enumerate}
    \item For the diffusion dominant problem, the number of inner loop iterations for convergence required by both AMG and DPC almost does not grow as spatial resolution increases. 
    \item For transport dominant problem and multiscale variable scattering problem, inner loop iterations grows as spatial resolution increases, regardless of using DPC and AMG.
    \item Using two-level DPC, each SI step requires comparable or even less inner loop iterations to solve \eqref{eq:SI} compared to the AMG inner loop preconditioner. However, DPC is less efficient due to additional iterations at its second level to  numerically solve the linear system determined by its diffusion operator. Moreover, for transport dominant and variable scattering problem, the iteration number of the AMG solver inside the second level of DPC has a sub-linear growth as the spatial resolution improves. 
\end{enumerate}

\begin{table}[htbp]
  \centering
 \medskip
 \begin{subtable}{\textwidth}
   \centering
    \begin{tabular}
    {|l|c|c|c|c|c|c|c|c|c|c|c|}
    \hline
 & \multicolumn{2}{|c|}{AMG} & \multicolumn{3}{|c|}{DPC} \\ \hline & & & & & \\[-0.9em]    
 $N_x\times N_y$   & $N_{\textrm{SI}}$ & $\bar{N}_{\textrm{Inner-CG}}$ & $N_{\textrm{SI}}$ & $\bar{N}_{\textrm{Inner-CG}}$ & $\bar{N}_{\textrm{IPC}}$ \\ \hline
 $32\times 32$ & $3$ & $2.41$ & $3$ & $3.70$ & $3.00$\\ \hline
 $48\times 48$ & $3$ & $2.63$ & $3$ & $3.87$ & $3.00$ \\ \hline
 $64\times 64$ & $3$ & $2.75$ & $3$ & $4.00$ & $3.00$ \\ \hline
 $80\times 80$ & $3$ & $2.76$ & $3$ & $4.14$ & $3.00$ \\ \hline
    \end{tabular}
      \caption{Diffusion dominant problem.}
\end{subtable}
 \medskip
 \begin{subtable}{\textwidth}
   \centering
    \begin{tabular}
    {|l|c|c|c|c|c|c|c|c|c|c|c|}
    \hline
 & \multicolumn{2}{|c|}{AMG} & \multicolumn{3}{|c|}{DPC} \\ \hline & & & & & \\[-0.9em]    
 $N_x\times N_y$   & $N_{\textrm{SI}}$ & $\bar{N}_{\textrm{Inner-CG}}$ & $N_{\textrm{SI}}$ & $\bar{N}_{\textrm{Inner-CG}}$ & $\bar{N}_{\textrm{IPC}}$ \\ \hline
 $32\times 32$ & $6$ & $30.41$ & $6$ & $35.48$ & $26.08$\\ \hline
 $48\times 48$ & $6$ & $45.60$ & $6$ & $43.69$ & $36.70$ \\ \hline
 $64\times 64$ & $6$ & $60.96$ & $6$ & $50.28$ & $47.34$ \\ \hline
 $80\times 80$ & $6$ & $76.20$ & $6$ & $55.91$ & $57.86$ \\ \hline
    \end{tabular}
      \caption{Transport dominant problem.}
\end{subtable}
 \medskip
 \begin{subtable}{\textwidth}
   \centering
    \begin{tabular}
    {|l|c|c|c|c|c|c|c|c|c|c|c|}
    \hline
 & \multicolumn{2}{|c|}{AMG} & \multicolumn{3}{|c|}{DPC} \\ \hline & & & & & \\[-0.9em]    
 $N_x\times N_y$   & $N_{\textrm{SI}}$ & $\bar{N}_{\textrm{Inner-CG}}$ & $N_{\textrm{SI}}$ & $\bar{N}_{\textrm{Inner-CG}}$ & $\bar{N}_{\textrm{IPC}}$ \\ \hline
 $32\times 32$ & $5$ & $14.92$ & $6$ & $21.12$ & $13.41$\\ \hline
 $48\times 48$ & $5$ & $22.59$ & $5$ & $30.41$ & $18.13$ \\ \hline
 $64\times 64$ & $5$ & $30.45$ & $5$ & $35.11$ & $22.37$ \\ \hline
 $80\times 80$ & $5$ & $38.06$ & $3$ & $38.75$ & $26.47$ \\ \hline
    \end{tabular}
      \caption{Variable scattering problem.}
\end{subtable}
\caption{Comparison between the full-rank AMG and DPC.\label{tab:pc-compare}}
\end{table}
\bibliographystyle{elsarticle-num} 
\bibliography{ref}

\end{document}